\documentclass{amsart}
\usepackage{epsf}
\usepackage{epic}
\usepackage{eepic}
\usepackage{amssymb}
\usepackage{amsfonts}
\usepackage{amscd}
\usepackage{verbatim}
\usepackage{epsfig}
\input xy
\xyoption{all}

\newtheorem{guess}{Theorem}[section]
\newcommand{\bth}{\begin{guess}}
\newcommand{\eeth}{\end{guess}}

\newtheorem{propo}[guess]{Proposition}
\newcommand{\bpropo}{\begin{propo}}
\newcommand{\epropo}{\end{propo}}

\newtheorem{lema}[guess]{Lemma}
\newcommand{\blem}{\begin{lema}}
\newcommand{\elem}{\end{lema}}

\newtheorem{defe}[guess]{Definition}
\newcommand{\bdefe}{\begin{defe}}
\newcommand{\edefe}{\end{defe}}

\newtheorem{coro}[guess]{Corollary}
\newcommand{\bcor}{\begin{coro}}
\newcommand{\ecor}{\end{coro}}

\newtheorem{rema}[guess]{Remark}
\newcommand{\brem}{\begin{rema}$\!\!\!$~\rm}
\newcommand{\erem}{\end{rema}}

\DeclareFontFamily{OT1}{rsfs}{}
\DeclareFontShape{OT1}{rsfs}{n}{it}{<-> rsfs10}{}
\DeclareMathAlphabet{\mathscr}{OT1}{rsfs}{n}{it}

\newtheorem{thm}{Theorem}[subsection]
\newtheorem{lem}[thm]{Lemma}
\newtheorem{prop}[thm]{Proposition}
\newtheorem{cor}[thm]{Corollary}
\newtheorem{rem}[thm]{Remark}
\newtheorem{rems}[thm]{Remarks}

\theoremstyle{definition}

  \newtheorem{defi}[thm]{Definition}
  \newtheorem{defis}[thm]{Definitions}

  \newtheorem{conv}[thm]{Convention}
  \newtheorem{convs}[thm]{Conventions}

\numberwithin{equation}{thm}

\newcommand{\Cref}[1]{Corollary~\textup{\ref{#1}}}
\newcommand{\Dref}[1]{Definition~\textup{\ref{#1}}}

\newcommand{\Lref}[1]{Lemma~\textup{\ref{#1}}}
\newcommand{\Pref}[1]{Proposition~\textup{\ref{#1}}}
\newcommand{\Rref}[1]{Remark~\textup{\ref{#1}}}

\newcommand{\Sref}[1]{Section~\textup{\ref{#1}}}
\newcommand{\Ssref}[1]{Subsection~\textup{\ref{#1}}}
\newcommand{\Tref}[1]{Theorem~\textup{\ref{#1}}}
\newcommand{\smcirc}%
   {{\raise.15ex\hbox to.7em{$\hss \scriptstyle\circ\hss$}}}

\def\bilap#1{\hbox to 0pt{\hss#1\hss}}
  \def\Rarrow#1{\bilap{\hbox to#1{\rightarrowfill}}}
  \def\Larrow#1{\bilap{\hbox to#1{\leftarrowfill}}}
  \def\Equals#1{\bilap
                   {\raise 4pt\hbox
                     {\vrule width#1 height.5pt}%
                    \kern-#1\raise 1pt\hbox
                     {\vrule width#1 height.5pt}%
                   }}

\newcommand{\EQAL}[1]%
{\,\begin{picture}(#1,0)%
\put(0,3){\line(1,0){#1}}%
\put(0,1){\line(1,0){#1}}%
\end{picture}\,}%

\newcommand{\vlto}[1]%
{\,\begin{picture}(#1,3)%
\put(0,2){\vector(1,0){#1}}%
\end{picture}\,}%

\newcommand{\vllarrow}[1]%
{\,\begin{picture}(#1,3)%
\put(#1,2){\vector(-1,0){#1}}%
\end{picture}\,}%

\newcommand{\dirlm}[1]%
   {
      {\lim\hskip-1.58em\lower.65ex
        \hbox{$
                 {}_{\stackrel{\lower1ex\hbox
                                         {$\scriptstyle -\!\!\!\longrightarrow$}
                                       }{\vbox to0pt{\vss\vskip.6ex
                                             \hbox{$\scriptstyle{}^{#1}$}\vss}}
                    }
             $}
      }
\:}

\newcommand{\subdirlm}[1]%
   {
      {\lim\hskip-1.5em\lower.6ex
        \hbox{$
                    {}_{\stackrel{\lower1ex\hbox
                                            {$\scriptstyle\longrightarrow$}
                                 }{ ^{#1} }
                       }
              $}
      }
\:}

\newcommand{\inlm}[1]%
    {
       {\lim\hskip-1.58em\lower.65ex
         \hbox{$
                  {}_{\stackrel{\lower1ex\hbox
                                         {$\scriptstyle \longleftarrow\!\!\<-$}
                               }{\vbox to0pt{\vss\vskip.6ex
                                             \hbox{$\scriptstyle{}^{#1}$}\vss}}
                     }
              $}
       }
\:}

\def\hz#1{{\hbox to 0pt{#1}}}

\def\>{\mspace {1mu}}
\def\<{\mspace{-1mu}}
\def\({{\textup(}}
\def\){{\textup)}}
\def\bigl#1{{\textup{\begin{large}#1\end{large}}}}
\def\bigr#1{{\textup{\begin{large}#1\end{large}}}}

\newcommand{\Lr}{\Longrightarrow}
\newcommand{\Llr}{\Longleftrightarrow}
\newcommand{\es}{\emptyset}
\renewcommand{\mod}{{\rm mod}\,}

\newcommand{\set}{\!:=}

\newcommand{\ra}{\rightarrow}

\newcommand{\lr}{\longrightarrow}
\newcommand{\cu}{{\mathcal U}}
\newcommand{\co}{{\mathcal O}}

\newcommand{\X}{{\mathscr X}}

\newcommand{\Y}{{\mathscr Y}}
\newcommand{\Z}{{\mathscr Z}}

\newcommand{\U}{{\mathscr U}}
\newcommand{\W}{{\mathscr W}}

\newcommand{\eL}{{\mathscr L}}
\newcommand{\ba}{{\mathbb A}}
\newcommand{\bg}{{\mathbb G}}
\newcommand{\bp}{{\mathbb P}}
\newcommand{\De}{\Delta}
\newcommand{\wh}{\widehat}
\newcommand{\wt}{\widetilde}
\newcommand{\ds}{\displaystyle}
\newcommand{\ol}{\overline}
\newcommand{\spec}{{\mathrm{Spec}\,}}
\newcommand{\proj}{{\mathrm{Proj}}\,}

\newcommand{\br}{{\mathbb R}}
\newcommand{\pic}{{\text{Pic}}\,}
\newcommand{\supp}{{\text{supp}}\,}
\newcommand{\hra}{\hookrightarrow}

\newcommand{\iso}%
{{\mkern8mu\longrightarrow \mkern-25.5mu{}^\sim\mkern17mu}}

\newcounter{savectr}
\title[Geometric Reductivity]{Geometric reductivity\\-- A quotient space approach}
\author[P.\,Sastry]{Pramathanath Sastry}%
\address{Chennai Mathematical Institute\\
Plot H4, SIPCOT IT Park\\
Padur Post, Siruseri\\
Kanchipuram District 603 103\\
India. }
\email{pramath@cmi.ac.in}
\author{C.S. Seshadri}
\address{Chennai Mathematical Institute\\
Plot H4, SIPCOT IT Park\\
Padur Post, Siruseri\\
Kanchipuram District 603 103\\
India. }
\email{css@cmi.ac.in}

\date{\today}

\begin{document}

\maketitle

\tableofcontents
\section{{\bf Introduction}}
Mumford's {\em Geometric Invariant Theory} or {\em GIT} is a major technique for
finding quotients of algebraic schemes acted upon by reductive algebraic groups.
It has been successful in finding solutions to moduli problems in the category of {\em algebraic schemes}. 
In the first edition (i.e., the 1965 edition) of  {\em Geometric Invariant Theory} \cite{M},
Mumford restricted himself to algebraic schemes over fields of characteristic zero.
In order to make his theory applicable over fields of arbitrary characteristic, he
made the following conjecture in the Preface to the first edition of {\em ibid}.\,(a conjecture 
subsequently proved by Haboush \cite{H} in 1975):
 
{\em Let $G$ be a reductive algebraic group over an algebraically closed field $k$. Then 
$G$ is {\sc geometrically reductive}, i.e., for
every finite-dimensional rational $G$-module $V$ and a $G$-invariant point $v\in V$, $v\neq 0$,
there is a $G$-invariant homogeneous polynomial $F$ on $V$ of positive degree such
that $F(v)\neq 0$.}
 
As a consequence, it can be shown that
if $X=\spec{A}$ is an algebraic scheme on which  a our reductive algebraic group
$G$ acts, then the affine scheme $Y=\spec{A^G}$ is an algebraic scheme, i.e,
the ring of invariants $A^G$ is finitely generated as a $k$-algebra (a result of Nagata
\cite{N}) and the canonical
morphism $f\colon X\to Y$ (induced by the injection $A^G\hra A$) is surjective.
Further, if $Z$ is a closed $G$-stable
subset of $X$, then $f(Z)$ is closed in $Y$, and if $f$ separates {\em disjoint closed $G$-stable}
subsets of $X$, i.e., given two disjoint closed $G$-stable subsets $Z_1$ and $Z_2$ of $X$,
then $f(Z_1)$ and $f(Z_2)$ are also disjoint. In other words, $f\colon X\to Y$ is what is called
a {\em good quotient} \cite{S2}. These results are in fact {\em equivalent} to the
conjecture.
 
A major consequence then is that Mumford's technique---as set out in \cite{M}---for 
constructing the quotient (in the category of algebraic schemes)
of the {\em semi-stable locus} of a projective algebraic scheme $X$ on which a reductive
group $G$ acts linearly\footnote{i.e., $G$ acts linearly on the ambient projective space in which
$X$ is embedded.}
works over fields of arbitrary characteristic. Recall that in the first 
edition of {\em ibid}.\,such quotients were only constructed when the underlying field was of
characteristic zero, and over such a field reductive algebraic groups have been known,
via Hermann Weyl's work (see \cite{W}), to be
linearly reductive, whence geometrically reductive. 
In fact in  characteristic zero geometric reductivity is 
equivalent to the complete reducibility of finite dimensional
$G$-modules.

Geometric reductivity (for our reductive algebraic group $G$) 
was first proved for the case of 
$SL(2)$ (hence $GL(2)$) in characteristic 2 by Oda \cite{O},  and in general
by the second author \cite{S1}.  
Haboush's proof \cite{H} uses in an essential way the irreducibity of the Steinberg representation.
There is also a different approach
to the problem due to Formanek and Procesi, \`a priori 
for the full linear group \cite{F-P}, but the general case can be
deduced from this (see \cite{S5}). 

Successful as the above approach has been in solving the problem of quotients (of algebraic
schemes by reductive agebraic groups), a direct attack on the quotient problem has an undoubted philosophical attraction. 
And success here would yield all the consequences of geometric reductivity
(e.g., Nagata's  result on finite generation of invariants), and almost as an after thought, 
also yield geometric reductivity. 

Before Haboush settled the conjecture, the second author, in \cite{S2},
made an attempt to solve this conjecture following the quotient
space approach and had partial success. 
 In {\em ibid}.\,it is also shown that
constructing the Mumford or GIT quotient is equivalent to constructing a 
quotient by a proper equivalence relation on a projective variety;
in fact,  proving the conjecture  
is equivalent to showing that a natural line bundle on this
projective variety, which is known to be {\it nef}, is in fact 
{\it semi-ample} (i.e. a suitable power of this line bundle has no base
points).  Recently Sean Keel \cite{SK} has given a very interesting criterion 
for a nef line bundle on a projective variety to be semi-ample
in characteristic $p$, $p>0$. Using this result of Keel and
strengthening the methods of \cite{S2}, we give a proof of Mumford's
Conjecture in this paper.
  
Geometric reductivity of a reductive group $G$ is equivalent
to showing that the set $Y$ of equivalence classes of
semi-stable points for a linear action of $G$ on a projective
scheme $X$ has a canonical structure of a projective scheme 
(see below for the definitions and notations).  Roughly speaking 
one can say that the proof given here (of this equivalent form of
Mumford's conjecture) consists in checking the 
Nakai-Moishezon criterion for ampleness for a natural choice of a 
line bundle $L$ on $Y$.  A principal tool is the Hilbert-Mumford
criterion i.e. a process of reduction to a maximal torus or
even a 1-dimensional torus (Chapter 2, \cite{M}), for checking stability, 
semi-stability etc. of points.  An essential 
difficulty in this approach is that it is not easy to see, \`a 
priori, any natural scheme theoretic structure on $Y$.  This has to be
built up in stages, first as a topological space on which suitable
notions of properness, morphisms etc. have to be introduced,
eventually culminating in a scheme structure on $Y$ (which is is shown
to be projective via $L$).
In the case when ``stable = semi-stable'' this process is
simpler; it is easier to show that $Y$ is a proper scheme and
the proof is, indeed, checking the Nakai-Moishezon criterion
for $L$ on $Y$ (\cite{S2}).  In the general case one shows that
there is a projective scheme which surjects onto $Y$; in fact,
one can find such a projective scheme $Q$ which is ``generically
finite'' over $Y$ and one works with such ``models'' for $Y$
($Y$ is the quotient of $Q$ by a proper equivalence relation, as 
mentioned above).  The Nakai-Moishezon criterion is to be
interpreted as checking that $L$ (i.e. the pull-back of $L$ on
$Y$) is ``big'' on $Q$ and this is done by refining the 
methods of \cite{S2}.  However, this does not suffice to complete
the proof of geometric reductivity for one cannot expect $L$ 
to be ample on $Q$ but only ``semi-ample'' and one requires
some analogue of the Nakai-Moishezon criterion for semi-ampleness.
This is achieved by the work of Sean Keel (\cite{SK}) and appealing
to this work, the semi-ampleness of $L$ on $Q$ follows.
With a little more work, the required structure of a projective 
scheme on $Y$ also follows.  A more comprehensive outline of
proof is given in \cite{S5}.

One knows that geometric reductivity for reductive algebraic schemes holds 
over a general {\it base scheme}
(see \cite{S5} or Appendix G to Chap.~1, \cite{M}), the proof being again based on
the irreducibility of the Steinberg representation.  One would also like to
prove this by the quotient space approach as in this paper. For this one
requires a suitable generalisation of Sean Keel's result and this seems
to pose difficulties.

\subsection{Conventions, Notations and Definitions}\label{ss:conventions}
We work throughout over an algebraically closed field of {\em positive characteristic}
$p>0$. Thus, for example, a {\em variety} or a {\em scheme} will mean a $k$-variety
or a $k$-scheme respectively. A variety means a {\em separated reduced finite type} 
scheme (over $k$).

We fix a semi-simple algebraic group $G$ over $k$ (except briefly in \Sref{s:s=ss} where
we allow $G$ to be reductive algebraic). The aim of
the paper is to show that $G$ is geometrically reductive.

As is standard, ${\mathbb{G}}_m$ denotes the multiplicative group scheme of ``non-zero
scalars" in $k$, i.e.,
\[ {\mathbb{G}}_m \set \spec k[T,\,T^{-1}].\]

We will largely be working in a setting where $G$ acts linearly on a projective scheme (in fact,
more often than not, on a projective variety). It is convenient to have in place a terminology which
will act as a shorthand for certain recurring situations.

\begin{defis}\label{def:G-pair} A {\em linear $G$-pair}---or simply a {\em $G$-pair}---is 
a pair $(X,L)$ with $X$ a projective scheme and $L$ an ample line bundle on $X$ 
such that $G$ acts on $X$ and this action lifts to a linear action on $L$. The $G$-pair
$(X,L)$ is said to be {\em reduced} if  $X$ is reduced. A {\em linear
$G$-triple}--or simply a {\em $G$-triple}---$(X,\,L,\,\bp(V))$ consists of a $G$-pair $(X,\,L)$,
a finite dimensional rational $G$-module $V$, and a $G$-equivariant closed embedding
$X\hookrightarrow \bp(V)$ such that $L$ is the restriction of the tautological ample bundle
$\co(1)$ on $\bp(V)$. The $G$-triple $(X,\,L,\bp(V))$ is said to be {\em reduced} if (as before)
$X$ is reduced.
\end{defis}

Here are some other conventions and notations not listed above:

1) We repeat that by a variety we mean a separated reduced irreducible scheme of finite type. Points
on a scheme of finite type mean closed (whence $k$-rational) points.

2) If $V$ is a $k$-vector space, then we identify $V$ with the scheme $\spec{S(V^*)}$, where
$V^*$ is the dual of $V$ and $S(V^*)$ is the symmetric algebra on $V^*$.

3) If $(X,\,L,\,\bp(V))$ is a $G$-triple, we often denote the tautological ample bundle $\co(1)$ on
$\bp(V)$ by the letter $L$. If we need to distinguish between $L$ on $\bp(V)$ and $L$ on $X$,
we use the symbols $L_{\bp(V)}$ for the former and $L_X$ for the latter.

4) If $(X,\,L,\,\bp(V))$ is a $G$-triple, then, as is standard, $\widehat{X}$ will denote the
cone in $V$ over $X$. Note that $\widehat{X}$ is a closed $G$-stable subscheme of $V$.
Note also that $\widehat{\bp(V)}=V$.

5) The {\em homothecy action} on $\widehat{X}$, for a $G$-triple $(X,\,L,\,\bp(V))$ is the action
of $\bg_m$ on $\widehat{X}$ given by scalar multiplication.

6) A {\em one parameter subgroup} $\lambda$ of $G$ is a map of algebraic groups $\lambda\colon
\bg_m\to G$. We use the abbreviation 1-PS for ``one parameter subgroup". If $G$ acts algebraically
on a $k$-scheme and $\lambda$ is a 1-PS of $G$, then we often call the resulting action
of $\bg_m$ on $X$ as the {\em action of $\lambda$ on $X$}. 

\section{\bf{The Main strategy}}\label{s:main}

Recall that if $(X,\,L,\,\bp(V))$ is a $G$-triple, the semistable locus $X^{ss}=X^{ss}(L)$ of $(X,\,L)$
is the locus of points $x\in X$
such that if $\wh{x}$ is a point on the cone $\wh{X}\hookrightarrow V$ over $X$ which represents
$x$, then the orbit of $\wh{x}$ does not contain the origin $0\in V$ in its closure. If
the orbit of $\wh{x}$ is closed, and $\dim{{\wh{x}}G}=\dim{G}$, we say $x$ is a stable point and 
denote the stable locus $X^s$ or $X^{s}(L)$. We refer the reader to \Sref{s:basic-ss} for a more
extended discussion. In particular, there it is shown that the relation $x\sim x'$  whose graph is
$\{(x,\,x')\in X^{ss}\times_kX^{ss}\,\vert\, {\overline{xG}}={\overline{x'G}}\}$ is an equivalence
relation, the so-called {\em semi-stable equivalence relation}. (See \Dref{def:ss-equiv}.)
Our focus, as we have pointed out earlier, 
is to extend the techniques of \cite{S2} to
prove that our semi-simple algebraic group $G$ is geometrically reductive. This implies that every
reductive group is geometrically reductive.
In \cite[p.\,550,\,Theorem\,7.1]{S2} it is proven 
that if $(X,\,L)$ is a $G$-pair with $X$ {\em normal} and {\em projective}, whose stable locus
$X^s$ agrees with its semi-stable locus $X^{ss}$ (both loci with respect to $L$),
and $X={\rm{Proj}}\,R$
where $R=\bigoplus_{n\ge 0}\Gamma(X,\,L^n)$, then $R^G$ is a finite type $k$-algebra. Moreover,
if $Y={\rm{Proj}}\,R^G$, and
$\pi\colon X\to Y$ the {\em rational} map induced by the inclusion $R^G\hookrightarrow R$,
then the map $\pi$ is regular on the semi-stable locus $X^{ss}$, and the resulting map
$X^{ss}\to Y$ is a geometric quotient. It is also proven in \cite[p.\,553,\,Theorem\,7.2]{S2}
that if $x\in X$ is a stable point ($X$ not necessarily normal) then there exists a $G$-invariant
non-constant homogeneous polynomial $p$ on $X$
such that $p(x)\neq 0$. The problem is to extend
this result to semi-stable points. In this section we flesh out the main strategy and reduce
the problem to finding a map of stacks $Q\to [X^{ss}/G]$ (from a normal projective variety
$Q$) with certain properties, the most
important amongst them being that the resulting map from $Q$ to a natural stratified space $Y$
associated with the $G$-action on $X^{ss}$ is generically finite and the line bundle $L_Q$ on
$Q$ induced by $L$ is nef and big.

\subsection{Preliminaries} We begin by giving a few definitions.

\begin{defi}\label{def:saturated} The $G$-invariant map $X^{ss}\hookrightarrow \bp(V)$
is said to be {\em saturated},  if $X^{ss}\neq\emptyset$, and for semi-stably equivalent points 
$v$ and $v'$  in $\bp(V)$,
 $v\in X^{ss}$ implies that $v'\in X^{ss}$. We often simply say $X^{ss}$ is saturated when
 the $G$-invariant embedding into $\bp(V)$ is clear.
\end{defi}

Recall that a line bundle $L$ on a projective algebraic scheme is {\em nef} if
$\deg{L\vert_C}\ge 0$ for every irreducible curve on $X$, and it is {\em big} if for
$n\gg 0$ the {\em regular} global sections of $L^n$ define a {\em rational  birational} map.
In case $X$ is reduced and projective, and $L$ is
nef, then it is big if and only if $L^{(r_i)}\vert_{X_i} >0$ for every irreducible component
$X_i$ of $X$, with $r_i=\dim{X_i}$. (See \cite[VI.2.15 and VI.2.16]{K2}.)
Finally $L$ is {\em semi-ample} if some positive
power of it is generated by global sections. In other words, given a point $x\in X$ there is
a section of a positive power of $L$ which does not vanish at $x$. If $L$ is semi-ample (say
$L^n$ is generated by global sections) and
big and then the regular morphism on $X$ induced by $L^n$ is birational on to its image.

\begin{defi}\label{def:G-s-ample} Let $(X,\,L,\,\bp(V))$ be a $G$-triple. We say $L$ is 
{\em $G$-semi-ample on $X^{ss}$}, or $L$ is {\em $G$-semi-ample on $X$}, if given $x\in X^{ss}$,
there is a postive integer $n$ and an element $s\in \Gamma(X^{ss},\,L^n)^G$ such that 
$s(x)\neq 0$. Equivalently, given $\hat{x}\in {\widehat{X}}^{ss}$ there is a $G$-invariant
regular non-constant homogenous function $F$ on ${\widehat{X}}^{ss}$ such that 
$F(\hat{x})\neq 0$.
\end{defi}

\begin{defi}\label{def:usual} We say $\De=(X,\,L,\,\bp(V),\,X^{ss}\xrightarrow{\alpha} Y)$ is a
{\em quotient data}  if $(X,\,L,\,\bp(V))$ is a $G$-triple,
$X^{ss}\set X^{ss}(L)$,
$Y$ the {\em topological space} obtained by quotienting $X^{ss}$ by the semi-stable equivalence
relation, and $\alpha$ the resulting quotient map 
(cf.~\Dref{def:top-quot} in \Sref{s:basic-ss} below). We say
$\De$ is {\em saturated} if $X^{ss}\hra \bp(V)^{ss}$ is saturated, {\em reduced} if
$X$ is reduced, and {\em irreducible} if $Y$
is irreducible. $\De$ is a {\em strong quotient data} (or simply a {\em strong data}) if it is
reduced, saturated, $\bp(V)^s\neq \emptyset$, and $L$ restricted to the fibres
of $\alpha$ is trivial.{\footnote{This can be achieved by
replacing $L$ by a suitable power of $L$, as we show in \Rref{rem:strata}.}}
(The importance
of the conditions $\De$ saturated and $\bp(V)^s\neq \emptyset$ is seen in \Pref{prop:nonempty}.)
We often write $(X^{ss}/\negmedspace/G)_{\text{top}}$ for the quotient $Y$.
\end{defi}

Note that if $\Y\set (\bp(V)^{ss}/\negmedspace/G)_{\text{top}}$ and 
$\tilde{\alpha}\colon \bp(V)^{ss}\to \Y$ 
is the resulting quotient map of topological spaces, then $X^{ss}\hookrightarrow \bp(V)$ is
saturated if and only if $\tilde{\alpha}^{-1}(Y)=X^{ss}$, where we regard $Y$ as a closed subspace
of $\Y$ in a natural way.

\begin{defi}\label{def:gen-quot} Let $\De=(X,\,L,\,\bp(V),\,X^{ss}\xrightarrow{\alpha} Y)$ be a 
quotient data. A {\em non-empty} scheme $U$ is said to be a {\em generic quotient}
$\De$ as above if the underlying topological space of $U$ is an open dense subset of $Y$
such that $\alpha^{-1}(U)\xrightarrow{{\text{via}}\,\alpha} U$ is a morphism of schemes
($\alpha^{-1}(U)$ having the canonical open $X^{ss}$-subscheme structure on the open set
$\alpha^{-1}(U)$), and such that $\co_U$ is the sheaf of $G$-invariant sections of the direct
image of $\co_{\alpha^{-1}(U)}$.
\end{defi}

We show later (see \Lref{lem:U}) that if $\De$ is a reduced irreducible quotient data, then a generic
quotient for $\De$ exists. A little thought shows that by working with each irreducible component,
and removing all points which lie in more than one irreducible component, a generic quotient
exists for all reduced quotient data, whether irreducible or not.

Given quotient data $\De$, we write $X_\De$, $L_\De$, $\bp(V_\De)$, $\alpha_\De$, $Y_\De$
etc. for the various datum comprising $\De$.

If $\De$ satisfies all the requirements of a strong quotient data except the requirement that $L$
is trivial on the fibres of $\alpha$, then, this requirement is easily achieved by replacing
$L$ by a suitable positive power of itself. Indeed, by \Lref{lem:U} this is achieved on a dense
open subset of $Y$, and working with the complement of $U$ in $Y$, and continuing the
process, by noetherian induction we achieve what we wish. One consequence is that
we have a line bundle $L_{/G}$ on the algebraic stack $[X^{ss}/G]$ (see
remarks in \ref{conv:L-power}). 
We point out that $X^{ss}(L^n)=X^{ss}(L)$ (resp. $X^s(L^n)=X^s(L)$) for
$n$ postive as is readily verified among other methods by the Hilbert-Mumford criterion
\cite[p.\,520,\,Theorem\,2.2]{S2} (see also [{\em Ibid.},\,p.\,519,\,Proposition\,2.1\,(2)]).

Given a reduced quotient data $\De$ it can be strengthened to a strong data with very little effort.
Indeed replace $V_\De$ by $V_\De\oplus W$ with $W$ a rational finite dimensional $G$ such
that $\bp(W)^s\neq \emptyset$, and $X_\De$ by the closure of the inverse image of $Y$ under
the semi-stable quotient map on $(\bp(V)^{ss}/\negmedspace/G)_{\text{top}}$ and finally replacing
$L_\De$ by suitable positive power of itself, we get a strong quotient data. 

\subsection{Zariski locally trivial principal $G$-bundles}\label{ss:zariski-triv} 
The problem of showing the geometric
reductivity of $G$ is equivalent to  showing that if $\De$ is a strong quotient data, then
$L_\De$ is $G$-semi-ample on $X_\De$. Indeed it is enough to show that $\co(1)_{\bp(V)}$
is $G$-semi-ample for a finite dimensional $G$-module $V$,  and we note that the quotient data 
$(\bp(V),\,\co(1)_{\bp(V)},\,\bp(V)^{ss}\xrightarrow{\wt{\alpha}}\Y)$ is strong.
In this section we work with strong {\em irreducible} quotient data $\De$
and reduce the problem of showing $L_\De$ is $G$-semi-ample to that of finding a Zariski
locally trivial principal $G$-bundle $P\to Q$ and a $G$-invariant map $P\to X^{ss}$ such that
$Q\to Y$ is ``generically finite" (see \Dref{def:gen-finite}), 
and such that if $f\colon Q\to [X^{ss}/G]$ is the classifying
map, then $f^*L_{/G}$ is nef and big. Here $L_{/G}$ is the line bundle on $[X^{ss}/G]$
defined by $L=L_\De$. Other technical hypotheses are required to be
satisfied by the principal bundle $P\to Q$. Here is what is needed.

Suppose $\De=(X,\,L,\,\bp(V),\,X^{ss}\xrightarrow{\alpha} Y)$ is a
strong {\em irreducible} quotient data. We quickly summarize what we
need from \Sref{s:strata} so that this section can be followed. Note that $Y$ is a stratified space. 
Indeed by \Lref{lem:U} we can find a generic quotient $U$
for $\De$, and by working with the complement $Y'=Y\setminus U$, we get a strong quotient data
$\De'$ where $X'=X_{De'}$ is the closure in $X$ of $\alpha^{-1}(Y')$. Therefore we can find
a generic quotient $U'$ of $\De'$ and $U'$ is an open subset of $Y$. Proceeding in this manner
we have stratified space (cf. \Rref{rem:strata}). The original open stratum is called the {\em big
stratum} for this stratification.
If $Q$ is an algebraic scheme and $f\colon Q\to [X^{ss}/G]$ 
is a map of stacks, then (as $Q$ is the base of a principal $G$-bundle) we have an obvious 
continuous map $q\colon Q\to Y$,and this map is a stratified map (roughly, over each stratum
we have a map of schemes, with the inverse image of a non-big stratum being given the reduced
structure). We can regard $Y$ itself as a stack (see \S\S\,\ref{ss:Ystack}, especially 
\eqref{eq:Ystack}), though it is not an algebraic stack as defined (lacking a smooth atlas)
whence the above considerations on $Q$ define a map $\gamma\colon [X^{ss}/G]\to Y$.
Finally since $\De$ is a strong data, $L$ is trivial on the fibres of $\alpha$, whence its pull
back to the principal bundle defined by $f\colon Q\to [X^{ss}/G]$ descends to $Q$. This means
we can talk of a line bundle $L_{/G}$ on $[X^{ss}/G]$ which should be regarded as the line
bundle to which $L\vert_{X^{ss}}$ descends.

In what follows, let $\wt{X}$ be the normalization of $X$, $\wt{L}$
the pull back of $L$ to $\wt{X}$, and ${\wt{X}}^{ss}={\wt{X}}^{ss}{\wt{L}}$.

We will show in \Sref{s:conclusion} that there
exists an {\em irreducible normal projective variety} $Q$ together with a map of
algebraic stacks
\stepcounter{thm}
\begin{align*}\label{map:classifying}\tag{\thethm}
f\colon Q\to [X^{ss}/G]
\end{align*}
such that:
\begin{enumerate}
\item \label{beta}The principal $G$-bundle $\beta\colon P\to Q$ corresponding to $f$ is Zariski 
locally trivial. 
\item The map $q\set\gamma\smcirc f\colon Q\to Y$ is ``generically finite". In other words, if
$U$ is a generic quotient for $\De$ then the map of schemes $q^{-1}U\to U$ is generically
finite (see also \Dref{def:gen-finite}). By \Cref{cor:proper} it follows that $Q\to Y$ is surjective.
\item Set $L_Q= f^*L_{/G}$. Then $L_Q$ is nef and big on $Q$. 
\item \label{positive} If $C$ is a closed integral curve in $Q$ such that $q\vert_C$ is non-constant, 
then $\deg(L_Q\vert_C)>0$. 
\item \label{gam-orb}Let $k(Y)\set k(U)$. Then $k(Q)$ is normal over $k(Y)$ and the finite group 
$\Gamma={\mathrm{Aut}}_{k(Y)}(k(Q))$ 
acts on $Q$, and $q\colon Q\to Y$ is 
$\Gamma$-invariant for the trivial action of $\Gamma$ on $Y$. There exists a generic quotient
$U$ for the data $\De$ such the fibres of $q$ over $U$ are $\Gamma$-orbits.
\item\label{lift} The action of $\Gamma$ lifts to $L_Q$, whence a positive power of $L_Q$ 
descends to ${\overline{Q}}\set Q/\Gamma$. For definiteness, suppose $r$ is positive and
$L_Q^r$ descends to $\eL$ on ${\overline{Q}}$.
\setcounter{savectr}{\value{enumi}}
\end{enumerate}
Note that $P$ is realized as the ``base change" $P\set Q\times_{[X^{ss}/G]}X^{ss}$, and
$\beta$ is the projection to the first factor. $P$, being smooth over the normal variety $W$, is itself
normal. We also have a second projection which is a
$G$-invariant map $\pi\colon P\to X^{ss}$. Note that  
$\Gamma({\overline{Q}},\eL^n)\iso \Gamma(Q,\,L_Q^{rn})^\Gamma$ for every postive integer $n$.
In addition to the conditions above, the map \eqref{map:classifying} also satisfies:
\begin{enumerate}
\setcounter{enumi}{\value{savectr}}
\item\label{G-gam-inv}The map $\pi\colon P\to X^{ss}$ takes values in an irreducible component
of $X^{ss}$, say $X_1^{ss}$ and the map $P\to X_1^{ss}$ is dominant. Moreover
${\widetilde{\pi}}\colon P\to 
{\widetilde{X}}^{ss}$ is the map induced by $\pi$, then the sections of $\eL^n$, up to a suitable
power of $p$, can be identified with the $\Gamma$-invariant sections
of the pull back under ${\widetilde{\pi}}\colon P\to {\widetilde{X}}^{ss}$ of $G$-invariant sections of 
${\wt{L}}^{rn}$. In other words, if  $t\in \Gamma({\widetilde{X}}^{ss},\,{\wt{L}}^{rn})^G$
then a suitable $p$-power of $t$, say $t^n$ (with $n=p^m$), is in the image of the  map
\[\Gamma({\widetilde{X}}^{ss},\,{\wt{L}}^{rn})^G\to \Gamma({\overline{Q}},\,\eL^n)\]
which is the composite
\begin{align*}
\Gamma({\widetilde{X}}^{ss},\,{\wt{L}}^{rn})^G\to
[\Gamma({\widetilde{X}}^{ss},\,\wt{\pi}_*{\wt{\pi}}^*{\wt{L}}^{rn})^G]^\Gamma  
& \lr \,[\Gamma(P,\,{\wt{\pi}}^*{\wt{L}}^{rn})^G]^\Gamma\\
&  {\phantom{X}}\set \,\Gamma(Q,\,L^{rn}_Q)^\Gamma \\
& \iso \Gamma({\overline{Q}},\,\eL^n).
\end{align*}
\end{enumerate}

\subsection{Geometric Reductivity} We now show, assuming the existence of the map
\eqref{map:classifying} for strong irreducible quotient data, that Geometric Reductivity of $G$ 
holds. More precisely, we show that if $\De$ is a strong quotient data then $L_\De$ is
$G$-semi-ample on $X_\De$.

Fix a strong quotient data $\De=(X,\,L,\,\bp(V),\,X^{ss}\xrightarrow{\alpha} Y)$.

\begin{lem}\label{lem:factorization1} Suppose $\De$ is irreducible.
Let $f\colon Q\to [X^{ss}/G]$ be the map \eqref{map:classifying} satisfying conditions (1)--(7)
of \S\S\ref{ss:zariski-triv}.
Assume $L_Q$ descends to a line bundle $L_{\overline{Q}}$ on ${\overline{Q}}$. If
$L_Q$ is semi-ample on $Q$ then
\begin{enumerate}
\item $L_{\overline{Q}}$ is semi-ample.
\item There is a factorization of $q\colon Q\to Y$ given by the commutative diagram
\stepcounter{thm}
\begin{equation*}\label{diag:factorization1}\tag{\thethm}
{\xymatrix{
Q \ar[d]_{q} \ar[r]^{\varphi} & {\overline{Q}} \ar[ld]_{\ol{q}} \ar[d]^{\psi} \\
Y & \ar[l]^{q_{{}_W}} W
}}
\end{equation*}
where $\varphi$ is the natural quotient map $Q\to Q/\Gamma$, $\psi$ is birational and $W$ is
projective and normal, 
with $L_{{}_{\overline{Q}}}^n$ descending to an ample line bundle $L_W$ on $W$ for
a suitable positive integer $n$,
and $q_W$ has finite fibres. Moreover, if $\bar{q}\colon {\overline{Q}}
\to Y$ is the composite $\bar{q}=q_W\smcirc \psi$, there exists a generic quotient $U$
of $\De$ such that the maps
$({\bar{q}})^{-1}(U)\to U$ and
$q_W^{-1}(U)\to U$ are bijective continuous maps, 
and $({\bar{q}})^{-1}(U)\to q_W^{-1}(U)$ is an isomorphism of schemes.
\end{enumerate}
\end{lem}

\proof Since $Q\to Y$ is $\Gamma$-
invariant for the trivial action of $\Gamma$ on $Y$, we have a map ${\bar{q}}\colon {\overline{Q}}
\to W$ such that ${\bar{q}}\smcirc \varphi= q$. According to property\,\eqref{gam-orb} in 
\S\S\ref{ss:zariski-triv}, we have a generic quotient $U$ for $\De$ such that
$q^{-1}(U)\to U$ has fibres which are $\Gamma$-orbits. It follows that
$({\bar{q}})^{-1}(U)\to U$ is bijective.

Now $L_{\overline{Q}}$ is semi-ample. Moreover it
is big and nef since $L_Q$ is big and nef. By replacing $L_{\overline{Q}}$ by a positive
power of itself, we may assume $L_{\overline{Q}}$ is base point free. Since $L_{\overline{Q}}$
is big, the projective map  induced by it is birational on to its image, whence we have
the birational map $\psi\colon {\overline{Q}}\to W$ and $L_{\overline{Q}}$ descends to an 
ample line bundle $L_W$ on $W$. In fact
\stepcounter{thm}
\begin{equation*}\label{eq:projS}\tag{\thethm}
W= {\mathrm{Proj}}\,S
\end{equation*}
where $S=\oplus_{n\ge 0} S_n$ is the graded ring given by 
$S_n=\Gamma({\overline{Q}},L^n_{\overline{Q}})$.
Moreover, $W$ is normal, since ${\overline{Q}}$
is. In particular, the fibres of $\psi\colon \ol{Q}\to W$ are connected. For $w\in W$, we claim that 
${\bar{q}}(\psi^{-1}(w))$ is a point in $Y$. Suppose not. Then we have points $a$ and $b$
in $\psi^{-1}(w)$ such that ${\bar{q}}(a)\neq {\bar{q}}(b)$. Let $C$ be a closed integral
curve in the connected space $\psi^{-1}(w)$ which passes through $a$ and $b$. Then
$C$ is generically finite on to its image $C'$ in $Y$. By condition\,\eqref{positive} satisfied
by $Q$ (see \S\S\,\ref{ss:zariski-triv}) 
we see that $\deg{L_{\overline{Q}}}\vert_C >0$. On the other hand $L_{\overline{Q}}=
\psi^*L_W$, and hence $L_{\overline{Q}}\vert_C$ is a trivial line bundle. This gives us a
contradiction. We therefore deduce a map $q_W\colon W\to Y$ such that the diagram above
commutes. It remains to show that $q_W$ has finite fibres. If not, we have a closed integral
curve $C$ in $W$ which maps to a point in $y$ in $Y$. Since $L$ is trivial on ${\alpha}^{-1}(y)$,
it follows that $L_{\overline{Q}}$ is trivial on the proper transform of $C$ in ${\overline{Q}}$
by property\,\eqref{G-gam-inv} enjoyed by the map $f$ in \eqref{map:classifying}, whence 
$\deg(L_W\vert_C)=0$.
On the other hand $L_W$ is ample. This is a contradiction.
The assertions that $(\bar{q})^{-1}(U)\to q_W^{-1}(U)$ is an isomorphism of schemes and
that $q_W^{-1}(U)\to U$ is bijective are obvious for the generic quotient $U$ asserted in
\eqref{gam-orb} of \S\S\ref{ss:zariski-triv}.
\qed

\begin{lem}\label{lem:proj} Let $\De$ be irreducible. Let $f\colon Q\to [X^{ss}/G]$ be the map 
\eqref{map:classifying} satisfying conditions (1)--(7)
of \S\S\ref{ss:zariski-triv}.
Assume $L_Q$ descends
to a line bundle $L_{\overline{Q}}$ on ${\overline{Q}}$. If $L$ is $G$-semi-ample on $X^{ss}$
then
\begin{enumerate}
\item $L_Q$ is semi-ample.
\item If $R=\oplus_{n\ge 0} R_n$ is the graded ring whose $n^{\rm{th}}$-graded piece is
given by $R_n=\Gamma(X^{ss},\,L^n)^G$, then
$R$ is a finitely generated $k$-algebra.
\item We have a morphism of schemes ${\overline{Q}}\to {\mathrm{Proj}}\,R$ and a bijective
continuous map ${\mathrm{Proj}}\,R \to Y$ fitting into a commutative diagram
\[
{\xymatrix{
Q \ar[r]^{\varphi} \ar[d]_{q} & {\overline{Q}}\ar[d] \\
Y & {\mathrm{Proj}}\,R \ar[l]
}}
\]
\end{enumerate}
\end{lem}

\proof The assertion that $L_Q$ is semi-ample is straightforward. Indeed, given a point $a\in Q$,
we can find point $b\in P=X^{ss}\times_{[X^{ss}/G]}Q$ a positive integer $n$ and a $G$-invariant 
section $s$ of $L^n\vert_{X^{ss}}$ such that $s(\pi(b))\neq 0$, whence the pull-back
$\pi^*s$ is non-vanishing on $b$. The corresponding section $\sigma$ of $L_Q^n$ is such that
$\sigma(a)\neq 0$. 
From \Lref{lem:factorization1} it follows that $L_{\overline{Q}}$ is semi-ample. By replacing $L$
by a suitable power if necessary, we assume that $L_{\overline{Q}}$ is actually base point free.
We have the commutative diagram \eqref{diag:factorization1} with $W={\mathrm{Proj}}\,S$
as in \eqref{eq:projS}. Note that $S_n=\Gamma(W,\,L_W^n)$. 

Now $R\hookrightarrow S$. In fact $R$ is a graded $k$-subalgebra of $S$. By our hypothesis,
we can find a finite number of $G-\Gamma$-invariant sections of 
a suitable power $L^n$ of $L$ such that their
pull backs to $P$ considered as sections of powers of $L_W^n$ have no base points. Let
$R'$ be the finitely generated $k$-algebra generated by these sections. Then $R'\to S$ is
finite and since $R'\hookrightarrow R \hookrightarrow S$, $R$ is a finite $R'$-algebra, whence
a finitely generated $k$-algebra. Since $q_W\colon W\to Y$ has finite fibres, therefore for
any $y\in Y$, $L_W\vert_{q_Q^{-1}(y)}$ is trivial, whence the map $W\to {\mathrm{Proj}}\,R$
sends $q_W^{-1}(y)$ to a single point in ${\mathrm{Proj}}\,R$. It follows that we have a 
bijective continuous map ${\mathrm{Proj}}\,R \to Y$.
\qed

\begin{rem}{\em Consider diagram \eqref{diag:factorization1}. Since ${\overline{Q}}=Q/\Gamma$
and $\Gamma$ is the group of $k(Y)$ automorphisms of $k(Q)$,
it is clear that $k(Y)\to k(W)=k({\overline{Q}})$ is a purely inseparable field extension.}
\end{rem}

\begin{lem}\label{lem:red-to-irr} Let $Y^{(i)}$ be the distinct irreducible
components of $Y$ and $X^{(i)} = \alpha^{-1} (Y^{(i)})$ (with its
reduced structure) $1 \leq i \leq r$.  Suppose that $L$ is $G$-semi-ample
on $X^{(i)}$ for $1 \leq i \leq r$.  Then $L$ is $G$-semi-ample on $X$.
\end{lem}

\proof
 Let us note the following general fact.  Let $S$ be
a projective variety and $M$ an ample line bundle on $S$.  Let $T$
be a closed subscheme of $S$ with ideal sheaf $I = I(T)$. 
Then given a section $s$ of $M|_T$, $s^r$ can be extended to a 
section of $M^r$ for $r \gg 0$.  This follows from the fact
$H^1 (S, M^r \otimes I) =0$ for $r \gg 0$ and writing the usual 
exact sequence.

The lemma is proved by induction on $r$ and hence we suppose that 
it is true for $X' = X^{(1)} \cup \cdots \cup X^{(r-1)}$ (scheme
theoretic union).  Let $X'' = (X')^{ss} \cap (X^{(r)})^{ss}$. 
Write $L_i$ for $L\vert_{X^{(i)}}$ and $L''$ for $L\vert_{X''}$.
Then the image of $X''$ in $Y^{(r)}$ is a proper closed subset
of $Y^{(r)}$.  Let $s$ be a $G$-invariant section, say of $L''$
and $s_1$ its restriction to $X''_{\mathrm{red}}$.  Let $R =R (X^{(r)})$ be
the graded ring such that $R_n =\Gamma(X^{(r)},\,L_r^n)^G$. Note
that $\proj R= Y^{(r)}$ as a topological space. Write $L_R$ for the natural ample
line bundle on $\proj R$ defined $L$.
Then the image of $X''$ in ${\mathrm{Proj}}\,R$ is a closed subscheme $T$ of $\proj R$
with the reduced structure.  Then we see that $s_1$ can be identified
with a section of  $L_R\vert_T$ (to do this apply \Lref{lem:proj}).  Then by the initial 
remark $s_1$ can be extended
to a section of $L_R$ on ${\mathrm{Proj}}\,R$ (we may have to take a power of $L$
and also raise $s$, to the same power) and we identify this 
extension as a $G$-invariant section $t$ of $L$ on $(X^{(r)})^{ss}$.
Thus we see that $s$ and $t$ coincide on $(X'')^{ss}_{{\mathrm{red}}}$ and 
then raising to the power $p^a$, $a \gg 0$, we can suppose by \Lref{lem:bijective1}
that  $s$ and $t$ coincide on $X''$.  Thus without loss of generality 
we can say that $s$ can be extended to a $G$-invariant section
on $X^{ss}$.  Besides, we can also achieve this extension so
that it does not vanish at a point of $Y^{(r)}$ outside $T$.
From this the lemma follows easily
\qed

We will need the following result of Keel (see \cite[p.\,254,\,Theorem\,0.2]{SK}):

\begin{thm}\label{thm:keel}Let $S$ be a projective variety of dimension $n$ 
over an algebraically closed field of characteristic $p>0$.{\footnote{The reader is 
reminded that this is the default assumption on our underlying field $k$ throughout the paper. 
We have re-stated this assumption in this statement, to emphasize that this is where the 
assumption of positive characteristic enters in a crucial way in this paper.}}
Let $L$ be a nef and big line bundle on $S$. Suppose that for any 
proper closed subset $T$ of $S$, $L\vert_T$ is semi-ample for the reduced scheme structure
on $T$. Then $L$ is semi-ample on $S$.
\end{thm}

\proof This is Keel's work already mentioned in the introduction. In greater detail, let $S$ be
any projective scheme and $L$ a nef line bundle on $S$. An irreducible subvariety $E$ os $S$
is called {\em exceptional} if $L\vert_E$ is {\em not} big. The {\em exceptional locus}
$E(L)$ of $L$ is defined to be the closure with the reduced structure of the union of all the
exceptional varieties. Then the result in \cite[p.\,254,\,Theorem\,0.2]{SK} is that $L$ is semi-ample 
on $S$ if and only if $L\vert_{E(L)}$ is semi-ample.
\qed

\begin{thm}\label{thm:main} $L$ is $G$-semi-ample on $X^{ss}$.
\end{thm}

\proof
By \Lref{lem:red-to-irr} it is enough to prove the theorem when $Y$ is irreducible. We prove this
by induction on $\dim{Y}$. In other words we suppose that if $\De^*$ is an irreducible strong 
quotient data with $\dim{Y_{\De^*}}<\dim{Y}$, then $L_{\De^*}$ is $G$-semi-ample on
$X_{\De^*}$. By \Lref{lem:red-to-irr} it follows that if $Y'$ is any proper closed subset of $Y$, and if
$X'=\alpha^{-1}(Y')$ is endowed with the reduced structure of a closed subscheme of $X$, then
then $L\vert_{X'}$ is $G$-semi-ample on $X'$.

Since $\De$ has been assumed to be irreducible, we have a normal irreducible projective variety
$Q$ and a map $f\colon Q\to [X^{ss}/G]$ satisfying (1)---(7) in \S\S\,\ref{ss:zariski-triv}.

By the induction hypothesis, it follows that $L_Q$ restricted to any proper closed subvariety of
of $Q$ is semi-ample (we need property \eqref{G-gam-inv} of \S\S\ref{ss:zariski-triv}). 
It follows then, by
\Tref{thm:keel} that $L_Q$ is semi-ample on $W$. Lemma \ref{lem:factorization1} can therefore
be applied. Consider the resulting commutative diagram 
\eqref{diag:factorization1}, namely:
\[
{\xymatrix{
Q \ar[d]_{q} \ar[r]^{\varphi} & {\overline{Q}} \ar[d]^{\psi} \\
Y & \ar[l]^{q_{{}_W}} W
}}
\]
According to {\em loc.cit.} we have a generic quotient $U$ for $\De$ such that, with
$q'=q_w\smcirc \psi$, $(q')^{-1}(U)\to U$ and $q_W^{-1}(U)\to U$ are bijective continuous maps.
Let $C=Y\setminus U$ and $Z=q_W^{-1}(C)$. Set $X^{ss}_{C}
=\alpha^{-1}(C)$ and let $X_C$ be its closure in $X$. Give $Z$ and
$X^{ss}_C$ their reduced scheme structures and write $L_Z=L_W\vert_Z$. According to 
\Lref{lem:proj} $C$ exists as a reduced scheme, and 
$L\vert_{X^{ss}_C}$ descends to an ample bundle $L_C$ on $C$. Clearly $L_Z$ descends
to $L_C$. Let
$q_Z\colon Z\to C$ be the resulting map induced by $q_W$.
Consider the diagram,
with the solid arrows already defined:
\stepcounter{thm}
\begin{equation*}\label{diag:CZdotted}
{\xymatrix{
Z \ar@{^(->}[r]^i  \ar[d]_{q_{{}_Z}} & W \ar@{.>}[d]^{g} \\
C \ar@{^{(}.>}[r]^{i'} & {\overline{W}}
}}
 \tag{\thethm}
\end{equation*}
We claim we can find a variety ${\overline{W}}$ whose underlying topological space is $Y$,
such that the dotted arrows in \eqref{diag:CZdotted} can be filled to make the diagram commute,
and such that $L_W$ descends to an ample line bundle $L_{\overline{W}}$ on ${\overline{W}}$
with the property that $L_{\overline{W}}\vert_C=L_C$.

The above claim is an easy consequence of the fact that for large positive $n$
 sections of $L^n_Z$ which come from
sections of $L^n_C$ can be extended to sections of $L^n_W$ such that these extensions separate
points of $W\setminus Z=q_W^{-1}(U)$. For, if such extensions exist, they give a base point free
linear system on $W$, which is a sub-linear system of the very ample linear system given by 
$L^n_W$. We now prove the existence of such extensions of such sections of $L_Z=q_Z^*L_C$.
Let $x_1, x_2\in q_W^{-1}(U)$, with $x_1\neq x_2$, and set 
\[Z'=Z\cup\{x_1,\,x_2\}\]
with $Z'$ being given the reduced scheme structure. Let $I$ be the ideal sheaf of the closed
subscheme $Z'$ of $W$. Let $L_{Z'}=L_W\vert_{Z'}$.
Then we have an exact sequence of coherent $\co_{W'}$-modules
\[0\to L^n_W\otimes I \to L^n_W \to L^n_{Z'}\to 0.\]
Let $\sigma\in \Gamma(C,\,L^n_C$. The section $q_Z^*(\sigma)\in\Gamma(Z,\,L^n_Z)$ can
be extended to a section $s$ of $L^n_{Z'}$ by setting $s(x_1)\neq 0$ and $s_(x_2)=0$
(for $L^n_{Z'}=L^n_Z\oplus k_{x_1}\oplus k_{x_2}$, where $k_{x_i}$ is the residue field
of $\co_{W,x_i}$ for $i=1,2$).  Since $H^1(W,\,L^n_W\otimes I)=0$ for $n\gg 0$, the
short exact sequence of sheaves displayed above shows that $s$ can be extended to all
of $W$ to give a section of $L_W$ provided $n$ is large enough.

Clearly the line bundle $L_{\overline{W}}$ is ample on ${\overline{W}}$ and the underlying
topological space of ${\overline{W}}$ is $Y$. It is now evident that the dotted arrows in
\eqref{diag:CZdotted} can be filled as required.

Since $P (\set Q\times_{[X^{ss}/G]}X^{ss})$ is irreducible therefore there exists a unique 
component $X_1^{ss}$ of $X^{ss}$ on which $p\colon P\to X^{ss}$ takes values. Let
${\wt{X}}_1^{ss}$ and ${\wt L}_1$ denote, respectively, the normalization of $X^{ss}_1$ and
the pull back of $L$ to ${\wt X}_1^{ss}$. Note that the map ${\wt p}\colon P\to {\wt X}^{ss}$
actually takes values in ${\wt X}_1^{ss}$. Let $s\in\Gamma({\overline{W}},\,L^n_{\overline{W}})$.
Then its pull-back $s^*\in \Gamma(W,\,L^n_W)$ arises from a $G$-invariant section
$\sigma_s\in \Gamma({\wt X}^{ss}_1,\,{\wt L}^n_1)^G$, and whence identifies with
a regular function on the normalization of the cone ${\wh X}_1^{ss}$ over $X_1^{ss}$. The graph
of this section defines a closed subset $\Gamma$ of ${\wh X}^{ss}_1\times_k{\mathbb A}^1$.
Since $\sigma_s$ arises from $s\in\Gamma({\overline{W}},\,L)$, which is topologically $Y$,
the projection $\Gamma\to {\wh X}_1^{ss}$ is proper and bijective. By \Lref{lem:bijective1}
suitable power of $p$ of this function on the normalization goes down to a regular function
on ${\wh X}_1^{ss}$. Yet another power of this section extends to $X^{ss}$. Such sections
achieve the required $G$-semi-ampleness.
\qed

\begin{cor}\label{cor:fin-gen} Let $S=\oplus_{n\ge 0}S_n$ be the graded ring defined by
$S_n=\Gamma(X^{ss},\,L^n)$. Then $R\set S^G$ is a finitely generated graded ring, and 
$Y$ acquires a canonical scheme structure $Y={\proj R}$.
\end{cor} 

\proof 
By \Lref{lem:proj}, this follows when $X$ is irreducible. Then by a devissage argument as in
\cite{S3}, the corollary follows. One formulates a more general statement that if $M$ is
a $G$-coherent module on $\bp(V)$ then $M^G$ is finitely generated. We leave the details
to the reader.
\qed

\section{\bf{Basic properties of Semi-stable equivalence}}\label{s:basic-ss}

Throughout this section, we fix a $G$-triple $(X,\,L,\,\bp(V))$. 
In this section we study the space $Y$ of semi-stable
equivalence classes as a topological space and show
that it has the expected properties which are consequences
of geometric reductivity.
For example, granting the notions of semi-stability, unstability etc., as given below,
if $X^{ss}$ denotes the semi-stable locus of $X$ with respect
to $L$, then in \Pref{prop:G-sep} and \Cref{cor:proper} we show,
respectively, that the graph of the semi-stable equivalence relation
in $X^{ss}\times_kX^{ss}$ is closed and that  it is ``proper'' i.e. ``$X^{ss} \mod
G$ is proper''. As a consequence, $Y$ is {\em separated}---i.e. the diagonal
map $Y\to Y\times Y$ is closed with respect to the topology on $Y\times Y$ induced by
the action of $G\times_kG$ on $X^{ss}\times_kX^{ss}$---and $Y$ is ``proper".
Moreover, in \Pref{prop:L-trivial} we show that a suitable power
of $L$, when restricted to a semi-stable 
equivalence class, is trivial. This is the first step towards
seeing that a power of $L$ ``descends to $Y$".  

The key result for all these is \Pref{prop:U(S)-closed}, i.e. that the ``$S$-unstable locus" 
$\cu (S)$ is closed in $X^{ss}$ (for suitable $G$-stable subsets $S$ of $X^{ss}$). This
corresponds to the expected result that for the quotient map $\alpha : X^{ss} \lr Y$,
$\alpha (S)$ is closed.  An important technical point concerns 
extensions of functions by taking their $p$-th powers \Lref{lem:bijective1}.  
This is probably well-known to
experts but we give a proof here as it is an important
point where ${\text{char}}\,k=p>0$ is used.

\subsection{Semi-stablity, polystability, stability, and unstability} Recall that if $(X,\,L,\,\bp(V))$
is a $G$-triple, then a point $\wh{x}$ of $\widehat{X}$ is said to be {\em unstable} if
the vertex $0\in {\widehat{X}}$ of the cone $\widehat{X}$ lies in the orbit closure
${\overline{\wh{x}G}}$. A point on $\widehat{X}$ which is not unstable is called {\em semi-stable}.
The locus of semi-stable points is denoted by ${\widehat{X}}^{ss}$. 
A point $x\in X$ is said to be semi-stable if for some 
$\wh{x}\in {\widehat{X}}\smallsetminus\{0\}$ lying over
$x$, $\wh{x}$ is semi-stable. Since the homothecy action on ${\widehat{X}}$ commutes with the
$G$-action on ${\widehat{X}}$, $x\in X$ is semi-stable if and only if every point $\wh{x}\in{\widehat{X}}
\smallsetminus\{0\}$ lying over $x$ is semi-stable. The semi-stable locus in $X$ is denoted
$X^{ss}(L)$, or, if the line bundle $L$ is understood from the context, simply $X^{ss}$.

A point $\wh{x}\in \widehat{X}$ is said to be {\em polystable} if $\wh{x}\neq 0$ and the $G$-orbit through
$\wh{x}$ is closed. We say $x\in X$ is polystable if for some (and hence all) $\wh{x}\in {\widehat{X}}
\smallsetminus\{0\}$ lying over $x$, $\wh{x}$ is polystable. 

We say $\wh{x}\in {\widehat{X}}$ is {\em stable}, or {\em properly stable}, if the $G$-orbit $\wh{x}G$ is
closed in ${\widehat{X}}$ and $\dim{\wh{x}G}=\dim{G}$. Equivalently, $\wh{x}\in {\widehat{X}}$ is stable
if the orbit morphism $G\to {\widehat{X}}$, $g\mapsto \wh{x}g$ is proper. We denote by 
${\widehat{X}}^s$
the $G$-stable locus of stable points in ${\widehat{X}}$. A point $x\in X$ is said to be stable if
for some (and hence all) $\wh{x}\in {\widehat{X}}\smallsetminus\{0\}$ lying over $x$, $\wh{x}$ is stable.
We denote by $X^s(L)$ (or, if $L$ is understood from the context, by simply $X^s$) the
locus of stable points in $X$. Note that a stable point is polystable.

The notion of an unstable point can be generalized as follows:

\begin{defi}\label{def:S-unstable} Let $(X,\,L,\,\bp(V))$ be a $G$-triple.
Let $S$ be a closed $G$-invariant subset of $\wh{X}$ (e.g. 
$\wh{X} =V$), which we can endow with the reduced structure.
Following Kempf \cite{GK}, we say that a point $\wh{x} \in \wh{X}$ is $S$-unstable if 
the orbit closure $\ol{\wh{x}G}$ meets
$S$.  We denote by $\cu(S, \wh{X})=\cu(S)$ the set of $S$-unstable points in $\wh{X}$.
If $\lambda\colon\bg_m\to G$ is a 1-PS of $G$, then the locus of $S$-unstable points
in $\wh{X}$ under the action of $\lambda$ on $\wh{X}$ is denoted $\cu(S,\,\lambda)=
\cu(S,\,\wh{X},\,\lambda)$. 
\end{defi}

\begin{rems} {\tiny{${}_.$}}
{\em
\begin{enumerate}
\item[(a)] Note that a point in $\wh{X}$ is $(0)$-unstable if and only if it is unstable.
\item[(b)] $\cu (S,\wh{X}) = \cu (S,V)\cap \wh{X})$.
\item[(c)] The set $\cu(S)$ is again $G$-stable, and $S\subset \cu(S)$.
\item[(d)] if $S$ is homogeneous (i.e., invariant under the the homothecy action), then
$\cu(S)$ is again homogeneous.
\end{enumerate}}
\end{rems}

\subsection{The $\mu$ function}
We shall now recall some basic facts from Geometric Invariant Theory \cite{M}. Critical to our
understanding stability, semi-stability, and unstability is the notion of the $\mu$-function. We
give two definitions : \Dref{def:extr-mu} as well as a more intrinsic one (applicable in more
general situations) in \Dref{def:intr-mu}.

Let $(X,\,L,\,\bp(V))$ be a $G$-triple such that the standard linear system on $\bp(V)$ arising
from $\co(1)$, when restricted to $X$, is a complete linear system (i.e., the {\em trace} on $X$---of
the complete linear system on $\bp(V)$---is a complete linear system on $X$). This is equivalent
to saying that $H^0(X,\,L)=V^*$.  Let $\lambda\colon \bg_m\to G$ be a 1-PS. Recall that
the action of $\lambda$ on $X$, $\bp(V)$, $\wh{X}$, $V$ is, by definition, the action of $\bg_m$
on these spaces induced by $\lambda$ and the $G$-action on them (see (6) of \Ssref{ss:conventions}).

Now, action of $\lambda$ on $V$ can be diagonalised, whence we can find a basis $\{e_i\}$ of
$V$ such that the one dimensional subspaces $ke_i$ are $\bg_m$-stable under the action
of $\lambda$, and hence give rise to characters, one for each subscript $i$:
\[\chi_{{}_i}\colon \bg_m\to \bg_m,\qquad t\mapsto t^{r_i} \]
with $r_i=r_i^L(\lambda)$ an integer. Now, suppose $x\in X$ and suppose 
\[\wh{x}=\sum x_ie_i\] 
is a point of $\wh{X}\smallsetminus\{0\}$ lying over $x$. 

\begin{defi}($\mu$-function) \label{def:extr-mu} With notatations as above, the $\mu$-value
of $x$ with respect to $L$ and $\lambda$ is
\[\mu^L(x,\,\lambda)\set \max_{j, x_j\neq 0} \{-r_j\} = -\min_{j, x_j\neq 0} \{r_j\}.\]
\end{defi} 
Clearly, the definition of $\mu^L(x,\,\lambda)$ does not depend on the choice of the point
$\wh{x}$ lying above $x$. If the role of $L$ is understood from the context, we will simply
write $\mu(x,\,\lambda)$ for $\mu^L(x,\,\lambda)$.

Note that we have:
\stepcounter{thm}
\begin{equation*}\label{eq:hilb-mum}\tag{\thethm}
\begin{split}\lim_{t\to 0} x\circ \lambda \,\, {\text{exists}} &\Longleftrightarrow \mu^L(x,\lambda) \ge 0 \\
\lim_{t\to 0} x\circ \lambda= 0  &\Longleftrightarrow \mu^L(x,\,\lambda) >0
\end{split}
\end{equation*}

We can define the function $\mu$ in a more intrinsic manner as follows: Let $G$ act on a {\em complete}
$k$-scheme $X$, and suppose the action lifts linearly to an action on a line bundle $L$ on $X$.
Note that we are not assuming that $X$ is projective, or that $L$ is ample. Let $x\in X$, $\lambda$
a 1-PS, and denote by $\Psi_x\colon \bg_m\to X$ the orbit morphism defined by 
$t\mapsto x\circ \lambda(t)$. Since $X$ is complete, the morphism $\Psi_x$ extends to a morphism
of $\bp^1$ into $X$, which we again denote by $\Psi_x$. Let $x_\circ=\Psi_x(0)$. Then $x_\circ$
is invariant under the action of $\lambda$. Let $L_\circ$ be the fibre of $L\to X$ over $x_\circ$.
Then the operation of $\bg_m$ on $L_\circ$ is defined by a character $\chi\colon \bg_m\to\bg_m$
defined by $t\mapsto t^r$.

\begin{defi}(Intrinsic definition of $\mu$)\label{def:intr-mu} Let $X$, $L$, $\lambda$, $x$, and $r$ be
as above. Then the $\mu$-value of $x$ with respect to $L$ and $\lambda$ is\[\mu^L(x,\,\lambda) (=\mu(x,\,\lambda)) =-r.\]
\end{defi}

\subsection{The Hilbert-Mumford criterion for $S$-unstablility} 
Let $(X,\,L,\,\bp(V))$ be a
$G$-triple. Suppose $S\subset V$ is a $G$-invariant closed subscheme. We claim that we can find
a finite dimensional rational $G$-module $W$ and a $G$-equivariant morphism
\stepcounter{thm}
\begin{equation*}\label{eq:f-V-W}\tag{\thethm}
f\colon V\to W
\end{equation*}
(now thinking of $V$ and $W$ as schemes) such that set-theoretically, or even scheme-theoretically,
$f^{-1}(0)=S$. Indeed, let $F_1,\dots, F_l$ be generators of the ideal $I_S\subset k[V]$ of $S$,
where $k[V]$ is the ring of functions $\Gamma(V,\,\co_V)$ on the affine scheme $V$. In other
words, $k[V]=S(V^*)$, the symmetric algebra on the dual linear space $V^*$ of $V$. By considering---if 
necessary---the $G$-span of the $k$-linear space spanned by $F_1,\,\dots,\,F_l$, we may assume that
the linear space spanned by the $F_i$ is $G$-stable and that $F_1,\,\dots,\,F_l$ are linearly 
independent over $k$. Let $W^*$ be this $l$-dimensional $G$-stable subspace of $k[V]$, 
and let $W$ be its
dual. Now, $S(W^*)$ maps in a {\em surjective} and {\em $G$-equivariant} manner
onto the subalgebra $k[F_1,\,\dots,\,F_l]$ of $k[V]$, and we
have a sequence of $G$-equivariant maps of $k$-algebras:
\[ k[W]= S(W^*) \twoheadrightarrow k[F_1,\dots,F_l] \subset k[V] \]
The resulting algebraic map of varieties $f\colon V\to W$ meets our requirements.

We would like to understand the locus of points in $\cu(S)\smallsetminus S$. Suppose 
$v\in \cu(S)\smallsetminus{S}$. We can find a discrete valuation ring (d.v.r) $A$ with residue field
$k$ and quotient field $K$, as well as $K$ valued point of $G$, $\theta\in G(K)$, such that
$v\circ\theta\in V(K)$ is really an $A$-valued point, i.e., $v\circ\theta\colon A\to V$, and the closed
point of $A$ maps to a $k$-valued point $v_\circ \in S$.  If $t$ denotes the 1-parameter represented
by $\spec{A}$, we write, as a shorthand,
\[\lim_{t\to 0}v\circ \theta(t) = v_\circ.\]
In this situation, the Hilbert-Mumford Theorem, as generalized by Kempf \cite{GK}  for the case
of $S$-instability, states that there exists in fact, a 1-PS $\lambda\colon\bg_m\to G$ such that
(in an obvious notation)
\[\lim_{t\to 0} v\circ\lambda(t) = v_\circ.\]
One can show:

\begin{thm}\label{thm:hilb-mum}{\em{(Hilbert-Mumford)}} Let $(X,\,L,\,\bp(V))$ be a
$G$-triple and $S\subset\wh{X}$ a closed $G$-invariant subset. Let $f\colon V\to W$
be as in \eqref{eq:f-V-W}, i.e., $f$ is a map of varieties such that $f^{-1}(0)=S$. Fix a Borel
subgroup $B\subset G$ and a maximal torus $T\subset B$ of $G$, and let $\Gamma(T)$
be the co-root lattice of $T$ and $C(B)$ the Weyl chamber in 
$\Gamma(T)\otimes_{\mathbb{Z}}{\mathbb{R}}$ corresponding to $B$, and $\ol{C(B)}$ its
closure in the Euclidean space $\Gamma(T)\otimes {\mathbb{R}}$.
\begin{enumerate}
\item[(a)] The following are equivalent:
\begin{itemize}
\item[(i)] $v\in \cu(S)\smallsetminus S$;
\item[(ii)] $f(v)\neq 0$ and there exists a one-parameter subgroup $\lambda$ of $G$ such that 
$\lim_{t\to 0} v\circ\lambda(t)$
exists and $\lim_{t\to 0} f(v)\circ\lambda(t)=0$;
\item[(iii)] $f(v)\neq 0$ and there exists $g\in G$ and a 1-PS $\lambda$ of $T$ in the chamber
$\ol{C(B)}$ such that $\lim_{t\to 0}(v\circ g)\lambda(t)$ exists and $\lim_{t\to 0} f(v\circ g)\lambda(t)=0$.
\end{itemize}
\item[(b)] The following are equivalent:
\begin{itemize}
\item[(iv)] $v\in \cu(S)$;
\item[(v)] There exists a one-parameter subgroup $\lambda$ of $G$ such that $v\in\cu(S,\,\lambda)$;
\item[(vi)] There exists $g\in G$ and a one-parameter subgroup
$\lambda$ of $T$ in the chamber
$\ol{C(B)}$ such that $v\circ g\in\cu(S,\,\lambda)$.
\end{itemize}
\end{enumerate}
\end{thm}

\proof In view of the comments made before the statement of the theorem, evidently, (i) and (ii)
are equivalent as are (iv) and (v). For the rest we only need to recall the well-known fact that
every 1-PS $\lambda$ of $G$ is conjugate to a 1-PS in $T$, in fact to and integral point in $\ol{C(B))}$.
\qed

\begin{prop}\label{prop:finite1PS}
Let the notations be as in the hypotheses of
\Tref{thm:hilb-mum}. There exist a finite number of 1-PS $\lambda_1,\dots, \lambda_N$
of $T$ in the chamber $\ol{C(B)}$ such that $\cu(S)=\bigcup_{1\le i\le N} \cu(S,\,\lambda_i)\cdot G$.
\end{prop}

\proof 
According to \Lref{lem:cones} in \Sref{s:s=ss} below, there are a finite number of closed convex cones
$C_\alpha$ in $\ol{C(B)}$---each $C_\alpha$ an intersection of a finite number closed 
half-spaces---such that for every $v\in V$ and $w\in W$, $\mu(v,\_)$ and $\mu(w,\_)$ are linear
in each $C_\alpha$. This is seen setting $d=2$, $X_1=\bp(V)$, and $X_2=\bp(W)$ in 
{\em loc.cit}. Let ${\mathfrak{S}}_\alpha$ be a finite set of generators over $\br^+$ for the cone
$C_\alpha$. We can choose ${\mathfrak{S}}_\alpha$ with integral coefficients, i.e., 
${\mathfrak{S}}_\alpha\subset \Gamma(T)$.
Let ${\mathfrak{S}}=\cup_\alpha {\mathfrak{S}}_\alpha$.
Since $\mu(v,\_)$ and $\mu(f(v),\_)$ are linear in each $C_\alpha$, we conclude that:\\
{\em  $\mu(v,\,\lambda)\ge 0$ and $\mu(f(v),\,\lambda)>0$  for every $\lambda\in \ol{C(B)}\setminus\{0\}$
if and only if $\mu(v,\,\lambda) \ge 0$ and $\mu(f(v),\,\lambda)>0$ for every $\lambda\in {\mathfrak{S}}$.}

Let ${\mathfrak{S}}=\{\lambda_1,\dots, \lambda_N\}$.  Then by \Tref{thm:hilb-mum}(a)---especially
(i)$\Leftrightarrow$(iii)---and the relations  \eqref{eq:hilb-mum}, 
it is evident that $\cu(S)$ is the union of $\cu(S,\,\lambda_i)\cdot G$ for
$i=1,\dots, N$.
 
 \qed

\subsection{The $S$-unstable locus $\cu(S)$ is closed}
Fix a $G$-triple $(X,\,L,\,\bp(V))$.
Recall the following properties of the function $\mu=\mu^L$ (see
\cite[Prop.\,3.1]{S2}): Let $\lambda\colon\bg_m\to T$ be a 1-PS and let $P(\lambda)(\supset B)$
be the parabolic subgroup of $G$ defined by
\[P(\lambda) \set \{g\in G\,\vert\, \lim_{t\to 0}\lambda(t)^{-1}g\lambda(t) \,\, {\text{exists in $G$}}\}.\]
Then, for every $v\in V$, we have:
\begin{enumerate}
\item[(a)](Change of coordinates) $\mu(v,\,\lambda) = \mu(v\circ g,\,g\lambda g^{-1})$, 
$g\in G$.
\item[(b)] $\mu(v,\,\lambda) = \mu (v,\,g\lambda g^{-1})$, $g\in P(\lambda)$.
\item[(c)] $\mu(v,\,\lambda)= \mu(v\circ g,\,\lambda)$, $g\in P(\lambda)$.
\item[(d)] $\mu(v,\,\lambda) = \mu(v\circ g,\, \lambda)$, $\lambda\in \ol{C(B)}$ and $g\in B$.
\end{enumerate}
Note that (c) follows from (a) and (b), and (d) is a special case of (c), since $\lambda\in\ol{C(B)}$
implies that $B\subset P(\lambda)$.

\begin{prop}\label{prop:U(S)-closed} (See \cite[p.\,524,\,Theorem\,3.1]{S2}) 
Let $S$ be a $G$-invariant closed subset of $\wh{X}$. Then the set $\cu(S,\,\wh{X})$ 
is a closed $G$-invariant subset of $\wh{X}$.
\end{prop}

\proof First note that for any 1-PS $\lambda$ of $G$, the set $\cu(S,\,\lambda)$ is closed. Indeed, since
$\bg_m$ is linearly reductive, the GIT quotient $V/\negthickspace/\lambda$ for the action
of $\lambda$ on $V$ exists as a variety and if $j\colon V\to V/\negthickspace/\lambda$ is the
canonical quotient morphism, then it is easily seen that $j(S)$ is closed in $V/\negthickspace/\lambda$.
Now, $\cu(S,\,\lambda)=j^{-1}(j(S))$, thus proving that $\cu(S,\,\lambda)$ is closed in $V$.
Let $\lambda_1,\dots,\lambda_N$ be the finite number of 1-PS of $T$ in $\ol{C(B)}$ guaranteed
by \Pref{prop:finite1PS} above, with the property that 
$\cu(S)=\bigcup_{1\le i\le N} \cu(S,\,\lambda_i)\cdot G$. Each $\cu(S,\,\lambda_i)$ is $B$-stable
and hence so is
\[U=\bigcup_i \cu(S,\,\lambda_i).\]
Now for any scheme $Z$ on which $B$ acts, one defines $Z\times^{B}G$ as the set of equivalence
classes of the equivalence relation on $Z\times G$ given by
$(z,\,g)\sim (v\circ b,\,b^{-1}g)$, $z\in Z$, $b\in B$, $g\in G$.\footnote{The equivalence relationship
is by a free group action of $B$, and hence the quotient $Z\times^BG$ exists as a scheme, and the
quotient map $Z\times_k B\to Z\times^BG$ is a principal $B$-bundle.}
The natural map $\pi_Z\colon Z\times^BG\to B\backslash G$ is a fibre bundle with fibre $Z$, and 
structure group $B$. The associated principal fibre bundle is the canonical quotient
$G\to B\backslash{G}$.
The subset $U\cdot G=\cu(S)$ of $V$ is the image
of the map $U\times G \xrightarrow{f} V$ given by $(u,\,g) \mapsto u\circ g$.
Let $p\colon U\times G\to U\times^BG$ be the natural quotient map. Then clearly the map
$U\times G \xrightarrow{f} V$ factors through $p$, giving us a map $f'\colon U\times^BG\to V$
such that $f=f'\circ p$, i.e., a we have a commutative diagram:
\[{\xymatrix{
U\times G \ar[dr]_{f} \ar@{->>}[r]^-{p} & U\times^BG\ar[d]^{f'} \\
& V
}}
\]

Next note that the map $V\times^BG \to V$ given by $(v,g)\mapsto v\circ g$ is isomorphic to the
trivial bundle $V\times B\backslash{G}\to V$, the explicit isomorphism 
$V\times^BG \iso V\times B\backslash{G}$ being $(v,\,g)\mapsto (v\circ g,\,Bg)$. If 
$U\times^BG \hookrightarrow V\times^BG$ is the closed immersion induced by the closed
immersion $U\hookrightarrow V$, we get a commutative diagram:
\[{\xymatrix{
U\times G \ar[dr]_{f} \ar@{->>}[r]^-{p} & U\times^BG{\phantom{X}} \ar[d]^{f'} \ar@{^{(}->}[r] & V\times^BG
\ar@{=}[d]\\
&  V & V\times B\backslash{G} \ar[l]^{\text{projection\phantom{XX}}}
}}
\]
Now $\cu(S)=U\cdot G= f(U\times_k G) = f'(U\times^BG)$---the last equality from the surjectivity of $p$. 
The map $f'$ is clearly proper, being the composite of a closed immersion 
followed by the proper map $V\times_k B\backslash{G}\to V$. Thus $\cu(S)$ is a closed subset 
of $V$.
\qed

\begin{cor}\label{cor:Xss-open}
The set $\wh{X}^{ss}$ (resp. $X^{ss}$) is open and $G$-invariant in
$\wh{X}$ (resp. $X$).  A point $x \in X^{ss}$ is stable i.e. 
$x \in X^s$ if and only if the orbit $O(x) = x G$ is closed in
$X^{ss}$ and $\dim O(x) = \dim G$ i.e. the orbit morphism $G \lr
X^{ss}$, $g \longmapsto x \circ g$, is proper.
\end{cor}
\proof Taking $S = (0)$ the first assertion follows. 
Recall that for $x \in \wh{X}^{ss}$, we have $\ol{O(x)}$ (closure in $\wh{X}$) 
$\subset \wh{X}^{ss}$.  The second assertion of the corollary reduces
to proving that for $x \in X^{ss}$, $O(x)$ closed in
$X^{ss} \Llr O(\wh{x})$ closed in $\wh{X}^{ss}$. But this is 
an immediate consequence of the following claim:

{\em  Let $x$ be a point in $X^{ss}$ and $\hat{x}$ a point in $\wh{X}^{ss}$ 
lying over~$x$. Then for any $y \in \ol{O(x)} \cap X^{ss}$, there exists a 
point $\hat{y} \in \ol{O(\hat{x})}$ lying over~$y$}. 

The claim is proved as follows:
Set $A \set k[[t]]$ and $K \set k(\!(t)\!)$. 
Let $C$ be an irreducible curve 
in~$X$ joining $x$ and $y$ and $C_1$ an irreducible curve 
in $\wh{X} \times G$ through $(\hat{x}, e)$ mapping dominantly 
to $C$ via $\wh{X} \times G \xrightarrow{\;\sigma\;} \wh{X} -\!\!\to X$. 
Projecting $C_1$ to $G$ results in a $\spec(K)$-valued point 
$g$ (or $g \in G(K)$), 
such that $xg \in X(A)$. As in the proof of \cite[p.\,520,\,Theorem\,2.2]{S2}, we may assume that $\wh{X} = V$ 
(i.e., affine $n$-space),
and that $g = U\lambda$ where $U \in G(A)$ 
satisfies $\lim_{t \to 0}U =$ identity matrix and 
$\lambda \in T(K)$ ($T$ a maximal
torus in $G$) is diagonal of the form $(t^{r_1}, \ldots,t^{r_n})$.  

It suffices to show that $\lim_{t \to 0}\hat{x}g$ 
exists, as then we may choose $\hat{y}$ to be this limit point 
(which is necessarily nonzero owing to semistability of~$x$). 
If the limit does not exist, then there is a unique integer
$s < 0$ such that 
$\hat{z} \set \lim_{t \to 0}t^{-s}\hat{x}g$ exists 
and is nonzero. Let $\hat{z}_i$ denote the $i$-th coordinate
of~$\hat{z}$. Since $\{ i \;|\: r_i \ge 0 \} \subset 
\{ i \;|\: \hat{z}_i = 0 \}$, we see that $\hat{z}$ is unstable via the 
action of~$\lambda^{-1}$. Since $\hat{z}$ lies over~$y$, this 
contradicts semi-stability of~$y$. \qed

\begin{defi}\label{def:S-unstable2}
Let $S$ be a closed $G$-stable subset of $X^{ss}$.
A point $x \in X^{ss}$ is said to be $S$-{\em unstable} if the closure
 $\ol{xG}$ in $X^{ss}$ of the orbit $xG$ meets $S$.  We denote by $\cu(S)=\cu(X,S)$ the 
$G$-invariant subset of $S$-unstable points in $X^{ss}$.
\end{defi}

\begin{cor}\label{cor:S-closed}
With above notations, $\cu (S)$ is a closed $G$-invariant subset of $X^{ss}$.  
\end{cor}
\proof Let $\wh{S}$ be the ``cone over $S$'' i.e. the homogeneous
closed subset of $\wh{X}$ defined by $S$.  In view of the above lemma,
the corollary is an
easy consequence of the proposition applied to $\wh{S}$.
\qed

\begin{defi}\label{def:ss-equiv} Let $v_1, v_2 \in V$ (or $\wh{X}$). We sat write $v_1\sim v_2$ if 
$\ol{O(v_1)} \cap \ol{O(v_2)} \neq \es$
(i.e., the orbit closures of $v_1$ and $v_2$ intersect). Similarly for $x_1, x_2\in X^{ss}$, we write 
$x_1\sim x_2$ if $\ol{O(x_1)} \cap \ol{O(x_2)} \neq \es$. One checks easily that $x_1\sim x_2$ 
if and only if there exist $v_1, v_2 \in V\setminus\{0\}$ lying over $x_1$ and $x_2$ respectively 
such that $v_1\sim v_2$. We call these relations {\em semi-stable (or orbit closure) equivalence 
relations}, on account of the following:
\end{defi}

\begin{cor}\label{cor:ss-equiv}
The relations defined in \Dref{def:ss-equiv} are equivalence relations.
\end{cor}
\proof To prove the equivalence relation property, we see
that it suffices to prove the transitivity property.  We first
observe that this property is equivalent to showing that there
is a unique closed orbit in $\ol{O(x)}$ (closure of $O(x)$ for
$x \in V$.  It suffices to consider this case). To see this, suppose 
that $\ol{O(x)}$ contains two distinct closed orbits $\ol{O(x_1)}$
and $\ol{O(x_2)}$.  Then we have $x_1 \sim x$ and $x \sim x_2$
but $x_1$ is not equivalent to $x_2$.  On the other hand if
every orbit closure has a unique closed orbit then we see that
$x_1 \sim x_2$ if and only if the unique closed orbits in 
$\ol{O(x_1)}$ and $\ol{O(x_2)}$ coincide.  Then the transitivity
property is immediate.  Let us now prove that every $\ol{O(x)}$ has 
a unique closed orbit.  Again suppose that there are two
distinct closed orbits $\ol{O(x_1)}$ and $\ol{O(x_2)}$ in
$\ol{O(x)}$.  Then we see that $x$ is 
$\ol{O(x_1)}$-unstable, which implies that $\ol{O(x)} \subset
\cu (\ol{O(x_1)})$ (by the proposition).  But $x_2 \not \in
\cu (\ol{O(x_1)})$, which leads to a contradiction.  This
proves that the relation is transitive. \qed

\begin{cor}\label{cor:polystable}
Let $\wh{X}^e$ and $\wh{X}^{ps,e}$ denote the subsets of 
$\wh{X}$ defined by:
\begin{align*}
& \wh{X}^e = \{ x \in \wh{X} \mid \dim O(x) \leq e \}\\
& \wh{X}^{ps,e}  = \{ x \in \wh{X} \mid \dim O(x) =e
{\text{{\em{ and $x$ is  polystable}}}} \} .
\end{align*}
Similarly, let $X^{ss,e}$ and $X^{ps,e}$ denote the subsets of $X^{ss}$ defined
by
\begin{align*}
& X^{ss,e} = \{ x \in X^{ss} \mid \dim O(x) \leq e \} \\
& X^{ps,e} = \{ x \in X^{ss} \mid x ~~{\text{\em{is polystable
and }}} \dim O(x) = e \}
\end{align*}
Then
\begin{enumerate}
\item $\wh{X}^e$ is closed and $G$-invariant in $\wh{X}$.  
\item $\wh{X}^{ps,e}$ is open and $G$-invariant in $\wh{X}^e$.  
\item $X^{ss,e}$ is closed and $G$-invariant in $X^{ss}$
\item $X^{ps,e}$ is open and $G$-invariant in $X^{ss,e}$.
\end{enumerate}
\end{cor}

\noindent{\em Observations:} Specializing to  $e=\dim{G}$ in (2), we have
$\wh{X}^s$ is open and $G$-invariant in $\wh{X}$. Note also that, with
$I_x$ the isotropy group at $x$, we have description:
\[{\wh{X}}^e = \{x\in X \mid
\dim I_x 
\geq \dim G-e\}. \]
 The sets $\wh{X}^s$ and $\wh{X}^{ps,e}$ could
be empty.

\proof
That $\wh{X}^e$ is closed is seen easily.  It is of course $G$-invariant.
$$ \wh{X}^e \setminus \wh{X}^{ps,e} = \cu (\wh{X}^{e-1})
~~ ({\rm in}~ \wh{X}^e)$$
and then the corollary follows from the proposition.  The proofs
of the other assertions are similar. \qed

\subsection{Base Change for $X^{ss}$}\label{ss:basechange}
Throughout this subsection we fix $(X,\,L,\,\bp(V))$, which is
a $G$-triple over $k$.
Let $k'$ be an extension of $k$ and say $k'$ is also algebraically closed.  
Let $X_{k'}$ be the base change of $X$ by $\spec k'
\lr \spec k$ (similarly $\wh{X}_{k'}$ etc.).  Then we claim that
$X_{k'}^{ss}$ is the base change of $X^{ss}$ by $\spec k' \lr 
\spec k$.  When the group $G$ is a torus, say a 1-PS subgroup,
this assertion is seen easily.  Then the general case follows 
from \Pref{prop:finite1PS}.  The point of this proof can also be 
stated as follows.  A point $x$ of $X_{k'}$ i.e. an element of
$X(k')$  ($k'$ valued points of $X$) is semi-stable  if and only
if $\exists$ a finite number of 1-PS $\lambda_1,\ldots ,\lambda_N$
{\it defined over $k$} in the chamber $C(B)$ and a point $g \in G(k')$
such that $x\circ g'$ is semi-stable with respect to $\lambda_i$.
This reduces the assertion to the case of a 1-PS.

Suppose that the base field $k$ is not algebraically closed, but
the schemes and actions are defined over $k$.  Let us define 
a geometric point $x$ of $X$ (i.e. $x\in X(\Omega )$, where
$\Omega$ is algebraically closed) to be semi-stable if $x$
is a semi-stable point of the base change $X_{\Omega}$.
Then the above argument, in fact, shows that there is
a $G$-stable open subscheme $X^{ss}$ of $X$ such that its 
geometric points are precisely the semi-stable points.
This should of course be stated for more general base schemes
$S$.  One should define semi-stability only for geometric points
and have an open subscheme $X^{ss}$ whose points are
precisely the semi-stable points as in \cite{S3}.

\begin{rem}\label{rem:dvrbase}  
{\em Let $A$ denote the discrete valuation ring $k [[t]]$
with quotient field $K$.  Let $x$ be a $K$-valued point 
of $X_K$, with $X_K$ denoting base change,
as above.  Let $H$ be an algebraic subgroup of $G$.
Then $H_K$ operates on $X_K$.   Let $Z_K$ be the closure in
$X_K$ of the $H_K$ orbit through $x$.  Let $Z_A$ be the closed 
subscheme of $X_A$, flat over $A$, determined by $Z_K$.
Then we see that the group scheme $H_A = H \times_k \spec A$
operates on $Z_A$ and the generic fibre $Z_K$ of $Z_A \lr
\spec A$ contains an open orbit under the action of $H_K$
(namely the $H_K$ orbit through $x$).}
\end{rem}

\begin{lem}
\cite[pp.\,528---529, Rmk.\,4.9]{S2} Let $x \in X^{ss}$ and $Z$ the $G$-orbit
through $x$ with its reduced structure.  Then there exists an integer $n > 0$
such that the restriction of $L^n$ to $Z$ is trivial.  
\end{lem}

\proof
Let $H$ be the isotropy subgroup of $G$ at the point $x$.  Then $Z$
is the homogeneous space $H \setminus G$ and the restriction of
the $G$-line bundle to $Z$ is defined by a character $\chi : H \lr \bg_n$.  
It suffices to prove that $\exists n >0$ such that $\chi^n$ is the
trivial character.  Suppose that this is not the case.  Then we 
see easily that we can find a 1-PS $\lambda$ of $H$ such that $(\chi 
\circ \lambda ) : \bg_m \lr \bg_m$ is surjective.  We see that 
$(\chi \circ \lambda )$ is the character defining the action of 
$\lambda$ on the fibre $L$ at $x$ and the surjectivity of
$\chi \circ \lambda$ implies that $\mu (x,\lambda ) \neq 0$
(see  \Dref{def:intr-mu} above).  Then either $\mu (x ,\lambda )$ or
$\mu (x,\lambda )$ is $<0$.  This contradicts the fact that
$x \in X^{ss}$. \qed

\subsection{Separatedness and properness properties of $(X^{ss}/\negmedspace/G)_{\text{top}}$}
\label{ss:G-sep}
In this section we investigate some basic topological properties of the semi-stable
equivalence relation defined in \Dref{def:ss-equiv} (see also \Cref{cor:ss-equiv}). A crucial
tool is \Lref{lem:bijective1} concerning proper bijective maps of algebraic schemes over fields
of positive characteristic.

For the remainder of this section we fix a quotient data
\stepcounter{thm}
\begin{equation*}\label{eq:q-data}\tag{\thethm}
\De=(X,\,L,\,\bp(V),\,X^{ss}\xrightarrow{\alpha} Y).
\end{equation*}
We may assume, without loss of generality, if the occasion demands, that
$V$ is the vector space dual of $H^0(X,\,L)$, i.e., that the linear system on $X$ induced
by the trace of the tautological linear system on $\bp(V)$ is complete.

\begin{defi}\label{def:top-quot}
The {\em product topology} on 
$Y\times Y$ is 
the quotient topology induced by the action of $G\times_kG$ on $X^{ss}\times_k X^{ss}$.
In other words, as toplogical spaces 
\[Y\times Y=
\bigl(X^{ss}\times_kX^{ss}/\negmedspace/(G\times_kG)\bigr)_{\text{top}}.\]
\end{defi}

We will be showing that $Y$ is ``separated", and ``complete"
in \Pref{prop:G-sep} and \Cref{cor:proper} respectively. 

\begin{prop}\label{prop:G-sep} 
The graphs of the semi-stable equivalence relations are closed. In fact, if $\Delta_X$ is the diagonal
in $\wh{X} \times \wh{X}$ (or $X^{ss} \times X^{ss}$) the
graph $\Gamma$ of the equivalence relation is given by
\[\Gamma = \cu (\ol{\Delta_X \cdot (G \times G)}).\]
In particular $Y$ is ``separated", i.e., the diagonal map $Y\to Y\times Y$ is a closed embedding,
with $Y\times Y$ having the product topology defined above.
\end{prop}

\proof Let us work with the case $\wh{X} = V$.  The proofs
in the other cases are immediate consequences.  Let us take a
d.v.r. $A$ and its quotient $K$ as in \Rref{rem:dvrbase}.  Note that
$\ol{\Delta_X \cdot (G \times G)}$ consists of points $(x_0,y_0)$
of the form:
{\renewcommand{\labelenumi}{{\rm (\roman{enumi})}}
\begin{enumerate}
\item $x$ and $y$ are $A$-valued points of $V$.
\item $\ds\lim_{t \ra 0} x(t) (=x)=x_0$, $\ds\lim_{t\ra 0} y(t) (=y) = y_0$.
\item There exists $g \in G(K)$ such that $y = x \cdot g$.
\item $x_0 \sim y_0$.
\end{enumerate}}
By Iwahori-Matsumsto \cite{I-M} (possibly by going to the integral closure
in a finite extension of $A$), we have $g = P \lambda Q$, where 
$P,Q$ are in $G(A)$ and $\lambda : \bg_m \lr G$ is a 1-PS defined
by a $K$-valued point of $G$.  Set $x' = x\cdot P$ and $y' = x \circ P 
\lambda = x' \circ \lambda$.  We see that $x_0' = \ds \lim_{t \ra 0} x'(t)$,
$y_0' = \ds \lim_{t \ra 0} y(t)$ are in the $G$-orbits of $x_0$ and
$y_0$ respectively (for $y_0 = y_0' \cdot Q_0$, $Q_0 = \ds \lim_{t\ra 0}
Q$) and thus to prove (iv) we can assume without loss of generality that
$y = x \circ \lambda$, so that we are reduced to the case $G$ is a torus
(of dimension one).  In this case we know that $x_0$ is semi-stably
equivalent to $y_0$ for the action of $\lambda$, hence \'a 
fortiori $x_0 \sim y_0$.

Thus we see that $\ol{\Delta_X \cdot (G \times G)} \subset \Gamma$.
Take a point $(x_0,y_0) \in \Gamma$.  Then $\ol{O(x_0)}$ and 
$\ol{O(y_0)}$ have a unique common closed orbit $O(z)$ and we can suppose
that $\exists ~ g= g_t$, $h=h_t$ with $g,h \in G(K)$ and $x_0 \circ
g$, $y_0 \circ h$ are $A$-valued points of $V$ with 
$z = \ds\lim_{t \ra 0} x_0 \circ g = \lim_{t \ra 0} y_0 \circ h$.
This implies that $(x_0,y_0) \in \cu (\ol{\Delta_X \cdot (G \times G)})$
so that $\Gamma \subset \cu (\ol{\Delta_X \cdot (G \times G)})$.
On the other hand let $(z,w) \in  \cu (\ol{\Delta_X \cdot (G \times G)})$.
This means that $\exists g = g_t$, $h=h_t$, such that if 
$z_0 = \ds\lim_{t \ra 0} z \circ g_t$, and $w_0 = \ds\lim_{t\ra 0}
w \circ h_t$ (in the sense as above) then $(z_0,w_0) \in$
($\ol{\Delta_X \cdot (G \times G)}$).  Then as shown above 
$z_0 \sim w_0$.  On the other hand we have obviously $z \sim z_0$,
$w \sim w_0$.  Hence by the transitivity of the semi-stable 
equivalence relation, we see that $z \sim w$ i.e. $\cu 
(\ol{\Delta_X \cdot (G \times G)}) \subset \Gamma$.  Thus we have
$\Gamma = \cu (\ol{\Delta_X \cdot (G \times G)})$. 
\qed

\smallskip

The next Proposition, viz.\, \Pref{prop:L-trivial} 
says that if the quotient $Y$ is
a single point, then $L$ is essentially trivial. More precisely, a power of $L$ is trivial. Morally
then, $L$ could be regarded as the pull back of a line bundle on 
$Y$, which, consisting of exactly one point, carries
up to isomorphism, only one line bundle. However, $Y$ not
(as yet) having a scheme structure, these observations right now come under the heading of
``gathering evidence" that $Y$ has a scheme structure.
The Proposition plays an essential role in showing that $Y$
has a natural scheme structure, making it the GIT quotient $Y$.

\begin{prop}\label{prop:L-trivial} Suppose the quotient data $\De$ is reduced
and $X^{ss}$ consists of {\it one} 
semi-stable equivalence  class, i.e., any two points of $X^{ss}$ are semi-stably
equivalent.   Suppose further that the closure of $X^{ss}$ in $\bp(V)$ is $X$.
Let $\wh{X}$ stand for the cone over $X$ with its reduced structure. Then a suitable power
of $L$ has a section $s$ over $X$ which is $G$-invariant and 
non-vanishing for every $x \in X^{ss}$.  In fact we can find $s$
which comes from a regular $G$-invariant function on $\wh{X}$.
\end{prop}
\proof We have a unique closed orbit in $X^{ss}$.  We
denote the closure of this in $X$ by $X_1$.  Then we observe
that 
$$ \cu (\wh{X}_1) = \wh{X} \leqno {\rm (i)}$$
i.e. the closure of the $G$-orbit through any point of
$\wh{X}$ meets $\wh{X}_1$.  Let us now suppose that
the proposition holds for $X_1$ (the proof will be given later after the proof of
\Lref{lem:bijective2}).  We first show that the proposition is true under our supposition.

We have then a $G$-invariant function $f$ on $\wh{X}_1$ such that
$f(x) \neq 0$ $\forall x \in \wh{X}_1^{ss}$ (of course 
$f(x) =0$ $\forall$ $x \in \wh{X} \setminus \wh{X}^{ss}$)
which is ``homogeneous'' i.e. $f(tx) = t^m f(x)$, where
multiplication by $t$ denotes the homothecy action and
$L$ is the $m$th power of the tautological ample line
bundle on $\bp (V)$.  Let us extend $f$ to a set theoretic
function $F$ ($F: \wh{X} \lr \ba^1)$ as follows.  Given $x \in
\wh{X}$, the orbit closure $\ol{O(x)}$ meets $\wh{X}_1$ (by (i)),
say at a point $y$.  Let us set $F(x) = F(y)$.  We claim that 
it is well-defined.  For, if $y_1,y_2$ are two such points, 
$y_1 \sim y_2$.  If $x \not \in \wh{X}^{ss}$, $y_1,y_2 \in
X_1 \setminus X_1^{ss}$ so that $F(x) = f(y_1) = f(y_2) =0$.
If $x \in X^{ss}$, $y_1,y_2$ are semi-stably equivalent, which
implies that $y_2 \in \ol{O(y_1)}$.  Since $f$ is $G$-invariant,
if follows that $f(y_1) = f(y_2)$.  Thus the set theoretic 
function $F$ is well-defined,
 $G$-invariant and is an extension of $f$.
We observe that $F$ is also ``homogeneous of the same degree''
as that of $f$.

Note that if the complete reducibility property holds for $G$
(e.g. $G$ a torus), then one knows that $f$ can be extended
to a $G$-invariant function on $\wh{X}$ and in our case
$F$ is uniquely determined.  Thus in this case $F : \wh{X}
\lr \ba ^1$ is a morphism.  We will use this fact below.

As we have seen in \Pref{prop:finite1PS} there exist 
1-PS $\lambda_i$, $1 \leq i \leq N$, of $T$ in $\ol{C(B)}$
such that
$$ \wh{X} = \bigcup_{1 \leq i \leq N} \cu (\wh{X}_1,\lambda_i)
\cdot G.  \leqno {\rm (ii)}$$
Set $\cu_i = \cu (\wh{X}_1,\lambda_i)$ i.e. $\cu_i$ is the
set of $\wh{X}_1$ unstable points for the 1-PS $\lambda_i$, 
then there is a unique $\lambda_i$-invariant regular function
$\theta_i$ on $\cu_i$ which extends $f$.  Now $F$ coincides
with $\theta_i$ on $\cu_i$ so that we see that the restriction
to each $\cu_i$ is regular.  Now $\theta_i-\theta_j$ vanishes set
theoretically on the scheme theoretic intersection $\cu_i
\cap \cu_j$.  Hence
\[ 
\left \{ \parbox{4.2in}{$(\theta_i-\theta_j)^q 
= (\theta_i^q-\theta_j^q)$
vanishes on the scheme intersection $\cu_i \cap \cu_j$, for
$q = p^r$, $r \gg 0$.} \right . \leqno {\rm (iii)}'
\]
Let then $\cu$ denote the (scheme theoretic) union $\bigcup \cu_i$.  
Then $\cu$ is reduced, $\cu_i$ being reduced.  We see for
$q=p^r$, $r \gg 0$, the restriction of $F^q$ to $\cu$ is a regular
function.  Thus without loss of generality, we can assume that 
the restriction of $F$ to $\cu$ is a regular function.

Let $F'$ be the regular function on $\cu \times G$ defined by
$F'(xg) = F(x)$; $x \in \cu$, $g \in G$.  As we have seen
in Proposition 1.1, the morphism $\cu \times G \lr \wh{X}$,
$(x,g) \longmapsto x \circ g$, factors as follows $\cu \times
G \lr \cu \times^B G \lr \wh{X}$.  We see that $F'$ goes down 
to the function $F$ on $\wh{X}$, in particular it goes down
to a function $F_1$ on $\cu \times^B G$. 
Since $\cu \times G \lr \cu \times^B G$ is a locally trivial
fibration, we see that $F_1$ is a regular function.  Since
$F_1$ goes down to $F$ and $\cu \times^B G \lr \wh{X}$ is a 
proper morphism, we see that if $\Gamma$ is the graph of $F$,
then the canonical morphism $\Gamma \lr \wh{X}$ is {\it proper
and bijective}.  Then from Lemmas \ref{lem:bijective1} and \ref{lem:bijective2} 
given below, 
it follows that $F^q$ is a regular function on $\wh{X}$
for $q = p^r$, $r \gg 0$.

\begin{lem}\label{lem:bijective1}
{\em{(}}Cf.~\cite[pp.\,260--261,\,Lemma\,1.4]{K1}{\em{)}} Let $j: P_1 \lr P_2$ be a proper, bijective
morphism of algebraic schemes.  Then the following hold:
\begin{enumerate}
\item If $f$ is a regular function 
on $P_1$, then $f^{p^r}$ is the pull-back of a regular function on
$P_2$ for $r \gg 0$.  
\item If $g_1$  and $g_2$ are regular
functions on $P_2$ such that $j^*(g_1) = j^*(g_2)$, then
$g_1^{p^r} = g_2^{p^r}$ for $r \gg 0$.
\end{enumerate}
\end{lem}
\proof  Since $(g_1 - g_2)^{p^r} = g_1^{p^r}-g_2^{p^r}$,
to prove the last assertion, it suffices to show that if 
$j^*(g) =0$ for a regular function $g$ on $P_2$, then
$g^{p^r} =0$ for $r \gg 0$.  The question is local so that
we can suppose that $P_2 = \spec A$, $P_1 = \spec
B$ and $j$ is given by a homomorphism $j^*: A \lr B$.
The hypothesis implies that $I = {\rm ker}~j^*$ is
nilpotent and since $g \in I$, $g^{p^r}=0$ for $r\gg 0$.
This proves the last assertion.

Now the question of proving the first assertion becomes local.
To see this let $\{ U_i\}$ be an open cover of $P_1$ and
$\{ U_i'\}$ the open cover of $P_2$ such that $U_i =
j^{-1} (U_i')$.  Let $f_i$ be the restriction of $f$
to $U_i$ and suppose that there exist $g_i$ in $U_i'$ such
that $f_i = j^*(g_i)$.  Then we see that for $r \gg 0$,
$g_i^{p^r}=g_j^{p^r}$ in $U_i'\cap U_j'$ so that 
$(g_i)^{p^r}$ define a regular function $g$ on $P_2$
and we have $f^{p^r} = j^* (g^{p^r})$.

We claim that we can suppose that $P_1$ is reduced.  For,
consider $(P_1)_{\rm{red}} \stackrel{i}{\lr} P_1 
\stackrel{j}{\lr} P_2$, where $i$ is the canonical
morphism $(P_1)_{\rm{red}} \lr P_1$.  Let $f' =i^* (f)$, 
$g$ a regular function on $P_2$ such that $(j \circ i)^*
(g) = f'$ and $f_1 = j^* (g)$.  Then $i^* (f_1) = i^* (f)$
so that by the above considerations we have $f_1^{p^r}
= f^{p^r}$ for $r\gg 0$, which implies that $j^*(g^{p^r})
= f^{p^r}$ for $r \gg 0$.  The claim follows.

We claim that we can also suppose that $P_2$ is reduced.  
Then we have the factorisation for $j$ ($P_1 = (P_1)_{\rm{red}}$)

\[
{\xymatrix{
P_1 \ar[dr]_{j'} \ar[r]^j & P_2 \ar[d]^i \\
& (P_2)_{\rm{red}}
}}
\]

Let $g'$ be a regular function on $(P_2)_{\rm{red}}$ such that
$f=(j')^*(g')$.  We see that there is a regular function
$g$ on $P_2$ such that $i^* (g) = g'$, for the considerations
being local (as observed above) we can suppose that $P_2
= \spec A$ so that $(P_2)_{\rm{red}} = \spec A/I$ and the existence
of $g$ follows.  We have then $j^*(g) =f$ and the claim 
follows.

We claim that $P_1$ (and hence $P_2$) can also be taken 
to be irreducible.  Let $(P_1)_i$ ($1\leq i \leq r)$ 
be the irreducible components of $P_1$ and
$(P_2)_i = j ((P_1)_i)$ the irreducible components
of $P_2$.  Let $f_i$ be the restriction of $f$ to
$(P_1)_i$ and $g_i$ a regular function on $(P_2)_i$
such that $j^* (g_i) = f_i$.  Then $g_i$ and $g_j$
coincide set theoretically on $(P_2)_i \cap (P_2)_j$
so that as we saw above in the proof \Pref{prop:L-trivial}, 
$(g_i)^{p^r}$ patch up to define a regular function $g$
on $P_2$ for $r \gg 0$ and then $f = j^*(g)$, which
proves the claim.

Thus we can suppose that $P_1$ and $P_2$ are both
reduced and irreducible.  The considerations being
local we can thus suppose that $P_1 = \spec B$, $P_2 = \spec A$,
where $A$ and $B$ are integral domains and $j^*: A\lr B$ is an
inclusion, $B$ being integral over $A$.  If $k(A)$ and $k(B)$
are the quotient fields of $A$ and $B$ respectively, 
our hypothesis makes $k(B)$ a purely inseparable 
extension over $k(A)$ of degree $q = p^r$.  Then
$B^q \subset k (A)$ and elements of $B^q$ are integral
over $A$.  Then if $B_1 = A \cdot B^q$, we have
$A \hra B_1 \hra B$,  $k(B_1) = k(A)$ and 
$\spec B_1 \lr \spec A$ is proper bijective and to
prove the lemma it suffices to consider this case.

Thus to prove the lemma we can also suppose that $k(A)
= k(B)$.  Then the conductor $C = \{ b \in B \mid 
b B \subset A\}$ is a non-zero ideal in $B$ (as well as
$A$) and the map of schemes
\[\spec B/C \lr \spec A/C\]
 induced by $A/C \hra B/C$
is a proper bijection and $\dim B/C < \dim B$. 
We now prove the lemma by induction on the dimension 
of $B$.  Let $\ol{f}$ be the image of $f \in B$ in 
$B|_C$.  Then by the induction hypothesis, there is
a $\ol{g} \in A|_C$ such that $\ol{g}$ maps to
$\ol{f}^q$, $q = p^r$, $r\gg 0$.  Let $g$ be a representative
of $\ol{g}$ in $A$.  Then we see that $g + \theta = f^q$
for $\theta \in C$ and this proves the lemma. \qed

Recall that a line bundle $L$ on a projective scheme $P$
is said to be {\it semi-ample} if given $x \in P$,
there is a section of a power of $M$ which does not vanish at $x$. 
\begin{lem}\label{lem:bijective2}
Let $j: P_1 \lr P_2$ be a proper, bijective morphism
of algebraic schemes, $L$ a line bundle on $P_2$
and $M$ the line bundle $j^*(L)$ on $P_1$.
Then given a section `$s$' of $M$, $s^q$ is the
pull-back of a section of $L^q$, $q = p^r$, $r\gg 0$.
In particular, we have
\[\mbox{$M$ semi-ample $\Longleftrightarrow L $ semi-ample}.\]
\end{lem}
\proof
The implication $\Longleftarrow$ is trivial.
The reverse implication is an easy consequence of \Lref{lem:bijective1}.
In greater detail, let $\theta_{ij}$ be the transition functions of $L$ 
with respect to an covering $\{ U_i\}$ of $P_2$
and $\{U_i'\}$ the open covering of $P_1$ given by
$U_i' = j^{-1} (U_i)$.  Then $\theta_{ij}' = j^* (\theta_{ij})$
are transition functions of $M$ with respect to the covering
$\{U_i'\}$.  The section `$s$' of $M$ is given by regular
functions $s_i$ in $U_i'$ such that 
$$s_j = s_i \theta_{ij}' ~~{\rm in}~~ U_i' \cap U_j'.$$
Then by the above lemma, $s_i^q$ is the pull-back of a 
regular function $t_i$ in $U_i$, $q=p^r$, $r\gg 0$.
Then the pull-backs of $t_j$ and $t_i \theta_{ij}^q$
coincide on $U_i \cap U_j$, since $s_j^q = s_i^q \theta_{ij}^q$.
Thus again by applying \Lref{lem:bijective1} (by taking a suitable $q$th
power), we can indeed suppose that $t_j = t_i \theta_{ij}^q$ in
$U_i \cap U_j$ i.e. $\{ t_i\}$ patch upto define a section
`$t$' of $L^q$ and we see that $s^q$ is the pull-back of the
section `$t$'.  \qed

\bigskip

It remains to prove the \Pref{prop:L-trivial} for the case $X_1$
(or $\wh{X}_1$).  First observe that if we have a (regular)
section $s$ of $L$ on $X_1$, then a suitable power of $s$
extends to the cone $\wh{X}_1$.  Thus it suffices to prove
the slightly weaker version that there exists a $G$-invariant
section $s$ on $X_1$ of a suitable power of $L$ such that 
$s(x) \neq 0$ $\forall x \in X_1^{ss}$.  Since $X_1^{ss}$ is
a single orbit, $X_1^{ss}$ is smooth.  Let $p: X_2 \lr X_1$
be the normalisation of $X_1$.  We have a canonical 
$G$-action on $X_2$ and $p$ is a $G$-morphism.  Then the 
pull-back $p^*(L)$ of $L$ is ample on $X_2$ and if $X_2^{ss}$
denotes the semi-stable points of $X_2$ for $p^*(L)$, one knows
that
$$ x \in X_2^{ss} \Llr p(x) \in X_1^{ss}. \leqno {\rm (iv)}$$
Now $p : X_2^{ss} \lr X_1^{ss}$ is an isomorphism and
$X_2^{ss}$ is a $G$-orbit.  Suppose that there is a $G$-invariant
section of $p^*(L)$ (or a suitable power) such that $s(x) \neq 0$
$\forall x \in X_2^{ss}$.  Then $s$ vanishes on $X_2 \setminus
X_2^{ss}$ and if $J$ denotes the ideal sheaf on $X_2$ obtained 
as the ideal sheaf on $X_1 \setminus X_1^{ss}$ on $X_1$, then a 
suitable power of $s$ belongs to $J$ (the support of $J$ is
$X_2 \setminus X_2^{ss}$).  Then we see that this power of $s$ 
comes from a section of a power of $L$ on $X_1$, having the
required properties.

Thus without loss of generality, we can suppose that $X_1$
is {\it normal}.

Since $X_1^{ss}$ is a $G$-orbit, $X_1^{ss} \simeq H \setminus G$,
the restriction of $L^n$ to $X_1^{ss}$ is trivial ($n \gg 0$). 
Hence without loss of generality we can suppose that there is a
regular $G$-invariant section of $L$ on $X_1^{ss}$.  We shall
now show that $s$ (or a suitable power of $s$) extends to a 
regular section on $X_1$.  This assertion is true if $G$ happens
to be a torus, in particular $\bg_m$.  (For given $x \in X_1^{ss}$,
we can find a $G$-invariant regular section $\theta$ of $L^r$ (for some
$r$) such that $\theta (x) \neq 0$.  Then $s^r$ and a constant multiple
of $\theta$ coincide on $X_1^{ss}$ and hence they coincide
everywhere).  Suppose that $s$ is not regular on $X_1$.  Then
$s$ has a {\it pole} at some $x_0 \in X_1 \setminus X_1^{ss}$
(i.e. $s$ is not regular at $x_0$ and has no indeterminacy at
$x_0$).  Choose some $x \in X_1^{ss}$.  Then there is a $K$-valued
point $g \in G(K)$ (or possibly by going to a finite extension
of $K$, $K$ being the quotient field of the d.v.r. $A$ as in
\Rref{rem:dvrbase}) such that $\ds \lim_{t \ra 0} x \circ g = x_0$.
Now by the Iwahori-Matsumoto theorem, we have $g = P \lambda Q$,
where $P, Q \in G(A)$ and $\lambda : \bg_m \lr G$ is a non-trivial
1-PS, so that the image of $\lambda$ is a subgroup $H$ of $G$,
$H \simeq G_m$ and $\lambda$ defines a $K$-valued point of $G$.
We see that $\lim_{t \ra 0} x \circ P \lambda = x_0' \in X_1
\setminus X_1^{ss}$ and $s$ has a pole at $x_0'$.  Let $P_0
= \lim_{t \ra 0} P$.  Then $P_0 \in G$ and $x \cdot P_0 \in 
X_1^{ss}$.  We have $x \circ P \lambda = x \cdot P_0 \cdot
P_0^{-1} P \lambda$ and $\lim_{t \ra 0} P_0^{-1} P = I$
(identity element of $G$).  Thus without loss of generality,
we can suppose that $\lim_{t\ra 0} x \circ \lambda = x_0 \in
X_1 \setminus X_1^{ss}$ and $s$ has a pole at $x_0$.  Let 
$Z$ denote the closure of the ``$\lambda$-orbit through $x$''
i.e. the closure of $xH$ in $X_1$ and take the restriction
$\theta$ of $s$ to $Z$.  Then $\theta$ is $H$-invariant and has a
pole at $x_0$ (which is in $Z$) and $Z^{ss} = xH$.  This
reduces the problem to $G = \bg_m$, where as we observed
above, the assertion is true i.e. a power of $\theta$ is
regular on $X_1$.  This leads to a contradiction of the
hypothesis that $s$ has a pole at $x_0$. This proves
\Pref{prop:L-trivial}. \qed

\begin{defi} Let $A=k[|t|]$ and let $K$ be its quotient field. We say ``$X^{ss}$ mod $G$ is
proper" if given $x\in X^{ss}(K)$ there exists a finite extension $K'$ of $K$ and $g\in G(K')$
such that $x\smcirc g \in X^{ss}(A')$ where $A'$ is the finite extension of $A$ associated with
$K'$.
\end{defi}
\begin{cor}\label{cor:proper}\cite[p.\,526,\,Thm.\,4.1]{S2}
Let $X$ be a closed $G$-invariant subscheme of
$\bp (V)$.  Then 
$X^{ss}$ mod $G$ is proper.

\end{cor}
{\proof As in \Rref{rem:dvrbase}, let $Z_K$ denote the closure
of the $G_K$ orbit $Z^0_K = G_K$ in $X_K$ and $Z_A$ the flat
closure of $Z_K$ in $X_A$.  Now $x$ is in $Z_K^{ss}$ (in the sense
of \Ssref{ss:basechange}).  Let $\ol{K}$ denote the algebraic closure of
$K$, then taking $Z_{\ol{K}}$, $X_{\ol{K}}$ etc. we see that
$Z_{\ol{K}}^0$ is an orbit under $G_{\ol{K}}$, $Z_{\ol{K}}$ is
the closure of $Z_{\ol{K}}^0$ and $Z_{\ol{K}}^0 \subset Z_{\ol{K}}^{ss}$.
Then by \Pref{prop:L-trivial}, there exists a (regular) section $s$ of 
$L_{\ol{K}}$ (or a suitable power of $L_{\ol{K}}$) on $Z_{\ol{K}}$
such that $s$ is $G_{\ol{K}}$ invariant and $s$ is non-vanishing 
at every point of $Z_{\ol{K}}^{ss}$.  Now $s$ is defined over a
finite extension of $K$ and thus without loss of generality we can
suppose that there exists a section $s$ of $L_K$ over $Z_K$
which is $G_K$ invariant and does not vanish at every closed or 
geometric point of $Z_K^{ss}$.

We claim that to prove the corollary, it suffices to prove that 
$Z_k^{ss}$ is not empty, where $Z_k$ denotes the closed fibre
of $Z_A \lr \spec A$ (we see that this is also necessary).
Suppose that $Z_k^{ss} \neq \es$, then we see easily that there
exists $y \in Z_A (A)$ (or possibly we may have to go a finite
extension of $K$, take the integral closure of $A$ in this 
extension etc.) such that its restriction to the closed fibre is 
$y_0$ and restriction to the generic fibre is in $Z_K^0$ (we
see that $y \in Z_A^{ss} (A) \subset X_A^{ss} (A))$.  We now
see that there exists $g \in G(K)$ such that $x \circ g = y$
(again we may have to go a finite extension of $K$ for this
assertion).

Let $p: \wt{Z}_A \lr Z_A$ be the normalisation of $Z_A$.
We take semi-stable points in $\wt{Z}_A$ with respect to the
pull-back of $L_A$ by $p$.  It suffices to prove that
$\wt{Z}_k^{ss} \neq \es$, where $\wt{Z}_k$ denotes the closed 
fibre of $\wt{Z}_A \lr \spec A$, for $\wt{Z}_k^{ss} = 
p^{-1}(Z_k^{ss})$.  Let $\wt{s}$ denote the pull-back of
the section $s$ of $L_K$ on $Z_K$.  Then by multiplying
$\wt{s}$ by a suitable power of the uniformising parameter
$\pi$ of $A$, we can suppose that $\wt{s}$ is a regular
section of $L_A$ on $\wt{Z}_A$.  The restriction
$\wt{s}_k$ of $\wt{s}$ to $\wt{Z}_k$ is obviously
$G (= G_k)$ invariant.  Then if $\wt{s}_k$ does not vanish
at $z \in Z_k$, then obviously $z \in \wt{Z}_k^{ss}$.  
Thus it suffices to prove that $\wt{s}_k$ is not identically
zero on $\wt{Z}_k$.  We claim that there exist $q,r$ such that
$q \geq 0$ and $r \geq 1$ such that $\pi^{-q}s^r$ is regular on 
$\wt{Z}_A$ and its restriction to $\wt{Z}_k$ is not identically zero.  
This would suffice.  To prove this claim let $m_i$
(resp. $n_i$) denote the order of vanishing of $\pi$ (resp.
$\wt{s}$) along the irreducible components $Z_k^i$ of
$\wt{Z}_k$ ($\wt{Z}_k$ is a divisor of the normal scheme
$\wt{Z}_A$).  Then the order of $\pi^{-q}s^r$ along $Z_k^i$
is $(rn_i-qm_i)$ (note that $m_i \geq 1$ and $n_i \geq 0$
$\forall i$).  To prove the claim we have only to show
that $(rn_i - qm_i) \geq 0$ $\forall i$ and $\exists$ $i_0$
such that $rn_{i_0} -q m_{i_0} =0$.  If we take $i_0$ such
that $\ds \frac{n_{i_0}}{m_{i_0}} = \min_i \frac{n_i}{m_i}$ and
$q,r$ such that $\ds \frac{q}{r} = \frac{n_{i_0}}{m_{i_0}}$, 
we are done.  This proves the corollary.\qed

\begin{rem}\label{rem:inv-extn}{\em
We see that in the proof of \Pref{prop:L-trivial}, the argument actually
shows that if $\wh{X}$ and $\wh{X}_1$ are closed $G$-invariant
subsets of $V$ such that $\wh{X} = \cu (\wh{X}_1)$ and $f$ is
a $G$-invariant function on $\wh{X}_1$, then $f$ extends
(uniquely) to a $G$-invariant function on $\wh{X}$ (and similar
assertion for the case $X$, $X_1$).  Again there is a more general 
assertion concerning the extension of a $G$-invariant section 
of $L$ on $X^{ss}$ to the whole of $X$ (see \cite[p.\,526,\,Theorem\,4.1]{S2}).}
\end{rem}

\section{{\bf{Stratified spaces}}}\label{s:strata}

Throughout this section we fix a {\em reduced and irreducible} quotient data
\stepcounter{thm}
\begin{equation*}\label{eq:De-red-irr}\tag{\thethm}
\De = (X,\,L,\bp(V),\,X^{ss}\xrightarrow{\alpha} Y).
\end{equation*}
This section consists largely of definitions which facilitate the proof of the main theorem.

\subsection{Stratification of $(X^{ss}/\negthickspace/G)_{\text{top}}$}
\begin{lem}\label{lem:U} 
There exists a non-empty open subset $U\subset Y$ such that
\begin{enumerate} 
\item[(i)] $U$ has the structure of a  $k$-variety;
\item[(ii)] The map $\alpha_{X,U}\colon \alpha_X^{-1}(U)\to U$ is a map of $k$-schemes;
\item[(iii)] The structure sheaf on $U$ is given by the $G$-invariant direct image of
$\co_{\alpha_X^{-1}}(U)$, i.e., by the sheaf $V\mapsto \Gamma(\alpha_X^{-1}(V), \co_X)^G$
{\em{(}}$V$ open in $U${\em{)}}.
\end{enumerate}
Moreover, a power of the  line bundle $L$ on $X^{ss}$ descends to a line bundle on $U$ and
by shrinking $U$, this descended line bundle may be assumed to be trivial.
\end{lem}

\proof 
Let $X$ be a closed $G$-invariant subset of $\bp (V)$  
and $\alpha$: $X^{ss} \lr Y$ (or $\alpha_X$) the map as above.  
Let us suppose that $Y$ is {\it irreducible}.  Let $e$ be the
maximum of the dimension of the closed orbits in $X^{ss}$.  Recall
(\Cref{cor:polystable}  to \Pref{prop:U(S)-closed}) that the subset $X^{ps,e}$ of closed
orbits of dimension $e$ is open in $X^{ss,e}$, the closed subset of orbits
of dimension $\leq e$ in $X^{ss}$ (we endow $X$, $X^{ss}$ etc. with
the reduced scheme structures).  The canonical map $X^{ss,e} \lr Y$ is
surjective and $X^{ss,e} \setminus X^{ps,e}$ maps onto a closed subset
of $Y$ which is not the whole of $Y$ so that $X^{ps,e}$ maps onto a
non-empty open subset of $Y$.  Then by the irreducibility of $Y$,
it is seen without much difficulty (noting that $G$ is irreducible)
that there exists a non-empty
{\it open irreducible $G$-invariant} subset of $X^{ps,e}$ which maps onto an
open subset ($\neq \es$) of $Y$ and we denote the closure of this 
subset in $X$ by $X_0$.  Thus we have a {\it closed, irreducible}
$G$-invariant subset $X_0$ of $X$ such that the canonical map
$\alpha$ ($=\alpha_{X_0}$): $X_0^{ss} \lr Y$ is surjective and we have
a non-empty, $G$-invariant open subset of $X_0^{ps,e}$ as well as of
$X_0^{ss}$, consisting of the orbits of dimension $e$, which are
closed in $X_0^{ss}$ and hence in $X^{ss}$ (polystable orbits of dimension
$e$, $X_0^{ps,e} = X^{ps,e} \cap X_0^{ss}$).  By the existence of 
{\it generic quotients}, we see easily that there is a $G$-stable open 
subset of $X_0^{ss}$, which is therefore of the form 
$\alpha_{X_0}^{-1} (U)$, where $U$ is open in $Y$ such that the 
geometric quotient $\alpha^{-1}_{X_0} (U)$~mod~$G$ exists and is 
a variety (for example by arguments in \cite{K1} or \cite{S4}).  Since 
$\alpha^{-1}_{X_0} (U)$~mod~$G$ identifies with
$U$ (set-theoretically and topologically) we can thus endow $U$
with the structure of a variety.  We denote this variety structure 
on $U$ by $U_0$ so that the structure sheaf on $U_0$ is given by the 
sheaf of $G$-invariant regular functions on $\alpha_{X_0}^{-1} (U)$.
Observe that the generic fibre of $\alpha_{X_0}^{-1} (U) \lr U$
is of the form ``$G$~mod~$H$'', ($H$ the isotropy group at a 
polystable point) and therefore by suitably shrinking $U$, we can
suppose that there exists a $G$-invariant regular section of a
(suitable power of $L$) on $\alpha_{X_0}^{-1} (U)$, which is
everywhere non-vanishing (see \Pref{prop:L-trivial} and 
\Ssref{ss:basechange} above), in particular,
(a suitable power of) $L$ on $\alpha^{-1}_{X_0} (U)$ descends to 
$U$.  Consider the map $\alpha_X : X^{ss} \lr Y$ and the map
$\alpha_X: \alpha_X^{-1} (U) \lr U$ induced by it.  We claim that 
given a $G$-invariant regular function on $\alpha_{X_0}^{-1} (U)$,
a suitable $p$th power of this function extends to a regular 
($G$-invariant) function on $\alpha_X^{-1} (U)$.  This follows by an
extension of the arguments in \Pref{prop:U(S)-closed} and \Pref{prop:L-trivial} , 
since we have
$$\cu (\alpha_{X_0}^{-1} (U)) = \alpha_X^{-1} (U)$$
where we define the left hand side by
$$ \cu (\alpha_{X_0}^{-1} (U)) = X^{ss} \setminus 
(X_0^{ss} \setminus \alpha_{X_0}^{-1} (U)).$$
Thus we have a variety structure on $U$ defined by $G$-invariant
regular functions on $\alpha_X^{-1} (U)$, which we denote
simply by $U$.  Then we have a morphism $U_0 \lr U$ which
is proper and bijective.  Now the generic fibre of 
$\alpha_X^{-1} (U) \lr U$ is a single semi-stable equivalence class
and by \Pref{prop:L-trivial} , by suitably shrinking $U$, we see that
there is a regular $G$-invariant section of $L$ on $\alpha_X^{-1}
(U)$, which is everywhere non-vanishing.  In particular, we see
that $L$ on $\alpha_X^{-1} (U)$ goes down to a line bundle on $U$
(we can assume it to be trivial).
\qed

\begin{defi}\label{def:dimY} We define the dimension of $Y$ to be $\dim{Y}=\dim{U}$, where 
$U$ is the non-empty
open set of \Lref{lem:U}. The function field of $Y$, $k(Y)$ is defined by the formula $k(Y)=k(U)$. 
Note that these notions are well-defined.
\end{defi}

\begin{rem}\label{rem:strata} 
{\em Note that $Y$ posseses
a stratification $\{U_\lambda\}$ by locally closed subsets $U_\lambda\subset Y$, 
each of which have the structure of a $k$-variety. Indeed, setting $Y'= Y\smallsetminus U$,
and $X'=\alpha_X^{-1}(Y')$, we note that $Y'= ((X')^{ss}(L)/\negthickspace/G)_{\text{top}}$. Since
$\dim{Y'} < \dim{Y}$, by induction we see the stratification of $Y$. In particular, by repeatedly
using the last statement of \Lref{lem:U}, we see that there is a positive integer $n$ such that
$L^n\vert_{\alpha_X^{-1}(y)}$ is trivial on every fibre $\alpha_X^{-1}(y)$ of $\alpha_X$.} 
\end{rem}

\begin{defi}\label{def:strat-map}
Let $\{U_\lambda\}$ be the stratification in \Rref{rem:strata}.  
Let $W$ be a $k$-scheme, together with locally closed 
subschemes $\{W_\lambda\}$ which give a stratification of $W$ (i.e., the $W_\lambda$
are disjoint and as sets their union equals $W$). A continuous map of 
$q\colon W\to Y$ is said to be a {\em stratified morphism} or a {\em stratified map} if, for each index 
$\lambda$,  $W_\lambda= q^{-1}(U_\lambda)$ and the map
\[ W_\lambda = q^{-1}(U_\lambda) \xrightarrow{q\vert_{W_\lambda}} U_\lambda \]
is a map of $k$-schemes.
\end{defi}

\begin{defi}\label{def:gen-finite} Let $\{U_\lambda\}$ 
a stratification of $Y$ as in
\Rref{rem:strata}, and $W$ a $k$-variety. A stratified map $q\colon W\to Y$ is said to be 
{\em generically finite} if $W$ is {\em complete} and
\[q^{-1}(U)\xrightarrow{\phantom{X}{\text{via}}\,q\phantom{X}} U\] 
is a generically finite map of $k$-varieties.
\end{defi}

Note that if $q\colon W\to Y$ is generically finite, then the field extension $k(Y)\to k(W)$ is finite.

\begin{lem}\label{lem:Gam-action} Let $W$ be
a {\em normal projective variety} and $q\colon W\to Y$ be a generically finite map. 
Then there exists a normal projective variety $W'$ such that
\begin{enumerate}
\item[(i)] There is a surjective map of varieties $W'\to W$ such that the composite
\[ W'\to W \xrightarrow{q} Y\]
is generically finite. We denote this composite $q'$.
\item[(ii)] The field $L\set k(W')$ is the least normal extension of $k(Y)$ which contains $k(W)$.
\item[(iii)] If $\Gamma\set {\text{Aut}}_{k(Y)}(L)$, then $\Gamma$ acts on $W'$ and the orbits
of $W'$ under $\Gamma$ are contained in the fibres of $q'$. In other words, as a continuous
map, $q'$ is $\Gamma$-invariant for the trivial action of $\Gamma$ on $Y$.
\item[(iv)] There exists a non-empty open subset $V$ of the ``big stratum" $U$ of $Y$ such that
for $v\in V$ the fibre $(q')^{-1}(v)$ is an orbit of $\Gamma$.
\end{enumerate}
\end{lem}

\proof Set $n=[k(W):k(Y)]_s$, the separable degree of $k(W)/k(Y)$. Let $R=k(W)^{\otimes n}$,
the tensor product being taken over $k(Y)$. Then $L$ can be expressed  as a quotient
of the $k(Y)$-algebra $R$. Let $I=\ker(R\to L)$. The symmetric group $S_n$ on $n$-letters
acts on $R$ in an obvious way, the group $\Gamma$ can be identified with the subgroup of
$S_n$ which leaves $I$ invariant.

The $n$-fold fibre-product $W_n \set W\times_Y\dots\times_YW$ makes sense as a closed
subscheme of the $n$-fold product $W\times_k\dots\times_kW$. Now, $S_n$ also acts
on $W_n$ and one can find a closed
subvariety $W'$ of $W_n$, such that, if one has $k(W')=R/I=L$. Through projection into one of
the factors we get a surjective morphism $W'\to W$ of projective varieties which at the generic
points is represented by the function field extension $K(W)\to L$. Let
The group $\Gamma$ acts on $W'$ and this action lifts to the normalisation
of $W'$. Replacing $W'$ by this normalisation, and composing $W'\to W$ with $q$ we get
a genrically finite map
\[q'\colon W'\to Y.\]
We claim that if $x,y\in W'$ with $y=\gamma x$, $\gamma\in\Gamma$ then $q'(x)=q'(y)$. This
is certainly true over $U$. Thus $q'$ and $q'\smcirc\gamma$ agree over $(q')^{-1}(U)$ which
is a $\Gamma$-stable set. Hence $q'$ and $q'\smcirc\gamma$ agree on all of $W'$. 
This proves (i)---(iii) and assertion (iv) is clear.
\qed

\medskip

The notion of a stratified map seems to depend on the chosen stratification of $Y$. However
there are maps $q\colon W\to Y$ which are stratified, whatever be the stratification of $Y$ into
a finite number of locally closed subsets $U_\lambda$ with a $k$-variety structure such that
 $\co_{U_\lambda}$ is the sheaf of $G$-invariant sections on $\alpha_X^{-1}(U_\lambda)$.
 This happens when $W$ carries a principal $G$-bundle $Z\to W$ and $Z$ has an equivariant
 map into $X^{ss}$. For this it is perhaps best to state matters in terms of stacks. 
 
 \subsection{Stack theoretic interpretation}\label{ss:Ystack}
 Recall that in this section our fixed quotient data \eqref{eq:De-red-irr} is 
 reduced and irreducible. There is a distinguished
 index $\lambda=\circ$ such that $U_\circ$ is an open subset of $Y$.
 Set $X_\lambda\set \alpha_X^{-1}(U_\lambda)$ and $\alpha_\lambda\colon X_\lambda\to U_\lambda$
 equal to
 the restriction of $\alpha_X$ to (the locally closed) subvariety $X_\lambda$ of $X^{ss}$. Since
 the structure sheaf $\co_{U_\lambda}$ is given by the sheaf of $G$-invariant sections of
 $X_\lambda$, we see that $\alpha_X\colon X^{ss}\to Y$ is a stratified map, whatever be the
 stratification on $Y$ induced by \Lref{lem:U} (i.e., by the choices of the open subsets predicted
 by \Lref{lem:U}).

 Recall that the stack $[X^{ss}/G]$ is the fibered category over $k$-schemes such that for any
 $k$-scheme $Q$, $[X^{ss}/G](Q)$ is the category whose objects  are pairs $(P\xrightarrow{\delta} Q,\, P\xrightarrow{j} X^{ss})$
 with $\delta\colon P\to Q$ a principal $G$-bundle and $j\colon P\to X^{ss}$ a
 $G$-equivariant map. Morphisms in $[X^{ss}/G](Q)$ are isomorphisms of such pairs.
 $[X^{ss}/G]$ is a stack with respect to either the fppf or the \'etale topologies
 on the category of $k$-varieties.
 One could restrict to (and we will do so) base $k$-schemes $Q$ which
 are $k$-varieties.
 
 Given the data $(Y,\,\{U_\lambda\})$, one can define another stack (also denoted $Y$)
 and a map of stacks $\gamma\colon[X^{ss}/G]\to Y$ as follows.
 For a $k$-variety $W$, let 
 \stepcounter{thm}
 \begin{equation*}\label{eq:Ystack}\tag{\thethm}
 {\underline{Y}}(W)=\{q\colon W\to Y\,\vert\, q\,\, {\text{is a stratified map}}\}.
 \end{equation*}
 Let $W'\to W$ be a faithfully flat and finitely presented map of varieties, $W''\set W'\times_WW'$ the
 two-fold product of $W'$ with itself over $W$, and 
 $p_i\colon W''\to W'$, $i=1,2$ the two projections. It is easy to see that if we have a stratified
 map $q'\colon W'\to Y$ such that $q'\circ p_1=q'\circ p_2$, then there is a unique stratified
 map $q\colon W\to Y$ such that $q'$ is the composite $W'\to W\xrightarrow{q} Y$. In fact
 descent works on each stratum of $W'$, giving us, set-theoretically a map $q\colon W\to Y$,
 which is a map of varieties on each $q^{-1}(U_\lambda)$. Now $q'\colon W'\to Y$ is continuous
 and $W'\to W$ being fppf, $W$ has the quotient topology from $W'$. It follows that $q$ is
 continuous. Thus $\underline{Y}$ is also a stack --- since it is a sheaf of sets on
 the fppf site over $k$-varieties. 
 
 Given a principal $G$-bundle $\delta\colon P\to Q$ and a $G$-equivariant map $j\colon P\to X^{ss}$,
 we have, clearly a continuous map 
 \stepcounter{thm}
 \begin{equation*}\label{eq:gamma1}\tag{\thethm}
 q=q(\delta,j)\colon Q\to Y.
 \end{equation*}
 We point out that $q(\delta,\,j)$ does not depend on the stratification $\{U_\lambda\}$ that
 of $Y$ that we've fixed. We now proceed to show that $q$ is stratified (and hence
 is stratified for every stratification as in \Rref{rem:strata}).
  If the locally closed
 subset $Q_\lambda=q^{-1}(U_\lambda)$ of $Q$ is given the reduced topology, and $P_\lambda
 =\delta^{-1}(Q_\lambda)$, then we have a $G$-equivariant map of varieties $j_\lambda\colon
 P_\lambda\to X_\lambda$.  The map $\delta_\lambda\colon P_\lambda\to Q_\lambda$ (given by
 restricting $\delta$ to $P_\lambda$) is a principal bundle. It follows
 (from the universal properties of the geometric quotient $Q_\lambda$ of $P_\lambda$ by $G$) that
 we have a map of varieties $q_\lambda\colon Q_\lambda\to Y$. Moreover, clearly,
 $q_\lambda=\gamma\vert_{Q_\lambda}$. Thus $q(\delta,\,j)$ is a stratified map.
 The maps $q(\delta,\,j)$ induce a map of stacks
 \stepcounter{thm}
 \begin{equation*}\label{eq:gamma2}\tag{\thethm}
 \gamma\colon [X^{ss}/G]\to {\underline{Y}}.
 \end{equation*}
 such that if $f\colon Q\to [X^{ss}/G]$ is the classifying map for the data $(\delta,\,j)$ then
 $\gamma\smcirc f$ is the map whose underlying stratified map is $q(\delta,\,j)$.
 
 From now on, we will identify the stack ${\underline{Y}}$ with $Y$. The main theorem (when proved)
 will show that $Y$ is actually a variety. But for the moment, it is not even clear it is an algebraic
 stack (i.e., that it has a presentation by a scheme which is smooth over it).
 
 Note that $[X^{ss}/G]$ can morally be regarded as a ``stratified  space", the ``stratification"
 being given by $\{[X_\lambda/G]\}$. Note also that morally
 $\gamma\colon [X^{ss}/G]\to Y$ can be regarded a ``stratified map".
 
 \begin{conv}\label{conv:L-power} We observed in \Rref{rem:strata} that there exists a positive
 integer $n$ such that $L^n$ is trivial on the fibres of $\alpha_X$. Since $X^{ss}(L)=X^{ss}(L^n)$,
 and $X^s(L)=X^s(L^n)$, we may assume, by replacing $L$ by $L^n$ if necessary, that $L$ is
 trivial on the fibres of $\alpha_X$. For the rest of our discussion we make this assumption.
\end{conv}
 
 In view of the above convention, if $(P\to Q, P\to X^{ss})$ is an object in $[X^{ss}/G](Q)$,
 then the pull back of $L$ to $P$ is trivial on the fibres of $P\to Q$ and hence descends to
 a line bundle on $Q$. In other words $L\vert_{X^{ss}}$ gives us a line bundle on
 $[X^{ss}/G]$. We make the following definition:
 
 \begin{defi}\label{def:L/G} The line bundle on $[X^{ss}/G]$ induced by $L$ on $X^{ss}$ via descent
 will be denoted $L_{/G}$.
 \end{defi}

Now suppose we have a $Q$-valued point of $[X^{ss}/G]$ where $Q$ is a $k$-scheme. In other
words suppose we have a principal $G$-bundle $\delta\colon P \to Q$ as well as a
$G$-equivariant map $j\colon P\to X^{ss}$. The line bundle $j^*L$ descends to a line bundle $L_Q$
on $Q$. In fact the line bundle $L_Q$ can be regarded as $f^*L_{/G}$ where $f\colon Q\to [X^{ss}/G$
is the natural map induced by the data $(\delta,\,j)$. There is also a resulting stratified map
$i\colon Q\to (X^{ss}/\negthickspace/G)_{\text{top}}$. In fact $i=\gamma\circ f$. 
It turns out that $L_Q$ depends only
on $i$ rather than on $(\delta,\,j)$, if $Q$ is a normal variety. In other words, if $Q$ is normal and
the base of another principal $G$-bundle
with an equivariant map to $X^{ss}$ such that the resulting stratified map
$Q\to (X^{ss}/\negthickspace/G)_{\text{top}}$ is again $i$, then, the line bundle induced on $Q$ from
this principal bundle is again $L_Q$. One can say it better in the following way:

\begin{prop}\label{prop:Pindep} Let $(X,\,L)$ be a $G$-pair and suppose 
$Y= (X^{ss}/\negmedspace/G)_{\text{\em{top}}}$ is irreducible.
Suppose $Q$ is a normal $k$-variety 
and we have
two maps $Q\xrightarrow{f} [X^{ss}/G]$ and $Q\xrightarrow{g} [X^{ss}/G]$ such that
$\gamma\circ f =\gamma\circ g$ as stratified maps from $Q$ to $Y$, i.e., suppose the diagram
\[
{\xymatrix{
Q \ar[r]^{f} \ar[d]_g & [X^{ss}/G] \ar[d]^{\gamma} \\
[X^{ss}/G] \ar[r]_{\gamma} & Y
}}
\]
commutes. Then, with $L_{/G}$ as in \ref{def:L/G} {\em{(}}see also \ref{conv:L-power}{\em{)}}, we 
have
\[ f^*L_{/G} \simeq g^*L_{/G}. \]
\end{prop}

\proof
 Note that we have a canonical topology on the product space 
$Q \times Y$.  The closed subsets in $Q \times Y$  
are the images in $Q \times Y$ of the closed $G$-invariants
subsets of $Q \times X^{ss}$ (for the trivial action of $G$ on $Q$),
under the canonical map $Q \times X^{ss} \lr Q \times Y$ (or equivalently
closed subsets of $Q \times X^{ss}$ which are saturated for this map).
The closed subsets of $Q \times Y$ can also be defined as follows.
Observe that $Q \times Y$ has also a stratification by subspaces which
are schemes (we take the product of $Q$ with the subschemes defining
a stratification of $Y$).  We say that a subset $C$ of $Q \times Y$
is closed if the intersection $C_i$ of $C$ with every subscheme in
this stratification is closed and satisfies a property for limits 
expressed by a valuation criterion as follows.
Let $A$ be a d.v.r. and $K$ its quotient field.  An $A$-valued point
$\theta$ of $Q \times Y$ (written $\spec A \lr Q \times Y$) is one 
which can be ``lifted'' to an $A$-valued point of $Q \times X^{ss}$.

\[
{\xymatrix{
\spec\,A \ar[r] \ar[dr]_\theta & Q\times_k X^{ss} \ar[d]^{1\times\alpha}\\
& Q\times Y
}}
\]

We see that $\theta$ can be viewed as a set theoretic map of a 
neighbourhood of the closed point of the smooth curve defined by
$\spec A$, into $Q \times Y$.  Suppose that $\theta$ defines a
$K$-valued point of some $C_i$.  Then the condition to be imposed
is that $\theta$ maps the closed point of $\spec A$ to a point of
$C$. 

Let $f\colon Q\to [X^{ss}/G]$ be given by the data $(P\xrightarrow{\delta} Q,\, P\xrightarrow{j} X^{ss})$.
Let $i\colon Q\to Y$ be the resulting map (i.e. $i=\gamma \circ f$). We have a commutative
diagram:
\[
{\xymatrix{
P \ar[d]_{\delta} \ar[r]^{j} & X^{ss} \ar[d]^{\alpha}\\
Q \ar[r]_{i} & Y
}}
\]

The map $i$ is ``nice'' in the following sense.
In fact, it is not difficult to see that if $\Gamma_0$ is the 
graph of $j$ and we take the canonical action of $G \times G$
on $P \times X^{ss}$ and take the closure $\Gamma$ of $\Gamma_0$.
$(G \times G)$ in $P \times X^{ss}$, then the canonical image
of $\Gamma$ in $Q \times Y$ is the graph of $i$.  This can be
expressed more intutively as follows.  Let $\theta$ be a $K$-valued
point of $Q \times Y$, which is the image of a $K$-valued point
$\phi = (\phi_1,\phi_2)$ of the graph of $j$, $\phi_1$ being 
a $K$-valued point of $P$ and $\phi_2$ a $K$-valued point
of $X^{ss}$.  Since $P~{\rm mod}~G$ and $X^{ss}~{\rm mod}~G$
are proper, we see that there exist $K$-valued points $g_1,g_2$
of $G$ such that $\psi_1 = \phi_1 \circ g_1$ (resp. $\psi_2=\phi_2 
\circ g_2$) is an $A$-valued point of $P$ (resp. $X^{ss}$).
We see that $j(\psi_1) = \psi_2 \circ g$, where $g = g_2^{-1}g_1$
is a $K$-valued point of $G$ (using the fact that $j(\phi_1) =
\phi_2$ and $j$ is $G$-equivariant).  Now $j(\psi_1)$ is an 
$A$-valued point of $X^{ss}$ so that the closed points of $X^{ss}$
determined by $j(\psi_1)$ and $\psi_2$ are semi-stably equivalent
i.e. they are in the same fibre of $X^{ss} \lr Y$ and determnie
a point $y$ of $Y$.  Let $x$ be closed point of $Q$ determined by
$\psi$.  Then we see that $(x,y)$ is in the graph of $i$ and this 
essentially shows that the graph of $i$ is closed.

Consider the ``base change'' $X' = Q \times_Y X^{ss}$.
Then $X'$ is a closed $G$-stable subset of $Q \times_Y X^{ss}$
(since $Q \times X^{ss} \lr Y \times Y$ is continuous and the 
diagonal is closed in $Y \times Y$).  We endow $X'$ with the 
canonical reduced structure as a closed subscheme of $Q \times
X^{ss}$.  We have a canonical morphism $j'$ of $P$ into $X'$. In fact we
have a commutative diagram (with the parallelogram being ``cartesian"):
\stepcounter{thm}
\begin{equation*}\label{diag:parallel}\tag{\thethm}
{\xymatrix{
P \ar[rr]^{j'} \ar[dr]_{\delta} \ar@/^2pc/[rrrr]^{j}
& & X' \ar[dl]_{\alpha'} \ar[rr]^{j''}  & & X^{ss} \ar[dl]^{\alpha}\\
& Q \ar[rr]_{i} & & Y
}}
\end{equation*}
Our strategy for establishing the Proposition is to prove that for a non-empty set $V$ of $Q$ we
have a natural isomorphism---behaving well with respect to restrictions to open subschemes of
$V$:
\stepcounter{thm}
\begin{equation*}\label{iso:delta-indep}\tag{\thethm}
\Gamma(V,\,f^*L_{/G}) \iso \Gamma((\alpha')^{-1}(V),\,(j'')^*L).
\end{equation*}
The right side depends only on the map $i\colon Q\to Y$ and not on $f$ (for
the space $X'$ depends only on $i$), whence establishing \eqref{iso:delta-indep} is equivalent
to establishing the Proposition.

We denote by $\ol{X}$ the closure of $j'(P)$ in $X'$ and endow it
with the reduced subscheme structure.  Then we have a commutative
diagram:
\[
{\xymatrix{
P \ar[rr]^{\ol{j}}  \ar [dr]_{\delta} & & \ol{X} \ar[dl]^{\ol{\alpha}} \\
& Q &
}}
\]
where all the maps are morphisms and $\ol{j}$ is $G$-equivariant.
Observe that if $x \in Q$, $y =i (x)$ and $\ol{X}_x$ denotes the
fibre of $\ol{\alpha}$ over $x$, then $\ol{X}_x$ can be identified with
a closed $G$-stable subset of $X_y^{ss}$.  We denote by the same
$L_{\ol{X}}$ the line bundle on $\ol{X}$, obtained as the pull-back of $L$
on $X^{ss}$ by the canonical morphism $\ol{X} \lr X^{ss}$.
Then  $\ol{j}^*L_{\ol{X}}=j^*L$.

We can assume without loss of generality that $i(Q)=Y$.  Hence
$V_0 = i^{-1} (U)$ is a non-empty open subset of $Q$, where
$U$ is the open subset of $Y$ as in \Lref{lem:U}.  By {\em loc.cit.}
a suitable power of $L$ descends to $U$, and this descended bundle
can be 
assumed to be trivial.  Without loss of generality, we assume that $L$ itself descends
to the trivial bundle $\co_U$ on $U$. Hence the restriction of $L_{\ol{X}}$ to
$(\ol{\alpha})^{-1} (V_0)$ (as well as $j^*L\vert_{\delta^{-1} (V_0)}$) can
be assumed to be trivial.  We see, in particular, that the 
$G$-invariant sections of $j^*L$ 
on $\delta^{-1} (V_0)$ identify with the $G$-invariant
sections of $L_{\ol{X}}$ on $(\ol{\alpha})^{-1} (V_0)$. In other words, we have
a canonical isomorphism:
\stepcounter{thm}
\begin{equation*}\label{iso:delta-alpha1}\tag{\thethm}
\Gamma(V_0,\,f^*L_{/G}) = \Gamma(\delta^{-1}(V_0),\,j^*L)^G \iso 
\Gamma((\ol{\alpha})^{-1} (V_0),\,L_{\ol{X}})^G.
\end{equation*}
We point out that for any open set $V$ of $Q$, by definition of $f^*L_{/G}$ we have 
$\Gamma(V,\,f^*L_{/G}) = \Gamma(\delta^{-1}(V),\,j^*L)^G$. For an open subscheme $W$
of $\ol{X}$, let $W_n$ denote its normalisation, and let $\varphi\colon\ol{X}\to X$ denote
the normalisation map.
We claim that for every
non-empty open set $V$ of $Q$ the following holds:
\stepcounter{thm}
\begin{equation*}\label{eq:LQ}\tag{\thethm}
\Gamma(V,\,f^*L_{/G})  \iso 
\Gamma([({\ol{\alpha}})^{-1}(V)]_n,\,\varphi^*L_{\ol{X}})^G
\end{equation*}
In other words, the claim is that the sheaf of sections of $f^*L_{/G}$  on $Q$ 
identifies with the sheaf
of $G$-invariant regular sections on the normalisation of
$\ol{X}$. This is to be thought of as the first step towards proving \eqref{iso:delta-indep}, 
which as we have noted, is sufficient for the Proposition.

Since a $G$-invariant section of $\ol{j}^*L$ on $(\ol{\alpha})^{-1} (V)$
maps by $j^*$ to a $G$-invariant section of $j^*L$ on $\delta^{-1}
(V)$, to prove \eqref{eq:LQ} we have only to show that if `$s$' 
is a section of $L$ on $V$, then it comes from a $G$-invariant 
section of $L$ on the normalisation of $(\ol{\alpha})^{-1}(V)$.
Now the restriction of `$s$' to $(V \cap V_0)$ identifies with
a $G$-invariant section of $L$ on $(\ol{\alpha})^{-1} (V \cap V_0)$, 
so that we can consider `$s$' to be a $G$-invariant {\it rational}
section of $L$ on $(\ol{\alpha})^{-1} (V)$.  Suppose that it is 
{\it not} regular on the normalisation of $(\ol{\alpha})^{-1} (V)$.
Then we have a non-empty {\it polar divisor} $D$ for `$s$' so that
there is an $x_1 \in D$ such that the ``value'' of $s$ at $x_1$
is $\infty$ (i.e. $x_1 \in D$ and is not a point of indeterminacy).
Let $y = \ol{\alpha}(x_1)$.  We shall now prove the following:
\stepcounter{thm}
\begin{equation*}\label{eq:subset}\tag{\thethm}
\ol{X}_y \subset D. \,\, {\text{In fact if
$z \in \ol{X}_y$ the value of `$s$' at $z$ is $\infty$.}}
\end{equation*}
Now \eqref{eq:subset} $\Lr$ \eqref{eq:LQ}, for if $P_y$ is the fibre of $P$ over $y$,
$P_y \subset \ol{X}_y$, and `$s$' considered as a section of $L$
over $\delta^{-1} (V)$ would have a polar divisor containing 
$P_y$, so that we have a contradiction since `$s$' has been
supposed to be regular over $\delta^{-1} (V)$.  Thus `$s$'
is regular on the normalisation of $(\ol{\alpha})^{-1}(V)$.

To prove \eqref{eq:subset}, let us first observe that the generic fibre of
$\ol{\alpha} : \ol{X} \lr Q$ has a dense $G$-orbit, since $\ol{X}$ 
has been defined as the closure of the image of $P$ in
$Q \times_Y X^{ss} = X'$.  Suppose now $x_2 \in \ol{X}_y$
and the ``value'' of $s$ at $x_2$ is {\it not} $\infty$.
Hence either `$s$' is regular at $x_2$ or it is a point of
indeterminacy.  We then see easily that there is an $A$-valued
point $\theta_2$ ($A$ a.d.v.r. as usual) of $\ol{X}$ such that
the closed point of $\spec A$ maps to $x_2$, the $K$-valued
point ($K$-the quotient field of $A$) of $\ol{X}$ determined by
$\theta_2$ is in the dense $G$-orbit of the generic fibre 
of $\ol{\alpha}$ and the ``restriction of $s$ to $\theta_2$''
is regular, i.e. if we set $s(x_2) = \lim_{t \ra 0} s|_{\theta_2}$,
then $s(x_2) \neq \infty$.  Now we can find another $A$-valued point
$\theta_1$ of $\ol{X}$ such that the closed point maps to $x_1$ and the
$K$-valued point defined by $\theta_1$ is in the dense $G$-orbit
defined by $\theta_1$.  We see that there is a $K$-valued point
$g$ of $G$ such that $\theta_1 \cdot g = \theta_2$ (we may have to
go to a finite extension of $K$).  Sicne `$s$' is $G$-invariant
and $s(x_1) = \infty$, it follows that $s(x_2) = \infty$, which
is a contradiction and the assertion \eqref{eq:subset} follows, whence
so does \eqref{eq:LQ}

Let $\alpha': X' \lr Q$ denote the canonical morphism and 
$s$ a regular $G$-invariant section of $L$ on $V$.  Then
$(\ol{\alpha})^{-1} (V)$ is a closed $G$-invariant subset of
$(\alpha ')^{-1} (V)$ and by the arguments as in \Pref{prop:L-trivial} or,
more precisely \Rref{rem:inv-extn},
we see that $s$ raised to a suitable power of $p$ extends to 
a regular $G$-invariant section of $L$ on $(\alpha ')^{-1} (V)$.
But since $s$ is already a rational section of $L$ and $Q$ is normal, 
we see that $s$ can indeed be identified with a $G$-invariant regular
section $L$ on $(\alpha ')^{-1} (V)$ i.e. in \eqref{eq:LQ} we can
replace $\ol{\alpha}$ by $\alpha'$.  This establishes \eqref{iso:delta-indep},
whence the Proposition.
\qed

\begin{cor}\label{cor:Pindep} Let $H$ be a finite group of $Y$-automorphisms of $Q$. 
Then $H$ lifts to an action of $L$.
\end{cor}

\proof The action of $H$ on $Q$ extends to $Q\times_kX^{ss}$ (as well as the line bundle
$L'$ on $Q\times_kX^{ss}$ which is the pull back of $L$), by taking the trivial action on
$X^{ss}$. Since $H$ is a group of automorphisms, we see that the action of $H$ on $Q$
lifts to an action on $X'=Q\times_YX^{ss}$, i.e. $X'$ is as in the proof of the Proposition.
This action of $H$ clearly commutes with the action of $G$. Hence it acts on the sheaf
of invariant sections of $L$ on $X'$. The assertion now follows from the isomorphism in
\eqref{eq:LQ}. \qed

\section{\bf{Reduction to the case stable = semi-stable}}\label{s:s=ss}

In this section we revisit and modify certain technically crucial Lemmas and Propositions in \cite{S2},
namely Lemma\,3.2, Propositions\,5.1 and 5.3 of {\em ibid}.                                
In this section we allow our group $G$ to be a reductive algebraic group (relaxing our
requirement that $G$ be semi-simple). Fix
a maximal torus $T$ in $G$, and a Borel subgroup $B$ of $G$ with $B \supset T$. The 
notations we use are as follows:

\begin{itemize}
\item $\Gamma(T)$ will denote the co-root lattice of $T$, i.e., $\Gamma(T)$ will denote the
abelian group of {\em one parameter subgroups} (1-PS) $\lambda\colon {\mathbb{G}}_m\to T$
of $T$.
\item $E\set \Gamma(T)\otimes_{\mathbb{Z}}{\mathbb{R}}$.
\item $C(B)$ will denote the Weyl chamber in $E$ associated to $B$, and ${\overline{C(B)}}$ will
denote its closure in $E$.
\end{itemize}

Now for any projective algebraic scheme $X$ with a $G$-action which is linear with respect to
an ample bundle $L$, and for a fixed $x\in X$, the function $\mu^L(x,\lambda)$ is a integral valued 
function on $\Gamma(T)$ and we extend this function to an ${\mathbb{R}}$-valued function on $E$
by setting
\[ \mu^L(x,a\lambda) = a \mu(x,\lambda), \qquad a\in{\mathbb{R}}, \quad \lambda\in \Gamma(T). \]

\stepcounter{subsection}

\begin{lem}\label{lem:cones} Let $X_1,\dots,X_d$ be projective algebraic schemes on which
$G$ acts, such that for each $i=1,\dots,d$, the action in linear with respect to an ample line bundle
$L_i$ on $X_i$. Then there exist a finite number of closed convex cones $C_\alpha$ contained in
${\overline{C(B)}}$ (resp. $E$)---independent of $x_i\in X_i$---such that

\begin{enumerate}
\item each $C_\alpha$ is the intersection of a finite number of half spaces in $E$, the half spaces
being of the form $\{x\in E\,\vert\, \theta(z)\ge 0\}$, $\theta$ being a linear form on $E$, with 
integral coefficients (with respect to a given basis) and
\[{\overline{C(B)}} ({\text{resp. $E$}}) = \cup_\alpha C_\alpha; \]
\item in every $C_\alpha$,  $\mu^{L_i}(x_i,\_)$ is linear for fixed $x_i\in X_i$, $i=1, \dots, d$.
\end{enumerate}
\end{lem}

\begin{lem}\label{lem:S-exist1} Let $p\colon Z\to X$ be a $G$ invariant
morphism between projective
algebraic $G$-schemes, the action being linear with respect to the ample line bundle $L$ on $X$
and $M$ on $Z$. Write $aL+bM$ for the line bundle $p^*(L)^a\otimes M^b$ ($a,b\in {\mathbb{Z}}$).
Then there exists a finite set $S\subset \Gamma(T)\setminus \{0\}$ such that for every line bundle
$N$ of the form $N=aL+bM$, where $a$ and $b$ are positive integers, we have
\begin{enumerate}
\item $\mu^N(z,\lambda)\ge 0,\,\, \forall\, \lambda\in \Gamma(T) \Longleftrightarrow 
\mu^N(z,\lambda) \ge 0,\,\, \forall\, \lambda\in S$.
\item If $\mu^N(z,\lambda)\ge 0$, $\forall\, \lambda\in \Gamma(T)$ and $\mu^N(z,\lambda_\circ)=0$ for
some $\lambda_\circ\in \Gamma(T)\setminus \{0\}$, then for some 
$\lambda_\circ\in S$, we have $\mu^N(z, \lambda_\circ)=0$.
\end{enumerate}
\end{lem}

\proof 
By \Lref{lem:cones}, we can subdivide $E$ into a finite number of closed convex cones $C_\alpha$
(each $C_\alpha$ an intersection of closed half-spaces), independent of $x\in X$ or $z\in Z$ such that
$\mu^L(x,\_)$ and $\mu^M(z,\_)$ are linear on each $C_\alpha$.

Now, for $z\in Z$, and for $N=aL+bM$, with $a,b$ positive integers, we have
\[\mu^N(z,\_)=a\mu^L(x,\lambda)+b\mu^M(z,\lambda), \qquad \lambda\in\Gamma(T),\,x=p(z).\]

Let $S_\alpha$ be a finite set of generators (over ${\mathbb{R}}^+$) for the cone $C_\alpha$. We can
choose $S_\alpha$ with integral coordinates, i.e., $S_\alpha\subset\Gamma(T)$.
Since $\mu^L(x,\_)$ and $\mu^M(z,\_)$ are linear in each $C_\alpha$, therefore
$\mu^N(z,\_)$ is linear on each $C_\alpha$. It follows that the finite set $S=\cup_\alpha
S_\alpha$ satisfies the assertion of the Lemma.
\qed

\begin{rem}\label{rem:S-exist1} {\em It is worth pointing out that \Lref{lem:S-exist1} also implies that
if $\mu^N(z,\,\lambda)\le 0$ for some non-trivial 1-PS $\lambda$, then
$\mu^N(z,\,\lambda_\circ)\le 0$ for some $\lambda_\circ \in S$. Indeed, if $\mu^N(z,\,\lambda_\circ)>0$
for every $\lambda_\circ\in S$, then by part (1), we have $\mu^N(z,\,\lambda_\circ)\ge 0$ for every
non-trivial 1-PS $\lambda_\circ$. If further, $\mu^N(z,\,\lambda_\circ)=0$ for any non-trivial
1-PS $\lambda_\circ$
then by part (1), $\mu^N(z,\,\lambda_\circ)=0$ for some $\lambda_\circ \in S$, giving the required
contradiction.}
\end{rem}

\begin{prop}\label{prop:b/a-small1} Let $p\colon Z \to X$ be a G-morphism between projective
algebraic schemes, the action being linear with respect to an ample bundle $L$ on $X$ and
a relatively ample bundle (with respect to $p\colon Z\to X$) $M$ on $Z$. Write
$N(a,b)=aL+bM$ for positive integers $a,b$. Then for $\frac{b}{a}$ sufficiently small, we have
\[ p^{-1}(X^s(L))\subset Z^s(N(a,b)) \subset Z^{ss}(N(a,b)) \subset p^{-1}(X^{ss}(L))\]
(the second inclusion is obvious).
\end{prop}

\proof First, we can find $a_1$ and $b_1$ such that $M_1=a_1L+b_1M$ is ample on $Z$.
Replacing $M_1$ by $M$, we can suppose without loss of generality that $M$ is in fact ample
on $Z$ and not merely relatively ample. We are thus in the situation of \Lref{lem:S-exist1}.
Let $S\subset \Gamma(T)\setminus \{0\}$ be the finite set satisfying the conclusions
of {\em loc.cit}. We have, for $z\in Z$
\[\mu^{N(a,b)}(z,\,\lambda) = a \mu^L(x,\,\lambda) + b\mu^M(z,\,\lambda), x=p(z). \]
Choose positive integers $a,b$ such that
\[ \bigg\lvert \mu^M(z,\,\lambda)\frac{b}{a}\bigg\rvert <1 \,\,{\text{for all $\lambda\in S$ and $z\in Z$.}} \]
We can do this for the functions $\mu^M(\_,\lambda)$, $\lambda\in S$ are finite in number, each
continuous on the compact set $Z$. Now (writing $N=N(a,b)$)
\[\mu^N(z,\,\lambda) = a \biggl[\mu^L(x,\,\lambda) + \frac{b}{a}\mu^M(z,\,\lambda)\biggr]. \]
Since $a$ is positive, the sign of $\mu^N(z,\,\lambda)$ is the same as the sign of
$\mu^L(x,\,\lambda) +\frac{b}{a}\mu^M(z,\,\lambda)$. Since $\mu^L(x,\,\lambda)$ is an
integer and $\lvert (b/a)\mu^M(z,\,\lambda)\rvert <1$ for $\lambda\in S$, it follows that for
$\lambda\in S$
\begin{itemize}
\item[(a)] $\mu^N(z,\,\lambda) \ge 0 \Longrightarrow \mu^L(z,\,\lambda)\ge 0$;
\item[(b)] $\mu^N(z,\,\lambda)\le 0 \Longrightarrow \mu^L(z,\,\lambda)\le 0$.
\end{itemize}

We first show that $Z^{ss}\subset p^{-1}(X^{ss}(L))$. Suppose on the contrary, there exists
$z\in Z^{ss}(N)$ such that $x=p(z)\notin X^{ss}(L)$. Then there exists $x_1=x\circ g$, $x_1\in X$, 
$g\in G$ such that $\mu^L(x_1,\,\lambda_\circ)<0$ for some non-trivial $\lambda_\circ \in \Gamma(X)$.
By our choice of $S$ (cf.\, \Lref{lem:S-exist1}), we may suppose $\lambda_\circ\in S$. Let
$z_1=z\circ g$ so that $x_1=p(z_1)$. Since $z\in Z^{ss}(N)$, therefore $z_1\in Z^{ss}(N)$, whence
$\mu^N(z_1,\,\lambda_\circ)\ge 0$.
However, $\mu^L(x_1,\,\lambda_\circ)<0$, whence by (a) above, $\mu^N(z_1,\,\lambda_\circ)<0$,
giving the required contradiction.

Next we show that $p^{-1}(X^s(L)) \subset Z^s(N)$. Suppose we can find a $z\in p^{-1}(X^s(L))$ such
that $z\notin Z^s(N)$. Then $x=p(z)\in X^s(L)$, whence 
\begin{equation*}\tag{*}
\mu^L(x\circ g, \lambda) > 0\qquad {\text{($\forall\, g\in G$ and 
$\forall\, \lambda_\circ \in \Gamma(T)\setminus \{0\}$)}}.
\end{equation*}
 On the other hand, $z\notin Z^s(N)$, whence there exists a $g\in G$
and a $\lambda\in \Gamma(T)\setminus \{0\)$ such $\mu^N(z\circ g,\,\lambda)\le 0$.
By \Lref{lem:S-exist1} (see \Rref{rem:S-exist1})
we conclude that $\mu^N(z\circ g,\,\lambda_\circ)\le 0$ for some
$\lambda_\circ\in S$. By (b) above, this implies that for this $g$ and this $\lambda_\circ$
we have $\mu^L(x\circ g,\, \lambda_\circ) \le 0$, contradicting (*).
\qed

One can modify the proof of \Lref{lem:S-exist1} in an obvious way to get:

\begin{lem}\label{lem:S-exist2} Let $G$ act on the projective schemes $X_1,\dots, X_l$ and $Y$
and suppose the action on $X_i$ ($i=1,\dots,l$) is linear with respect to an ample line bundle 
$L_i$ on $X_i$.
Suppose further that $\pic(Y)$ is generated by ample line bundles $M_1,\dots, M_r$ and that the
action of $G$ on $Y$ is linear with respect to each $M_j$, $j=1,\dots, r$. Let $Z_i\set X_i\times Y$,
$i=1,\dots,l$. Note that $Z_i$ is a $G$-scheme with respect to the diagonal action
of $G$ and that this action is linear with respect to every line bundle of the form $aL_i+bM$,
$M\in\pic(Y)$, $a,b\in{\mathbb{Z}}$. Then there exists a finite subset $S\subset \Gamma(T)\setminus \{0\}$ of 1-PS
$\lambda$ of $T$---independent of $i=1,\dots,l$ and of $z_i\in Z_i$---such that for every line
bundle $N_i(a,b)$ of the form
\[N_i(a,b) (=N_i) = aL_i+bM \qquad (a,b\in{\mathbb{N}})\]
with $M$ an ample line bundle on $Y$, we have
\begin{enumerate}
\item $\mu^{N_i}(z_i,\,\lambda)\ge 0, \,\,\forall \lambda\in \Gamma(T) \Longleftrightarrow
\mu^{N_i}(z_i,\lambda)\ge 0, \,\, \forall \lambda\in S$;
\item If $\mu^{N_i}(z_i,\,\lambda)\ge 0$, $\forall \lambda\in\Gamma(T)$ and 
$\mu^{N_i}(z_i,\,\lambda_\circ)=0$ for some $0\neq \lambda_\circ\in\Gamma(T)$, then
$\mu^N(z_i,\,\lambda_\circ)=0$ for some $\lambda_\circ\in S$.
\end{enumerate}
\end{lem}

\proof By \Lref{lem:cones}, we can subdivide $E$ into a finite number of convex cones $C_\alpha$
(each $C_\alpha$ an intersection of closed half-spaces)---independent of $x_i\in X_i$, $i=1,\dots,l$
and $y\in Y$---on which $\mu^{L_i}(x_i,\_)$, $i=1,\dots,l$ and $\mu^{M_j}(y,\_)$, $j=1,\dots,r$ are
linear. 

Write $z_i=(x_i,y)$, $x_i\in X_i$, $y\in Y$. Now $M=\sum_jb_jM_j$ for $b_j\in{\mathbb{Z}}$.
The proof of the Lemma is
is almost identical to the proof of \Lref{lem:S-exist1} once one observes that
\[\mu^{N_i}(z_i,\,\lambda) = a \mu^{L_i}(x_i,\,\lambda) + b(\sum_jb_j\mu^{M_j}(y,\,\lambda)) 
\qquad (\lambda\in \Gamma(T)),\]
which implies that on each $C_\alpha$, $\mu^{N_i}(z_i,\_)$ is linear. 
\qed

\begin{convs}\label{conv:L+M} Let $L$ and $M$ be line bundles on a scheme $S$ and consider
the positive integral linear combination $N(a,b)=aL+bM$, i.e. $N(a,b)=
L^{\otimes}a\otimes M^{\otimes}b$. For discussions involving notions which
are stable under ``multiplication" of $N(a,b)$ by a poisitive integer (e.g., ampleness,
semi-ampleness, nefness, bigness of $N(a,b)$) we will often write 
$N(a,b)=N_{\frac{b}{a}}=L+\frac{b}{a}M$. In particular, for such discussions, if
$\epsilon$ is a positive rational number, the symbol
\[N_\epsilon=L+\epsilon M\]
is a convenient shorthand.  This shorthand can be extended to include
pairs of maps $f\colon S\to T$ and $g\colon S \to U$
with $L$ a line bundle on $T$ and $M$ a line bundle on $U$, so that
$L+\epsilon M$ represents $f^*L+ \epsilon g^*M$.
\end{convs}

\begin{defi}\label{def:epsilon} Let $Y$ be a projective $G$-scheme as in
\Lref{lem:S-exist2}. Let $(X,\,L)$ be a $G$-pair, $Z=X\times_kY$, and
$p\colon Z\to X$ the projection to the first factor. 

A line bundle $M$ on $Y$ is said to be {\em stablizing for $(X,\,L)$} if 
there exists $\epsilon_0>0$ such that
\begin{itemize}
\item[(a)] $M$ is ample.
\item[(b)] For $0<\epsilon\le\epsilon_0$, with $N_\epsilon\set L+\epsilon M$,
\begin{enumerate}
\item[(i)] $p^{-1}(X^s) \subset Z^s(N_\epsilon)$ and  $p(Z^{ss}(N_\epsilon))\subset X^{ss}$;
\item[(ii)] $Z^s(N_\epsilon)=Z^{ss}(N_\epsilon)$;
\item[(iii)] If $0<\epsilon'\le\epsilon_0$, then $Z^s(N_\epsilon)=Z^s(N_{\epsilon'})$.
\end{enumerate}
\end{itemize}
If $M$ is stabilizing for $(X,\,L)$ and $\epsilon_0$ is a positive rational number upper bound
as in (b) above, then we say $M$ is {\em $\epsilon_0$-stabilizing for $(X,\,L)$}.
\end{defi}

Just as \Lref{lem:S-exist1} is used to prove \Pref{prop:b/a-small1}, one can use \Lref{lem:S-exist2}
to prove the following proposition.

\begin{prop}\label{prop:b/a-small2} Let $X_1,\dots,X_l$, $Y$;  $L_1,\dots, L_l$;  $M_1,\dots, M_r$
satisfy the hypotheses of \Lref{lem:S-exist2}. As in {\em loc.cit.}, let $Z_i=X_i\times Y$. Suppose that
given a finite set $S\subset \Gamma(T)\setminus \{0\}$, there exists an ample line bundle $M$
on $Y$ such that $\mu^M(y,\,\lambda)\neq 0$ for every $y\in Y$ and $\lambda\in S$. Then
there exists an ample line bundle $M$ on $Y$ which is stabilizing for all the G-pairs
$(X_i,\,L_i)$, $i=1,\dots,l$.
\end{prop}

\proof 
Once we find an ample $M$ which satisfies (ii) in \Dref{def:epsilon} for $\epsilon$ sufficiently small, then by \Pref{prop:b/a-small1}, (i) of the definition also follows (for this $M$) for $\epsilon$ sufficiently small. We first prove show that (ii) is satisfied for $\epsilon$ sufficiently small.

We are in a situation where \Lref{lem:S-exist2} applies. Let $S$ be a finite set of non-trivial
elements of $\Gamma(T)$ satisfying the conclusions of \Lref{lem:S-exist2}. Note that $S$ is
independent of $i\in\{1,\dots,l\}$. By our hypotheses there exists an ample line bundle
$M$ on $Y$ such that $\mu^M(y,\,\lambda)\neq 0$ for every $y\in Y$ and $s\in S$. Let 
$N_i=N_i(a,b)\set aL_i+bM$, for positive integers $a,b$. For
$z_i\in Z_i$, we have, by \Lref{lem:S-exist2}

\begin{itemize}
\item[(i)] $\mu^{N_i}(z_i,\,\lambda)\ge 0$ for every $\lambda\in\Gamma(T)$ $\Longleftrightarrow$
$\mu^{N_i}(z_i,\,\lambda)\ge 0$ for every $\lambda\in S$.
\item[(ii)] If $\mu^{N_i}(z_i,\,\lambda)\ge 0$ for every $\lambda\in \Gamma(T)$ and 
$\mu^{N_i}(z_i,\,\lambda_\circ)=0$ for some $\lambda_\circ\in\Gamma(T)\setminus\{0\}$, then 
$\mu^{N_i}(z_i,\,\lambda_\circ)=0$ for some $\lambda_\circ\in S$.
\end{itemize}

Choose positive integers $a$, $b$ such that
\begin{equation*}\tag{*}
\Big\lvert\frac{b}{a}\mu^M(y,\,\lambda)\Big\rvert < 1 \qquad \forall\, \lambda\in S, \quad y\in Y.
\end{equation*}
This is possible because for fixed $\lambda\in S$, $\mu^M(\_,\,\lambda)$ is a continuous
function on the compact space $Y$, and $S$ is finite. Suppose, by way of contradiction, 
$z_i\in Z_i^{ss}(N_i)$ and $z_i\notin Z_i^s(N_i)$. Then there exists $z_i'\in Z_i^{ss}(N_i)$,
$z_i'=z_i\circ g$ for some $g\in G$, and non-trivial $\lambda\in \Gamma(T)$ such that
$\mu^{N_i}(z_i',\,\lambda)=0$. By (ii) above, we may assume that $\lambda_\circ\in S$. Now
\[\mu^{N_i}(z_i',\,\lambda_\circ) = a \Bigl[\mu^{L_i}(z_i',\,\lambda_\circ) + \frac{b}{a}\mu^M(y_i',\,\lambda_\circ)\Bigr] \qquad y'_i=p_i(z_i') \]
and hence, $\mu^{N_i}(z_i',\,\lambda_\circ)=0$ implies that $\mu^{L_i}(z_i',\,\lambda_\circ) +
(b/a)\mu^M(y'_i,\,\lambda_\circ)=0$. By (*) this means that $\mu^{L_i}(z'_i,\,\lambda_\circ)=0$
and $\mu^M(y'_i,\,\lambda_\circ)=0$. But $M$ has been chosen so that $\mu^M(y,\,\lambda)\neq 0$
for any $y\in Y$ and $\lambda\in S$. This gives the required contradiction. We have therefore
shown that the ample bundle $M$ s=chosen fulfills (i) and (ii) in \Dref{def:epsilon} .

It remains to that this choice of $M$ satisfies (iii). Suppose $\mu_1$ and $\mu_2$ are integers, and 
$C=\{\epsilon >0\,:\, \epsilon\vert\mu_2\vert < 1\}$. We claim that if $\mu_1+\epsilon\mu_2 >0$
for any $\epsilon \in C$, then $\mu_1+\epsilon'\mu_2>0$ for every element $\epsilon'\in C$.
Indeed, let $\epsilon\in C$ be such that $\mu_1+\epsilon\mu_2>0$. Clearly $\mu_1$ cannot
be negative by definition of $C$. If $\mu_1=0$, then, clearly $\mu_2>0$ whence
$\mu_1 + \epsilon'\mu_2 = \epsilon'\mu_2 >0$ for every $\epsilon'\in C$ (in fact for
every $\epsilon'>0$). If, $\mu_1>0$, then $\mu_1\ge 1$. In this case, since
$\mu_1+\epsilon'\vert\mu_2\vert <1$ for $\epsilon'\in C$, it follows that $\mu_1+\epsilon'\mu_2 >0$ 
for such $\epsilon'$. The same argument shows that $\mu_1+\epsilon\mu_2 <0$ for some $\epsilon
\in C$ is equivalent to $\mu_1+\epsilon'\mu_2<0$ for all $\epsilon'\in C$.

We will suppress the index $i\in \{1,\dots, l\}$ and, for example, write $X$ for $X_i$, $Z$ for $Z_i$, etc.
Let $z=(x,y)\in Z=X\times_k Y$. Suppose $z\notin Z^s(N(a,b))=Z^{ss}(N(a,b))$.
Then for some $g\in G$ and
$\lambda\in \Gamma(T)\setminus \{0\}$, we have $\mu^{aL+bM}(z\circ g\, \lambda) < 0$, and
by our choice of $M$, we can, and will, take this $\lambda$  to lie in $S$. Let $z'=z\circ g$, and $x'=
x\circ g$, $y'=y\circ g$. Then setting
 $\mu_1=\mu^L(x',\lambda)$ and $\mu_2=\mu^M(y',\lambda)$, the above argument gives 
 $\mu^{a'L+b'M}(z',\lambda)<0$ for any pair of positive integers
 $a'$, $b'$ with $b'/a'<\epsilon_\circ$. In other words, $z'\notin Z^s(N(a',b'))=Z^{ss}(N(a',b'))$.
 This proves that $M$ is stabilizing for all the $(X_i,\,L_i)$.
\qed

\begin{rem}\label{rem:indep-b/a}
{\em{In \cite[p.\,550,\,Thm.\,7.1]{S2} it is proven (without assuming geometric reductivity of $G$) 
that if $(X,\,L)$ is a pair on which $G$ acts linearly, with $L$ an ample line bundle on 
$X$, and $X$ a normal variety, satisfying $X^s(L)=X^{ss}(L)$, then the geometric quotient
of $X^s(L)$ with respect to $G$ exists as a normal projective variety $Y$ and $L$ descends
to an ample line bundle on $Y$. In fact $Y={\rm{Proj}}(R^G)$, where
$R=\bigoplus_{n\ge 0}\Gamma(X,\,L^n)$. In view of this, if the $Z_i^s(N_i(a,b))$ are non-empty,
and ${\wt{Z}}_i$ the normalisation of $Z_i$, we get by
part (2) of the Proposition the
existence of ${\wt{Z}}^s_i(N_i(a,b))/\negmedspace/G$ as a geometric quotient for each $i$.
Part (3) of the Proposition shows that the quotients 
$W_i=W_{i,M}\set {\wt{Z}}^s_i(N_i(a,b))/\negmedspace/{G}$ for
$b/a$ sufficiently small, do not depend on $(a,b)$, but only on $M$ (and $i$).}}
\end{rem}

\begin{rem} {\em{The hypotheses on $M$ in \Pref{lem:S-exist2} and \Pref{prop:b/a-small2}
may well be unnecessary. To begin with, note that $\mu^M(y,\,\lambda)$ makes sense for
any $M\in\pic(Y)$, whether ample or not, provided the action of $G$ on $Y$ lifts to $M$. Next
note that if $M$ and $M'$ are algebraically equivalent line bundles then $\mu^M(y,\,\lambda)
=\mu^{M'}(y,\,\lambda)$ for every $y\in Y$ and $\lambda\in\Gamma(T)$. Indeed, if $M$ is algebraically
equivalent to zero, then by reducing to the case of curves, it is easy to see that $\mu^M=0$.
Finally note that the Neron-Severi group ${\text{NS}}(Y)$ is finitely generated.}}
\end{rem}

\subsection{Applications}\label{ss:nonempty} Let
\[ \De = (X,\,L,\,\bp(V),\,X^{ss}\xrightarrow{\alpha} Y)\]
be a quotient data. We fix a
Borel subgroup $B\subset G$ and a maximal torus $T$ of $G$ such that $T\subset B$.
As in the beginning of \Sref{s:s=ss}, we denote $\Gamma(T)$ the group of 1-PS of $T$,
and $C(B)$ the positive Weyl chamber associated to $(B,T)$. 
Recall that $B\backslash G$ has the following properties \cite[pp.\,533--534,\,Prop.\,5.3]{S2}:
\begin{itemize}
\item ${\text{Pic}}\,(B\backslash G)$ is finitely generated by ample line bundles which are linear
with respect to the natural action of $G$ on $B\backslash G$.
\item Given a finite set $S\subset \Gamma(T)\smallsetminus \{0\}$, there exists an ample line
bundle $M$ on $B\backslash G$, linear with respect to $G$, such that $\mu^M(y,\,\lambda)\neq 0$
for every $y\in B\backslash G$ and $\lambda\in S$.
\end{itemize}
We are therefore in a position to apply \Pref{prop:b/a-small2} to $Z\set X\times_k B\backslash G$
and deduce the existence of an stabilizing line bundle  (see \Dref{def:epsilon}) 
$M\in {\mathrm{Pic}}\,(B\backslash G)$ for $(X,\,L)$. 

\begin{prop}\label{prop:nonempty} Suppose $\De$ above is a strong quotient data.
Then we can find an $\epsilon_0>0$ and an
$\epsilon_0$-stabilizing $M$ for $(X,\,L)$ 
{\em ($M\in {\mathrm{Pic}}\,(B\backslash G)$)} such that for  $0<\epsilon\le\epsilon_0$ and
with $N_\epsilon\set L+\epsilon M$ the following hold:
\begin{enumerate}
\item[(i)] $Z^s(N_\epsilon)\neq\emptyset$.
\item[(ii)] $Z^s(N_\epsilon)/\negmedspace/G$ exists as a geometric quotient and the natural map
\[ q\colon Z^s(N_\epsilon)/\negmedspace/G \to Y\] 
is surjective.
\item[(iii)] Let $p\colon Z\to X$ be the natural projection. Then, given $x\in X^{ss}$, there exists
$x'\in p(Z^s(N_\epsilon))$ such that $x'$ is semi-stably equivalent to $x$.
\end{enumerate}
\end{prop}

\proof
Write
$\X=\bp(V)$, $\Y=(\bp(V)^{ss}/\negmedspace/G)_{\text{top}}$, $\Z=\X\times_k B\backslash G$.
Let $\tilde{\alpha}\colon \X\to \Y$ be the quotient map for semi-stable equivalence. Let 
${\mathcal{L}}=\co_{\bp(V)}(1)$. 

By \Pref{prop:b/a-small2}
we can find an $\epsilon_0$-stabilizing bundle $M\in {\mathrm{Pic}}\,(B\backslash G)$ for 
$(\X,\,{\mathcal{L}})$ for some $\epsilon_0>0$. Fix a positive rational number
$\epsilon$ which is less than or equal to $\epsilon_0$. To lighten notation we write $\Z^s$ for 
$\Z^s({\mathcal{N}}_\epsilon)$
and $Z^s$ for $Z^s(N_\epsilon)$ where ${\mathcal{N}}_\epsilon$ is the line bundle
${\mathcal{N}}_\epsilon\set{\mathcal{L}}+\epsilon M$ on $\Z$. Let 
$\tilde{p}\colon \X\times_kB\backslash G\to \X$ be the projection. Since 
${\tilde{p}}^{-1}(\X^s)\subset
\Z^s$ and $\X^s\neq \emptyset$, it follows that $\Z^s\neq \emptyset$. Now $\Z$ is normal and
$\Z^s=\Z^{ss}$ whence by \Rref{rem:indep-b/a} we have
a geometric quotient $\W=\Z^s/\negmedspace/G$ which is projective and on to which 
${\mathcal{N}}_\epsilon$ descends as an ample line bundle. 
Let $\tilde{\beta}\colon \Z^s\to \W$ be the quotient map.
There is a natural continuous map $\tilde{q}\colon
\W \to \Y$ such that the diagram
\[
{\xymatrix{
\Z^s \ar[d]_{\tilde{\beta}} \ar[r]^{\tilde{p}} & \X^{ss} \ar[d]^{\tilde{\alpha}}\\
\W \ar[r]_{\tilde{q}} & \Y
}}
\]
commutes. According \Lref{lem:U} there is a non-empty open subset $\U\subset \Y$ which has a 
scheme structure such that
\[\tilde{\alpha}^{-1}(\U)\xrightarrow{\phantom{X}{\text{via}}\,\,\tilde{\alpha}\phantom{X}} \U \]
is a map of schemes, and the structure sheaf
on $\U$ is the sheaf of $G$-invariants of $\co_{\tilde{\alpha}^{-1}(\U)}$. Since $\Y$ is irreducible,
$\U\cap \tilde{\alpha}(\X^s)$ is a non-empty open set and clearly
$\tilde{q}(\W)\supset \U\cap \tilde{\alpha}(\X^s)$. This immediately implies that
$\W \to \Y$ is surjective,  for $\X^{ss}$ mod $G$ is proper by Corollary\,\ref{cor:proper}.  
Now $Y$ is a non-empty closed subset of $\Y$. As a topological space,  
$Z^s= \Z^s\cap Z= \Z^s\cap {\tilde{p}}^{-1}(X)=\Z^s\cap {\tilde{p}}^{-1}(X^{ss})$. The last equality
uses the fact that $\tilde{p}(Z^{s})=\tilde{p}(Z^{ss})\subset X^{ss}$. Since
$X^{ss}$ is saturated, $X^{ss}={\tilde{\alpha}}^{-1}(Y)$, whence,
$Z^s =\Z^s\cap {\tilde{p}}^{-1}(Y)= {\tilde{\beta}}^{-1}({\tilde{q}}^{-1}(Y))$.
Since $\tilde{\beta}$ and $\tilde{q}$ are both surjective, $Z^s\neq \emptyset$. This proves (i).
 
Since ${\tilde{\beta}}$ maps distinct $G$-orbits into 
distinct points of $\W$, the ``orbit space" of $Z^s$ can be identified, 
as a topological space, can be identified with ${\tilde{q}}^{-1}(Y)$. It is not difficult
to show that this implies the existence of a geometric quotient $\beta\colon Z^s\to W$
such that the natural map $W\to {\tilde{q}}^{-1}(Y)$ is bijective (as a map of sets).
This gives the required surjectivity of $q\colon W\to Y$ proving (ii).

Part (iii) follows from (ii). In greater detail 
first pick $z\in \Z^s$ such that $\tilde{q}(\tilde{\beta})(z)=\tilde{\alpha}(x)$ using the 
surjectivity of $\tilde{q}\smcirc\tilde{\beta}$. Then $z$ actually lies in $Z^s$ since
$\tilde{\alpha}(x)\in Y$. Now set $x'=\tilde{p}(z)=p(z)$.
\qed

\subsection{Elimination of finite isotropies}\label{ss:elimination} 
In this subsection we summarize the
results in 
\cite[pp.\,536--544,\,\S\,6]{S2}.
Let $(X,L)$ be a $G$-pair with $X$ {\em normal} and $X^s\neq \emptyset$. According
to [ibid.,\,Theorem\,6.1 and Remark\,6.2] we can find a {\em normal} $G$-variety $Z$ and
a {\em finite surjective} $G$-morphism
$p\colon Z\to X$
such that (with $X^s=X^s(L)$, and $Z^s=Z^s(p^*(L))$):
\begin{enumerate}
\item[(i)] $G$ operates {\em freely} on $Z^s$ and the geometric quotient $W=Z^s/G$
exists as a normal variety, and the quotient map
$\beta\colon Z^s\to W$
is a principal bundle, locally trivial in the {\em Zariski topology}. 
\item[(ii)] If $k(X)$ and $k(Z)$ denote the function fields of $X$ and $Z$ respectively, the
extension $k(X)\to k(Z)$ if finite and normal.
\item[(iv)] If $\Gamma$ is the group of $k(X)$-automorphisms of $k(Z)$, the canonical action
of $\Gamma$ on $Z^s$ commutes with that of $G$.
\item[(v)] If $W$ is quasi-projective, $X^s/\negmedspace/G$ exists as a quasi-projective variety.
(See [ibid.,\,p.\,543,\,Remark\,6.1].)
\end{enumerate}

\section{\bf{Big line bundles}}\label{s:big}

In this section we will show that there is a normal projective variety $Q$ mapping to 
the stack $[X^{ss}/G]$  which is generically finite and dominant
over $Y=(X^{ss}/\negmedspace/G)_{\text{top}}$ and
such that the pull-back of the line bundle $L_{/G}$ to $Q$ is big. Here
$L_{/G}$ is the line bundle on $[X^{ss}/G]$  as in Definition\,\ref{def:L/G}. 

\subsection{Basic lemmas on bigness and nefness} The following two lemmas are the basic
tools for proving bigness and nefness of bundles. 

\begin{lem}\label{lem:nef} Suppose $L$ and $M$ are line bundles on an algebraic
scheme $W$ such that $N=N_\epsilon= L+\epsilon M$
is ample for sufficiently small positive. Then $L$ is nef.
\end{lem}

\proof Let $C\hookrightarrow W$ be a closed irreducible and reduced curve. Since $N_\epsilon$
is ample (for sufficiently small $\epsilon >0$) $\deg N_\epsilon\vert_C>0$. 
This means that for $\epsilon=\frac{b}{a}$ small,
$a\deg{L\vert_C} + b\deg M\vert_C > 0$, i.e.,
\[\deg L\vert_C > -\epsilon  \deg M\vert_C.\]
Letting $\epsilon$ approach zero, we conclude that $\deg L\vert_C \ge 0$.\qed

Here is a criterion for bigness of a line bundle  on
a variety in terms of the nefness of associated bundles on a blow-up of the ambient variety.

\begin{lem}\label{prop:big} Let $\psi\colon W'\to W$ be the blow up of an irreducible projective
variety $W$ by a coherent
ideal sheaf $I$ of $\co_W$ whose support is a finite number of points on $W$, and suppose
$r\set \dim{W} \ge 2$. Let $J=I\co_{W'}$ be the ideal sheaf of the exceptional divisor of the
blow-up and $L$ a nef line bundle on $W$ such that $\psi^*(L)+\epsilon J$ is nef on
$W'$ for sufficiently small $\epsilon$. Then $L^{(r)}>0$, i.e., $L$ is big on $W$.
\end{lem}

\proof Since $L$ is nef on $W$, $L^{(r)}\ge 0$. It suffices to show that $L^{(r)}=0$ leads to
a contradiction. The restriction of $L$ to the finite number of points which are blown up
 is trivial, whence using
the ``projection formula" (see \cite[p.\,296, Proposition\,2.11]{K2}) we conclude that
\[ \psi^*(L)^k\cdot J^{(l)} = 0 \qquad (k+l=r, \, k\ge 1).\]
It follows that
 \stepcounter{thm}
 \begin{equation*}\label{eq:psi-L-J}\tag{\thethm}
 (\psi^*(L)+\epsilon{J})^{(r)} = \psi^*(L)^{(r)}+\epsilon^rJ^{(r)}.
 \end{equation*}
Since $L^{(r)}=0$, one sees that $\psi^*(L)^{(r)}=0$. The LHS of \eqref{eq:psi-L-J} is
non-negative since $\psi^*(L)+\epsilon{J}$ is nef. Hence to get a contradiction it suffices to
show that $J^{(r)} < 0$. Let $E$ be the exceptional divisor of the blow-up $\psi$. One knows that
$J\vert_E$ is ample on $E$. By the Asymptotic Riemann-Roch (see 
\cite[p.\,208,\,Corollary\,2.14]{K2}) we know that if $M$ is a line bundle on a projective
variety $X$ and $F$ a coherent sheaf on $X$ with $\dim{\supp{F}}=l$, then
\stepcounter{thm}
 \begin{equation*}\label{eq:asRR}\tag{\thethm}
\chi(X,\,M^n\otimes F) = \dfrac{(M^{(l)}\cdot F)}{l!}n^l + O(n^{l-1}).
 \end{equation*}
One needs to set $X=\supp{F}$ in {\em loc.cit.~}to get the above formula, and this can
clearly be done without loss of generality. The LHS of \eqref{eq:asRR} is a polynomial function
in $n$ and $l$. Next consider the exact sequence
\[0\to J \to \co_{W'} \to \co_E \to 0.\]
We have (for $n\gg 0$) 
\[\chi(W',\,J^{n+1}) = \dfrac{J^{(r)}}{r!}(n+1)^r + O((n+1)^{r-1}), \]
and 
\[\chi(W',\,J^{n}) = \dfrac{J^{(r)}}{r!}n^r + O(n^{r-1}). \]
Taking the difference we get
\stepcounter{thm}
\begin{align*}\label{eq:asym-diff}
\chi(W',\,J^{n+1})-\chi(W',\,J^{n}) & = \dfrac{J^{(r)}}{r!}[(n+1)^r-n^r]+\dotsb+O(n^{r-2}) \\ \tag{\thethm}
& = \dfrac{J^{(r)}}{r!}r\cdot n^{(r-1)} + O(n^{r-2}) \\ 
& = \dfrac{J^{(r)}}{(r-1)!}\cdot n^{(r-1)} + O(n^{r-2}).
\end{align*}

On the other hand, we have the exact sequence 
\[ 0 \to J^{n+1} \to J^n \to J^n\vert_E \to 0\]
whence
\[-\chi(W',\,J^{n+1})+\chi(W',\,J^n) = \chi(E,\,J^n\vert_E).\]
Since $J\vert_E$ is ample,
\[\chi(E,\,J^n\vert_E) = \dfrac{b}{(r-1)!}n^{r-1} + O(n^{r-2})\]
with $b>0$. Comparing this with the asymptotic formula in \eqref{eq:asym-diff} above
we get $b=-J^{(r)}$. In other words $J^{(r)}<0$, giving the sought for contradiction.
\qed

\subsection{Equivariant Blow-ups}\label{ss:bl-up}

In this sub-section we fix an {\em irreducible standard quotient data} (see \Dref{def:usual})
\stepcounter{thm}
\begin{equation*}\label{eq:strng-irr}\tag{\thethm}
\De=(X,\,L,\,\bp(V),\,X^{ss}\xrightarrow{\alpha}Y)
\end{equation*}

For the rest of this section we deal with the following situation:
Let $u_0$ be a point in $Y$, and consider the reduced closed subscheme $C$
of $X$ which is the closure of $\alpha^{-1}(u_0)$ in $X$. 
Then $C$ is a closed $G$-invariant subscheme of $X$. Consider the blow-up
\[\theta\colon X'\to X\]
of $X$ along $C$. If $E\hookrightarrow X'$ is the exceptional divisor of the blow-up,
and $I(E)$ the resulting ideal sheaf of $E$ in $\co_{X'}$, then $E$ is $G$-invariant, and
$I(E)$ is a $G$-invariant invertible $\co_{X'}$-module which is relatively ample for the
map $\theta\colon X'\to X$. We then have

\begin{lem}\label{lem:bl-up1} There exists $\epsilon_0>0$ such that:
\begin{enumerate}
\item $L'_\epsilon=\theta^*(L)+\epsilon I(E)$ is ample on 
$X'$ for $0<\epsilon\le\epsilon_0$. 
\item $(X')^{ss}(L'_{\epsilon_{{}_1}})=X^{ss}(L'_{\epsilon_{{}_2}})$ for 
$0< \epsilon_1, \epsilon_2 \le \epsilon_0$. Let the common variety be denoted $X^{ss}$.
\item $(X')^{ss}\subset \theta^{-1}(X^{ss})$. Equivalently, the inverse image under $\theta$ of the
unstable locus of $(X,L)$ is contained in the unstable locus of $(X',\,N_\epsilon)$.
\item $(X')^{ss}\setminus E = X^{ss}\setminus C$. In particular $(X')^{ss}\neq \emptyset$.
\end{enumerate}
\end{lem}

\proof Since $I(E)$ is relatively ample, the assertion about the ampleness of $L'_\epsilon$ is clear.
Parts (1) and (2) follows from \Pref{prop:b/a-small2}. Part (3) is immediate from 
\cite[p.\,352,\,Theorem\,2.3(a)]{R}.\qed

It is, at this point, convenient for us to extend the notion of {\em stabilizing} as
well as the notion of {\em $\epsilon_0$-stabilizing} defined in \Dref{def:epsilon}  to include 
schemes $Z$ which are not necessarily
the product of $X$ and $Y$ but are close to that. This is primarily because we wish to work
with Zariski locally trivial principal $G$-bundles, and such fibrations will be obtained
by replacing
$X\times_k(B\backslash G)$ by its normalization and then by eliminating finite isotropies. To
that end we make the following definition (and remind the reader about the conventions
in \ref{conv:L+M}).

\begin{defi}\label{def:epsilon2} Let $Y$ be a projective $G$-scheme as in
\Lref{lem:S-exist2}, i.e., ${\mathrm{Pic}}\,(Y)$
is generated by a finite number of ample line bundles such that the action of $G$ on $Y$ is
linear with respect to all these generators.
Let $Z$ be a scheme which admits $G$-invariant morphisms $p\colon Z\to X$ and 
$\pi\colon Z\to Y$.
We say that a line bundle $M$ on $Y$ is {\em stablizing for $(X,\,L,\,Z)$} if 
the universal
map $(p,\,\pi)\colon Z\to X\times_kY$ is {\em finite and surjective} and 
(a) and (b) in \Dref{def:epsilon} are satisfied for this $Z$ and this $p$. 
Equivalently $M$ is stablizing for $(X,\,L,\,Z)$ if the universal map $(p,\pi)\colon Z\to
X\times_k Y$ is finite surjective and $M$ is stabilizing for $(X,L)$. Note that $M$ is stablizing
for $(X,\,L)$ if and only if it is so for $(X,\,L,\,X\times_kY)$.
We will say $M$ is $\epsilon_0$-stabilizing for $(X,\,L,\,Z)$ if $Z\to X\times_kY$ is finite surjective
and $M$ is  $\epsilon_0$-stabilizing for $(X,\,L)$.
\end{defi}

Let $\epsilon_0$ be as in the conclusion of \Lref{lem:bl-up1}. For the rest of this discussion
we fix a positive rational number $\epsilon$ such that $\epsilon\le \epsilon_0$, and 
as in {\em loc.cit.} use $L'_\epsilon$ for $\theta^*L+\epsilon I(E)$.

\begin{lem}\label{lem:bl-up2} There exists a positive rational number $\eta_0$, a line bundle
$M$ on $B\backslash G$, a normal projective $G$-variety $Z$ together with $G$-invariant
maps $p\colon Z\to X$, $\pi\colon Z\to B\backslash G$ such that
\begin{enumerate}
\item[(i)] $M$ is $\eta_0$-stablizing for $(X,\,L,\,Z)$ and $(X',\,L'_\epsilon)$. Fix $\eta>0$ with
$\eta\le \eta_0$ and write $N=N_\eta$ for the line-bundle $L+\eta M\set p^*(L)+\pi^*(M)$ on $Z$.
\item[(ii)] $Z^s\neq\emptyset$, where we write  $Z^s=Z^s(N)$.
\item[(iii)] $G$ acts freely on $Z^s$ and we have a (Zariski locally trivial) principal $G$-bundle 
\[\beta\colon Z^s\to W=W_M\]
where $W$ is a normal projective variety on to which the line bundle $N\vert_{Z^s}$ descends
as a line ample line bundle $N_W$.
\item[(iv)] A power of the line bundle $p^*(L)\vert_{Z^s}$ descends to a line bundle $L_W$ on
$W$ and $L_W$ is nef on $W$.
\end{enumerate}
\end{lem}

\proof Applying \Pref{prop:b/a-small2} and \Pref{prop:nonempty} to $(X_1,\,X_2;\,L_1,\,L_2;\,Y)$
with $X_1=X$, $X_2=X'$, $L_1=L$, $Y=L'_\epsilon$, and $Y=B\backslash G$, we deduce
the existence of $\eta_0>0$ and a line bundle $M$ on $B\backslash G$ such that $M$ is
$\eta_0$-stabilizing for $(X,\,L)$ and $(X',\,L'_\epsilon)$ and such that
$[X\times_k (B\backslash G)]^s(L +\eta M)\neq \emptyset$ for $0<\eta<\eta_0$
(\Pref{prop:nonempty} is needed for the last assertion).

Recall that if $G$ acts on an algebraic scheme $S$, then it acts on $S_{\mathrm{red}}$ as well
as the nornalization $\tilde{S}$ of $S_{\mathrm{red}}$ in a canonical way so that
the maps $S_{\mathrm{red}}\to S$ and $\tilde{S}\to S_{\mathrm{red}}$ are $G$-invariant.

The scheme $Z$ is obtained in two steps. First normalize 
$[X\times_k(B\backslash G)]_{\mathrm{red}}$ and then apply the technique of elimination of
finite isotropies described in \S\S\,\ref{ss:elimination}. Assertions (i)---(iii) follow easily
from the results quoted in \S\S\,\ref{ss:elimination}.

The first part of part\,(iv) is straightforward and well-known. Let $L_W$ be such a line bundle.
This means that for a suitable $n\ge 1$, $(p^*L^n)\vert_{Z^s}=\beta^*L_W$. We say
{\em $L$ descends to $\frac{1}{n}L_W$} to describe this situation. The second part is
a direct application of \Lref{lem:nef}. In fact since $N\vert_{Z^s}= \beta^*N_W$ and 
$(p^*L^n)\vert_{Z^s}=\beta^*L_W$ for a suitable $n\ge 1$, therefore a positive power of
$(\pi^*M)\vert_{Z^s}$ descends to a line bundle $M_W$ on $W$. A little thought shows that
\[N_W= L_W+ a\eta M_W\]
where $a$ is a positive rational number. 
This $a$ is independent of $\eta$ and hence \Lref{lem:nef} applies.
\qed

\bigskip

More can be said. 
Note that $L$ can be replaced by a positive power of itself without affecting the loci
$X^{ss}(L)$, $X^{s}(L)$, $Z^s(N)$ etc. According to \ref{conv:L-power}, there is a postive
integer $n$ such that $L^n$ is trivial on the fibres of $\alpha$. We replace $L$ by
$L^n$ for the rest of this section and assume without loss of generality that $L$ is trivial
on the fibres of $\alpha$. We therefore have a line bundle $L_{/G}$ on the stack $[X^{ss}/G]$
as in \Dref{def:L/G}.
The map $q\colon W\to Y$ in the proof of \Lref{lem:bl-up2} then has the following interpretation.
Since $\beta\colon Z^s\to W$ is a principal bundle and $p\colon Z^s\to X^{ss}$ is
$G$-invariant, by definition, we have a classifying map $f_\beta\colon W\to [X^{ss}/G]$. As in \Sref{s:strata}
we have a map $\gamma\colon [X^{ss}/G]\to Y$ and a continuous map
$q=\gamma\smcirc f_\beta\colon{W}\to Y$ ($q$ is a map of stratified spaces) 
such that the diagram
\stepcounter{thm}
\begin{equation*}\label{diag:beta-gamma}
{\xymatrix{
Z^s \ar[rr]^p \ar[d]_\beta & & X^{ss} \ar[d] \ar@/^/[ddr]^{\alpha} &\\
W \ar@/_/[drrr]_q \ar[rr]^{f_\beta} & & [X^{ss}/G] \ar[dr]^{\gamma} &  \\
& & & Y
}}\tag{\thethm}
\end{equation*}
commutes. According to \Pref{prop:nonempty} $q\colon W\to Y$ is surjective. Note that the line
bundle $L_W$ (to which $p^*L\vert_{Z^s}$ descends) is given by the formula
\[L_W=f_\beta^*L_{/G}.\]

At this point we draw the reader's attention to \Lref{lem:bl-up1}, which deals with the blow up
of $X$ along the closure of $\alpha^{-1}(u_0)$ in $X$, and remind the reader that
we have fixed an $\epsilon$ in the interval $(0,\epsilon_0]$ where $\epsilon_0$ is as in the
conclusion of the Lemma.

\begin{prop}\label{prop:bl-up} Let $W$ be the scheme in part\,{\em{(}}iii{\em{)}} of \Lref{lem:bl-up2}. 
Let $C_W\subset W$ be the reduced closed subscheme of $W$ given by
$C_W=q^{-1}(u_0)$, where $q\colon W\to Y$ is the map in 
Diagram\,\eqref{diag:beta-gamma}. 
If $W'' \xrightarrow{\psi} W$ is the 
blow up of $W$ along $C_W$ and $J$ the {\em{(}}coherent invertible{\em{)}} $\co_{W''}$-ideal of 
the exceptional divisor $E''$ of $W''\xrightarrow{\psi} W$,
then there exists a positive rational number $a$, independent of $\epsilon<\epsilon_0$,
such that $\psi^*(L_W)+a\epsilon J$ is nef on $W''$.
\end{prop}

\proof Recall that the line bundle $M$ on $B\backslash G$ is $\eta_0$-stabilizing for
$(X',\,L'_\epsilon)$ (see \Lref{lem:bl-up1}\,(i)). Set
\[ Z'\set Z\times_kX'.\]
Then we have a finite surjective $G$-invariant map $Z'\to X'\times_kB\backslash G$ and
hence $M$ is $\eta_0$ stablizing for $(X',\,L'_\epsilon,\,Z')$. We have a commutative diagram
\stepcounter{thm}
\begin{equation*}\label{diag:Z'Z} \tag{\thethm}
{\xymatrix{
Z' \ar[d]_{\theta'} \ar[r]^{p'} & X' \ar[d]^{\theta}\\
Z \ar[r]_{p} & X
}}
\end{equation*}
Let $\pi'\colon Z'\to B\backslash G$ be the map $\pi\smcirc\theta'$.
Fix $\eta\in (0,\,\eta_0]$ and write
\[ N'(=N'_\eta = N'_{\eta,\epsilon})\set {p'}^*L'_\epsilon +\eta{\pi'}^*M = L'_\epsilon+\eta M.\]

We assume $Z'$ is {\em normal and $G$ acts freely on} $(Z')^s(N')$ by first replacing
$Z'$ by its normalization and then eliminating finite isotropies (see \S\S\,\ref{ss:elimination}).
$M$ remains $\eta_0$-stabilizing for $(X',\,L'_\epsilon,\,Z')$ through all these modifications.
We have a finite surjective map 
\[\varphi\colon Z'\to Z\times_k X' \qquad (\varphi=(\theta',\,p')).\]
Set
\begin{align*}
(Z')^s & = (Z')^s(N') \\
(Z^s)' & = \varphi^{-1}(Z^s\times_{X^{ss}}(X')^{ss}).
\end{align*}
(Recall $(X')^{ss}\set (X')^{ss}(L'_\epsilon)$. See \Lref{lem:bl-up1}.) Note the difference
between $(Z^s)'$ and $(Z')^s$---they need not be the same, though both are $G$-invariant
open subvarieties of $Z'$. After all these modifications, we still have the commutative diagram
\eqref{diag:Z'Z}. The map $\theta'\colon Z'\to Z$ is more or less that blow up of $Z$ along
$p^{-1}(C)$ where, recall, $C$ is the centre of the blow-up $\theta\colon X'\to X$. 

Since $M$ is $\eta_0$-stablizing for $(X',\,L'_\epsilon,\,Z')$, we have $(Z')^s(N')=(Z')^{ss}(N')$.
Further, $Z'$ is normal and $G$ operates freely on $(Z')^s$. We therefore have (Zariski locally
trivial) principal $G$-bundle
\[\beta'\colon (Z')^s\to W'\]
with $W'$ normal projective, and such that $N'$ descends to an ample line bundle
$N'_{W'}$ on $W'$. Note that $N'_{W'}$ depends on $\epsilon$ and $\eta$, and if we
wish to draw attention to this, we write $N'_{W',\eta}$, or $N'_{W',\eta,\epsilon}$. Since
$L$ is trivial on the fibres of $\alpha\colon X^{ss}\to Y$, the line bundle ${p'}^*\theta^*L$
descends to a line bundle $L_{W'}$ on $W'$. In fact
\[L_{W'} = f_{\beta'}^*L_{/G}\]
where $f_{\beta'}\colon W'\to [X^{ss}/G]$ is the classifying map for the data
\[((Z')^s\xrightarrow{\beta'} W', (Z')^s\to X^{ss})\] 
consisting of a principal bundle and a $G$-invariant map, the second map being induced by
the composite $\theta\smcirc p'$.

The exceptional divisor $E$ of $\theta\colon X'\to X$ pulls back to locally principal (i.e., Cartier)
effective divisor $E'={p'}^{-1}(E)$ on $Z'$. This divisor $E'$ can be related to the
centre $C_W$ of the blow-up $\psi\colon W''\to W$ in the following way: We have a composite
$(Z^s)'\to Z^s\xrightarrow{\beta} W$. The inverse image of $C_W$ in $(Z^s)'$ under this
composite is then $E'\vert_{(Z^s)'}$. Since the latter is an effective Cartier divisor, by the universal
property of blow-ups, we have a unique map 
\[\beta''\colon (Z^s)'\to W''\]
such that $(\beta'')^{-1}(E'')=E'\vert_{(Z^s)'}$.

The spaces $(Z^s)'$ and $(Z')^s$ are open subsets of $Z'$ and their intersection $\tilde{Z}$
in $Z'$ is a $G$-invariant variety. The geometric quotient 
$\widetilde{W}=\widetilde{Z}/G$ exists as
a normal quasi-projective variety, being an open subvariety of  $W'=(Z')^s/G$. Since the
composite $\widetilde{Z}\to (Z^s)'\xrightarrow{\beta''} W''$ is $G$-invariant ($G$ acting
trivially on $W''$) and $\widetilde{W}$ is a categorical quotient we have a
(dominant) map $\widetilde{W}\to W''$ such that the composite $\widetilde{Z}\to (Z^s)'\to W''$
agrees with the composite $\widetilde{Z}\to\widetilde{W}\to W''$.

Let ${\overline{W}}$  be the scheme theoretic closure of the locally closed embedding
\[\widetilde{W} \xrightarrow{\text{diag}} W''\times_k W'.\]

Let $\mu\colon {\overline{W}}\to W''$, $\lambda\colon {\overline{W}}\to W'$, be the projections.

The situation is best described by the  diagram below where 
$Y'$ is the topological space consisting of semi-stable
equivalence classes on $(X')^{ss}$ and $\alpha'\colon (X')^{ss}\to Y'$ the
resulting map. The other arrows, not defined earlier, are as described after the diagram:
\[
\xymatrix{
& Z^s  \ar[ld]_p \ar'[d][ddd] _\beta & & & (Z^s)' \ar[lll]_{\theta'} \ar[ld]_{p'}
\ar'[d][ddd]_{\beta''} & & &  \\
X^{ss} \ar[ddd]_{\alpha} & & & (X')^{ss} \ar[ddd]^{\alpha'} \ar[lll]_{\theta} 
& & &  (Z')^s\ar[ddd]^{\beta'} \ar[lll]_{p'} & \\
&&&&&&& \\
& W \ar[ld]_q & & & W'' \ar'[l][lll]_{\psi} & & & {\overline{W}}  
\ar[ld]_\lambda \ar'[l][lll]_{\mu} \\
Y &&& Y' \ar[lll]_{\phi} & & & W' \ar[lll]_{q'} &
}
\]

\medskip

The map $q'\colon W'\to Y'$ is the obvious continuous map (the exact analogue
of $q\colon W\to Y$), $\phi\colon Y'\to Y$ the map induced by the equivariant
$G$-map $\theta$. 

The bottom rectangle in the diagram commutes for the following reason: First, we have
a commutative diagram
\[
{\xymatrix{
Z^s \ar[d]_{p} & (Z^s)' \ar[l]_{\theta'} \ar[d]^{p'} & {\widetilde{Z}} \ar@{_{(}-}[l] \ar@{_{(}-}[d] \\
X^{ss} & (X')^{ss}\ar[l]^{\theta} & (Z')^s \ar[l]^{p'}
}}
\]
from which it follows that
\[q\smcirc \psi\smcirc\mu\vert_{\widetilde{W}}=\phi\smcirc q'\smcirc \lambda\vert_{\widetilde{W}}.\]
with
${\widetilde{W}}\subset {\overline{W}}$ as above 
(i.e. ${\widetilde{W}}= (Z^s)'\cap(Z')^s/\negmedspace/G$). 
Since the diagonal map $Y\to Y\times Y$ is a closed immersion (see \Pref{prop:G-sep}), 
the locus in ${\overline{W}}$ where
the two maps $q\circ\psi\circ\mu$ and $\phi\circ q'\circ \lambda$ are equal is a closed subset
of ${\overline{W}}$. It follows, since ${\widetilde{W}}$ is dense in ${\overline{W}}$, that 
the bottom rectangle of the 3-dimensional diagram above commutes. 

The remaining sub-rectangles of the diagram clearly commute.

Now using $P_1\set {\overline{W}}\times_WZ^s$ and 
$P_2\set {\overline{W}}\times_{W'}(Z')^s$, we have
principal $G$-bundles $P_1\to {\overline{W}}$ and $P_2\to {\overline{W}}$ and $P_1$ and 
$P_2$ have $G$-invariant maps to $X^{ss}$, whence we have classifying maps 
\[{\xymatrix{
{\overline{W}} \ar @{->}^{f_1{\phantom{XX}}} @< 2pt> [r]
          \ar @{->}_{f_2{\phantom{XX}}} @<-2pt> [r] & [X^{ss}/G].
}}
\]
If $\gamma\colon [X^{ss}/G]\to Y$ is the map in \eqref{eq:gamma2} then commutativity of
the bottom rectangle of the 3-dimensional diagram amounts to saying that
\[\gamma\smcirc f_1= \gamma\smcirc f_2.\]
By \Pref{prop:Pindep} it then follows that $f^*_1L_{/G} = f_2^*L_{/G}$.
This is the same as saying:
\[\mu^*(\psi^*(L_W))\cong \lambda^*L_{W'}.\]

Recall that the ideal sheaf of the exceptional divisor $E$ of $\theta\colon X'\to X$ was denoted $I$.
Let $I'={p'}^*I\vert_{(Z')^s}$. A positive power, say $(I')^l$ descends to an invertible ideal
sheaf $I_{W'}$ of $\co_{W'}$. We express this by saying $I'$ descends to $\frac{1}{l}I_{W'}$.
Similarly the pull back of $M$ to $(Z')^s$ descends to, say $\frac{1}{m}M_{W'}$. 
We point out that $L_W'+\frac{\epsilon}{l}I_{W'}$ is nef since
 $N'_{W',\eta,\epsilon}=L_{W'}+\frac{\epsilon}{l}I_{W'}+
\frac{\eta}{m}M_{W'}$ is ample, whence \Lref{lem:nef} applies. This in turn means that
$\lambda^*(L_{W'} + \frac{\epsilon}{l}I_{W'})$ is nef.

We have, therefore, two coherent ideal sheaves, $\mu^*J$ and $\lambda^*I_{W'}$,  of 
$\co_{\overline{W}}$. Moreover, the closed topological subspaces of ${\overline{W}}$ underlying
the closed subschemes defines by these two ideals are the same. Since $\mu^*J$ and
$\lambda^*I_{W'}$ are invertible ideal sheaves, whence locally principal, and ${\overline{W}}$
is normal, we have positive integers $r$ and $s$ such that $\mu^*J^r=\lambda^*I_{W'}^{s}$.
Choosing $a=\frac{r}{sl}$ we see that 
\[\mu^*(\psi^*L_W + a\epsilon J) = \lambda^*(L_{W'} + \frac{\epsilon}{l}I_{W'}).\]
We have argued that the right side is nef. Hence so is $\psi^*L_W+ a\epsilon J$. 
\qed

\begin{cor}\label{cor:curve} If $C$ is a closed integral curve in $W$ such that $q\vert_C$ is not
a constant, then $\deg(L_W\vert_C)>0$.
\end{cor}

\proof Without loss of generality we can assume $\dim{W}>1$ 
by embedding $W$ in a higher dimensional variety (by enlarging $V$ for example, i.e., by 
considering the embedding of $X$ via the complete linear system
$\Gamma(X,\,L^n)$ for $n\gg 0$). Thus $C$ is
finite over $Y$ and $\dim{W}  > 1$. Pick $u_0\in Y$ and consider the blow-up $\psi\colon
W''\to W$ of the Proposition. Let $C'\hra W''$ be the proper transform of $C$. By the Proposition 
we have $\deg(\psi^*(L_W+\epsilon J)\vert_{C'})\ge 0$. Since $\deg(J\vert_{C'})<0$, this implies
that $\deg(\psi^*L_W\vert_{C'}) >0$, whence $\deg(L_W\vert_C) > 0$.
\qed

\section{\bf{Geometric Reductivity}}\label{s:conclusion}

Let $\De=(X,\,L,\,\bp(V),\,X^{ss}\xrightarrow{\alpha} Y)$ be an irreducible standard quotient
data.

All that remains to complete the proof of \Tref{thm:main} (i.e., to show $G$ is geometrically 
reductive) is to find a map $f\colon Q\to [X^{ss}/G]$ as in \eqref{map:classifying} satisfying
conditions (1)---(7) in \S\S\ref{ss:zariski-triv}. To that end pick a generic quotient $U$ for
the data $\De$ and reconsider Diagram\,\eqref{diag:beta-gamma} above, namely:
\[
{\xymatrix{
Z^s \ar[rr]^p \ar[d]_\beta & & X^{ss} \ar[d] \ar@/^/[ddr]^{\alpha} &\\
W \ar@/_/[drrr]_q \ar[rr]^{f_\beta} & & [X^{ss}/G] \ar[dr]^{\gamma} &  \\
& & & Y
}}
\]
The map $q^{-1}(U)\to U$ is a morphism of schemes, for the pull-back
by $p$ of a $G$-invariant regular function on $X^{ss}$ descends to a regular
function on $q^{-1}(U)$. The map $q^{-1}(U)\to U$ is proper since $W$ is projective and 
by \Cref{cor:proper} $X^{ss}$ mod $G$ is proper. We can certainly find a closed subvariety
of $q^{-1}(U)$ which is generically finite over $U$. By normalizing the closure
of this variety in $W$ and then applying
\Lref{lem:Gam-action} we get a a normal projective variety $Q$ and a map $f\colon Q\to
[X^{ss}/G]$ such that resulting map $Q\to Y$ is generically finite, and satisfies 
properties (1), (2) and (5) required of the map \eqref{map:classifying}. We continue to
use the symbol $q$ for the map $Q\to Y$.

Write $L_Q$ for $L_W\vert_Q=f^*L_{/G}$.  By \Lref{lem:bl-up2}\,(iv) we see that
$L_Q$ is nef on $Q$. Suppose $\dim Q >1$.  Pick a point $u_0\in U$ such that 
$q^{-1}(u_0) \subset Q$ is finite. Blow up $Q$ along this inverse image. According to
\Pref{prop:bl-up}, if $J$ is the ideal of the exceptional divisor of the blow-up $g\colon Q'\to Q$
along this inverse image, then $g^*L_Q+ \epsilon J$ is nef for sufficently small positive values
of $\epsilon$, whence by  \Pref{prop:big} $L_Q$ is nef and big on $Q$. If $\dim{Q}=1$, then 
according to \Cref{cor:curve} $\deg{L_Q}>0$, whence $L_Q$ is nef and big (in fact ample). So
in every case $L_Q$ is nef and big on $Q$. This establishes (3) of \S\S\,\ref{ss:zariski-triv}.
Property (6) follows  from the manner in which $Z$ and $Z^s$ were
found. As for property (4), this follows directly from \Cref{cor:curve}.

It remains to prove (7) of \S\S\,\ref{ss:zariski-triv} for our map $f\colon Q\to [X^{ss}/G]$. By replacing
$L$ by a positive power of itself if necessary, we may assume $L_Q$ descends to a line
bundle $L_{\ol{Q}}$ on $\ol{Q}\set Q/\Gamma$. We have to show that if 
$t\in \Gamma(\ol{Q},\,L_{\ol{Q}})$, then for some $n\gg 0$, with $n=p^m$, $t^n$ comes from a
$G$-invariant section on the pull-back of $L^n$ on the normalisation $\wt{X}^{ss}$ of 
$X^{ss}$. The proof of this runs along the same lines as that of \eqref{eq:LQ} in the proof of 
\Pref{prop:Pindep}. Let $g\colon {\wt{X}}^{ss}\to X^{ss}$ be the normalisation map and
let $\wt{L}=g^*L$.
Let $s$ be the $G$-invariant meromorphic section of $\wt{L}$ induced by $t$. (We point out
that $t$ can be regarded as a $(G,\,\Gamma)$-invariant section of 
${\wt{\pi}}_*\wt{\pi}^*\wt{L}$.) Let $D$ be the polar divisor of $s$ and let $x_1\in D$ be such
that $s(x_1)=\infty$. Let $y=\alpha g(x_1)$. Then as in \eqref{eq:subset} in the proof of
\Pref{prop:Pindep} we see that for any $x\in (\alpha\smcirc g)^{-1}(y)$, $s(x)=\infty$. 
Now on the open set
$(\alpha\smcirc \pi)^{-1}(U) = (q\smcirc \beta)^{-1}(U)$, the section $\wt{\pi}^*(s)$ coincides with
the section $\beta^*\varphi^*(t)$, where $\varphi\colon Q\to \ol{Q}$ is the natural map. 
Further since
$q\colon W \to Y$ is surjective $\wt{\pi}(P)$ meets $(\alpha\smcirc g)^{-1}(y)$. Thus we have
$x_2\in \wt{\pi}(P)$, with $\alpha(g(x_2))=y$. Since $\beta^*\varphi^*(t)$ is regular, we see that
$s(x_2)$ cannot be $\infty$. This leads to a contradiction and we done.

\end{document}